\documentclass[11pt]{article}
\usepackage[T1]{fontenc}
\usepackage{lmodern,microtype}
\usepackage[hmargin=24mm,vmargin=20mm]{geometry}
\usepackage{amsmath,amssymb,amsthm}
\usepackage{aliascnt}
\usepackage{booktabs,tabularx,array,enumitem}
\usepackage[dvipsnames]{xcolor}
\usepackage{tikz}
\usetikzlibrary{arrows.meta}

\usepackage[colorlinks=true,linkcolor=blue!50!black,citecolor=green!40!black,urlcolor=blue!50!black]{hyperref}
\usepackage[nameinlink,capitalise,noabbrev]{cleveref}
\newtheorem{theorem}{Theorem}[section]
\newaliascnt{lemma}{theorem}
\newtheorem{lemma}[lemma]{Lemma}
\aliascntresetthe{lemma}
\newaliascnt{proposition}{theorem}
\newtheorem{proposition}[proposition]{Proposition}
\aliascntresetthe{proposition}
\newaliascnt{corollary}{theorem}

\aliascntresetthe{corollary}
\DeclareMathOperator*{\argmin}{arg\,min}
\newcommand{\R}{\mathbb R}
\newcommand{\ip}[2]{\langle #1,#2\rangle}
\newcommand{\norm}[1]{\lVert #1\rVert}
\numberwithin{equation}{section}
\newcommand{\tocblock}{\addtocontents{toc}{\protect\vskip 0pt plus 7\protect\baselineskip
 \protect\penalty-200 \protect\vskip 0pt plus -7\protect\baselineskip\relax}}
\let\plainsection\section
\renewcommand{\section}{\tocblock\plainsection}
\setlist[enumerate]{label=\textup{(\alph*)},itemsep=3pt,topsep=4pt}
\allowdisplaybreaks[2]
\title{Silver rate is sharp for gradient descent}
\author{%
 Yechan Jeong\\
 {\normalsize Department of Mathematical Sciences}\\
 {\normalsize Seoul National University}\\
 {\normalsize\href{mailto:yechan@snu.ac.kr}{\texttt{yechan@snu.ac.kr}}}
 \and
 Youngjoon Hong\\
 {\normalsize Department of Mathematical Sciences}\\
 {\normalsize Seoul National University}\\
 {\normalsize\href{mailto:hongyj@snu.ac.kr}{\texttt{hongyj@snu.ac.kr}}}}
\date{}
\hypersetup{pdftitle={Silver rate is sharp for gradient descent},pdfauthor={Yechan Jeong, Youngjoon Hong}}
\begin{document}
\maketitle
\begin{abstract}
We prove sharp lower bounds for gradient descent with predetermined
real stepsizes on smooth convex and smooth strongly convex functions.
Let $p=\log_2(1+\sqrt{2})\approx1.2716$. For every $L>0$ and every
schedule $h\in\mathbb{R}^n$ of length $n\ge1$, there exist a convex
$L$-smooth function $f$ with a minimizer $x_\star$ and an initial point
$x_0\ne x_\star$ for which the last iterate $x_n$ satisfies
$(f(x_n)-f(x_\star))/(L\|x_0-x_\star\|^2)\ge c\,n^{-p}$ with an absolute
constant $c>0$. Fix $L>0$, $\kappa\ge\kappa_0>1$ and $0<\varepsilon<1$.
Any schedule that guarantees
$\|x_n-x_\star\|^2\le\varepsilon\|x_0-x_\star\|^2$ for every initial
point and every $L$-smooth and $L/\kappa$-strongly convex function
requires at least $c(\kappa_0)\,\kappa^{1/p}\log(1/\varepsilon)$ steps,
where $c(\kappa_0)>0$ depends only on $\kappa_0$. The two lower bounds
match the rates of the silver stepsize schedules up to constant factors
in the convex case and for $0<\varepsilon\le1/2$ in the strongly convex
case. We thus prove the conjecture of Altschuler and Parrilo that the
silver rates are optimal. We construct hard instances backward from the
last iterate and derive a single lower bound on the ratio of final
distance to initial distance across four consecutive intervals between
large stepsizes. We then average over stepsize thresholds using a
quantitative stability estimate for an inequality with equality at the
silver ratio $1+\sqrt{2}$.
\end{abstract}
\clearpage
{\small\setlength{\parskip}{0pt}\tableofcontents}
\clearpage
\section{Introduction}\label{sec:intro}

Gradient descent is a basic first-order method. For a convex function
$f:\R^d\to\R$ with $L$-Lipschitz gradient and a minimizer $x_\star$,
the constant stepsize $1/L$ guarantees that $f(x_n)-f(x_\star)$ is at
most an absolute constant times $L\norm{x_0-x_\star}^2/n$ after $n$
steps.
If $f$ is also $L/\kappa$-strongly convex, the stepsize $1/L$ gives
$\norm{x_n-x_\star}^2\le\varepsilon\norm{x_0-x_\star}^2$ after
$O(\kappa\log(1/\varepsilon))$ steps for $0<\varepsilon\le1/2$
\cite{Nes18}. Faster methods have long been known but modify the
update. Nesterov's accelerated method incorporates momentum. With
momentum $f(x_n)-f(x_\star)$ is at most an absolute constant times
$L\norm{x_0-x_\star}^2/n^2$ and
$\norm{x_n-x_\star}^2\le\varepsilon\norm{x_0-x_\star}^2$ holds after
$O(\sqrt\kappa\log(1/\varepsilon))$ steps for $0<\varepsilon\le1/2$
\cite{Nes83,Nes18}. The order $n^{-2}$ matches the classical worst-case
lower bound for first-order methods in unrestricted dimension
\cite{NY83}. For sufficiently large $\kappa$ and $0<\varepsilon\le1/2$
the order $\sqrt\kappa\log(1/\varepsilon)$ matches the classical lower
bound for smooth strongly convex functions \cite{NY83}. Recent work
shows that varying the stepsizes gives part of the speedup of momentum
without changing the update. How fast can gradient descent converge through the choice
of stepsizes alone?

Throughout, \emph{gradient descent} is the iteration
\[
 x_i=x_{i-1}-h_i\nabla f(x_{i-1}),\qquad 1\le i\le n,
\]
with a predetermined stepsize schedule $h=(h_1,\ldots,h_n)\in\R^n$. The
stepsizes may depend on $n$, $L$ and $\kappa$ but not on $f$, on $x_0$
or on the gradients observed during the iteration. Negative stepsizes
are allowed.

For quadratic functions it is classical that stepsizes alone attain
the order $n^{-2}$ in the convex case and the order
$\sqrt\kappa\log(1/\varepsilon)$ for sufficiently large $\kappa$ and
$0<\varepsilon\le1/2$ in the strongly convex case. The stepsizes of
Young \cite{You53} are derived from Chebyshev polynomials. Hence a
lower bound for gradient descent beyond the classical lower bounds
requires functions that are not quadratic. Beyond quadratic functions,
rates faster than $n^{-1}$ were established concurrently by Altschuler
and Parrilo \cite{AP-I,altschuler2025silver} and by Grimmer, Shu
and Wang \cite{GSW23}. The rate of \cite{GSW23} holds for the best
iterate on functions with bounded sublevel set
$\{x:f(x)\le f(x_0)\}$. See \cite{AP26} for an overview.

Altschuler and Parrilo introduced the silver stepsize schedule
in \cite{AP-I}. For smooth convex functions with $x_0\ne x_\star$ the
silver schedule attains the rate
\[
 \frac{f(x_n)-f(x_\star)}{L\norm{x_0-x_\star}^2}=O(n^{-p}),
 \qquad p=\log_2(1+\sqrt2)\approx1.2716,
\]
by \cite{altschuler2025silver}. For $L$-smooth, $L/\kappa$-strongly
convex functions with $\kappa$ sufficiently large the silver schedule
gives $\norm{x_n-x_\star}^2\le\varepsilon\norm{x_0-x_\star}^2$ after
$O(\kappa^{1/p}\log(1/\varepsilon))$ steps for
$0<\varepsilon\le1/2$ \cite{AP-I}. Here $1/p\approx0.7864$. In the
convex case the silver schedule is constructed recursively by joining
two copies of a shorter schedule with a larger intervening step. The
rate $n^{-p}$ follows from progress over several steps even though
individual steps need not decrease $f$. The constructions of
\cite{altschuler2025silver,GSW25} give the rate $n^{-p}$ when $n+1$ is
a power of two. Later constructions attain the rate $n^{-p}$ for
every $n$ \cite{GSW24,ZJ26} and the order
$\kappa^{1/p}\log(1/\varepsilon)$ for every $\kappa>1$ and
$0<\varepsilon\le1/2$ \cite{ZJ26sc}.
Altschuler and Parrilo conjectured that the rates $n^{-p}$ and
$\kappa^{1/p}\log(1/\varepsilon)$ are optimal among all stepsize
schedules \cite{AP-I,altschuler2025silver}.

Ma and Chen \cite{MC26} initiated a sequence of lower bounds for
gradient descent with predetermined nonnegative stepsizes on smooth
convex functions. The most recent results reach the exponent $p$ and
leave a subpolynomial gap to the rate $n^{-p}$. For nonnegative
schedules Ye and Liu \cite{ye2026almost} proved a lower bound of order
$n^{-p-o(1)}$ and Bell \cite{bell2026predetermined} obtained a lower
bound of order $n^{-p}(\log n)^{-0.40998}$. Ye and Liu
\cite{ye2026almost} asked whether lower bounds with the exponent $p$
extend to schedules with negative stepsizes. He and Cheng
\cite{he2026optimal} proved a lower bound of order $n^{-p-o(1)}$ for
real schedules. For strongly convex functions Liu and Ye \cite{LY26}
proved that nonnegative schedules need at least
$\kappa^{1/p-o(1)}\log(1/\varepsilon)$ steps. Liu and Ye asked whether
the subpolynomial factor can be removed, whether the lower bound of
\cite{LY26} extends to schedules with negative stepsizes and whether
the convex and the strongly convex lower bounds can be connected. Thus
the exponent $p$ was known for real schedules and the exponent $1/p$
for nonnegative schedules. Lower bounds that match the silver rates up to a constant
factor remained open.

In this paper we prove the two lower bounds conjectured by Altschuler
and Parrilo. \Cref{thm:main} shows that for every predetermined real
schedule of length $n$, the worst-case value of
$(f(x_n)-f(x_\star))/(L\norm{x_0-x_\star}^2)$ over convex $L$-smooth
functions is at least $cn^{-p}$ with an absolute constant $c>0$.
\Cref{thm:sc-complexity} shows that for $\kappa\ge\kappa_0>1$, every
real schedule that guarantees
$\norm{x_n-x_\star}^2\le\varepsilon\norm{x_0-x_\star}^2$ on $L$-smooth,
$L/\kappa$-strongly convex functions has at least
$c(\kappa_0)\kappa^{1/p}\log(1/\varepsilon)$ steps. The lower bounds of
\cref{thm:main,thm:sc-complexity} match the upper bounds of
\cite{GSW24,ZJ26sc} up to constant factors in the convex case and for
$0<\varepsilon\le1/2$ and $\kappa\ge\kappa_0$ in the strongly convex
case. \Cref{thm:main,thm:sc-complexity} remove the subpolynomial
factors of the lower bounds of
\cite{bell2026predetermined,he2026optimal,LY26,ye2026almost}.
\Cref{thm:main,thm:sc-complexity} hold for every $n$ with no
restriction on the magnitudes, the signs or the order of the stepsizes.
In particular negative stepsizes cannot improve the rates of the
silver schedules by more than a constant factor. \Cref{thm:sc-complexity} answers the questions
of \cite{LY26} on the subpolynomial factor and on negative stepsizes.

\Cref{thm:main,thm:sc-complexity} follow from
\cref{prop:sc-distance}, one estimate for strongly convex functions. We
construct hard instances backward from the last iterate and derive a
single lower bound on the ratio of final distance to initial distance
across four consecutive intervals between large stepsizes. We then
average over stepsize thresholds using \cref{lem:max-plus-excess}, a
quantitative stability estimate for an inequality with equality at the
silver ratio $1+\sqrt2$. \Cref{sec:main-results} states
\cref{thm:main,thm:sc-complexity} and lists the earlier lower bounds in
\cref{tab:lower-bounds,tab:sc-bounds}. \Cref{sec:outline} gives an
overview of the proof and \cref{sec:previous} discusses related work.

\subsection{Main results}\label{sec:main-results}

We state \cref{thm:main,thm:sc-complexity} with the matching upper
bounds and the earlier lower bounds. All norms and inner products are
Euclidean, $\log$ is the natural logarithm and for $x\in\R$ we write
$x^+=\max\{x,0\}$ and $x^-=\max\{-x,0\}$.

Let $L>0$ and let $d\ge1$ be an integer. A differentiable function
$f:\R^d\to\R$ is \emph{$L$-smooth} if
\[
 \norm{\nabla f(x)-\nabla f(y)}\le L\norm{x-y}\qquad(x,y\in\R^d).
\]
For $\lambda>0$ a function $f:\R^d\to\R$ is
\emph{$\lambda$-strongly convex} if $f-\lambda\norm{\cdot}^2/2$ is
convex. For differentiable $\lambda$-strongly convex $f$, the inequality
$f(x)\ge f(0)+\ip{\nabla f(0)}x+\lambda\norm{x}^2/2$ gives
$f(x)\to\infty$ as $\norm x\to\infty$. Hence $f$ attains a minimum.
Strict convexity gives a unique minimizer $x_\star$. We write
$\argmin f$ for the set of minimizers of $f$. For an integer $n\ge1$,
a schedule $h=(h_1,\ldots,h_n)\in\R^n$ and a point $x_0\in\R^d$ define
$x_1,\ldots,x_n\in\R^d$ by
\begin{equation}\label{eq:gd}
 x_i=x_{i-1}-h_i\nabla f(x_{i-1})\qquad(1\le i\le n).
\end{equation}
The points $x_0,\ldots,x_n$ form the \emph{trajectory} of $f$ from
$x_0$ with steps $h$.

\begin{theorem}\label{thm:main}
Let $p=\log_2(1+\sqrt2)$. There is an absolute constant $c>0$ such that
for every $L>0$, every integer $n\ge1$ and every $h\in\R^n$,
\begin{equation}\label{eq:main}
 \sup_{f,x_\star,x_0}
 \frac{2\bigl(f(x_n)-f(x_\star)\bigr)}{L\norm{x_0-x_\star}^2}
 \ge c\,n^{-p},
\end{equation}
where the supremum is over all integers $d\ge1$, all convex $L$-smooth
$f:\R^d\to\R$ with $\argmin f\ne\varnothing$, all $x_\star\in\argmin f$ and all
$x_0\in\R^d\setminus\{x_\star\}$.
\end{theorem}

With the factor $2$ the quotient in \eqref{eq:main} equals
$\norm{x_n}^2/\norm{x_0}^2$ for $f(x)=L\norm{x}^2/2$ and $x_\star=0$.
When $n+1$ is a power of two and $h$ is the silver schedule of
Altschuler and Parrilo~\cite{altschuler2025silver}, the left side of
\eqref{eq:main} with $L=1$ is at most an absolute constant times
$n^{-p}$. By~\cite{GSW24}, for every integer $n\ge1$ there is a
schedule $h\in(0,\infty)^n$ for which the left side of \eqref{eq:main}
with $L=1$ is at most an absolute constant times $(n+1)^{-p}$. By
\cref{lem:scaling} and \eqref{eq:ratio-scaling} the left side of
\eqref{eq:main} does not change when $1$ is replaced by $L>0$ and $h$
by $h/L$. Hence for every $L>0$ and every integer $n\ge1$,
\begin{equation}\label{eq:minimax}
 c\,n^{-p}\le\inf_{h\in\R^n}\,\sup_{f,x_\star,x_0}
 \frac{2\bigl(f(x_n)-f(x_\star)\bigr)}{L\norm{x_0-x_\star}^2}
 \le C\,n^{-p},
\end{equation}
where the supremum is over the $f$, $x_\star$ and $x_0$ of
\cref{thm:main} and $C$ is an absolute constant. The rate $n^{-p}$ is therefore sharp up to an
absolute constant factor. In particular $p$ is the optimal exponent
among all predetermined real stepsize schedules, as conjectured by
Altschuler and Parrilo~\cite{altschuler2025silver}.

The iterates of \eqref{eq:gd} satisfy
$x_i-x_0\in\operatorname{span}\{\nabla f(x_0),\ldots,\nabla f(x_{i-1})\}$.
Hence the lower bound for first-order methods in \cite{NY83} applies.
For $L=1$, all integers $n\ge1$ and all $h\in\R^n$ the left side of
\eqref{eq:main} is at least $c(n+1)^{-2}$ with an absolute constant
$c>0$. \Cref{tab:lower-bounds} lists later lower
bounds. The lower bounds of Ma and Chen \cite{MC26} were sharpened in
\cite{JCY26,Tsai26note,ye2026beyond}. Jung, Cho and Yun \cite{JCY26}
obtained the exponent $\log_2(1+\sqrt3)\approx1.450$ for positive
schedules. For real schedules the lower bounds in
\cite{bell2026predetermined,ye2026beyond} have exponents strictly
larger than $p$. The factor $\exp(-C\sqrt{\log n\,\log\log n})$ in
\cref{tab:lower-bounds} is of the form $n^{-o(1)}$.

\begin{table}[!htbp]
\centering\small
\renewcommand{\arraystretch}{1.28}
\begin{tabularx}{\linewidth}{@{}>{\raggedright\arraybackslash}p{.29\linewidth}
 >{\centering\arraybackslash}p{.18\linewidth}X@{}}
\toprule
Reference & Schedules $h$ & Lower bound\\
\midrule
Ma--Chen \cite{MC26}
 & $[0,\infty)^n$ & $\Omega(n^{-\sqrt{2+\sqrt3}-\delta})$
 for every fixed $0<\delta<2-\sqrt{2+\sqrt3}$\\
Tsai \cite{Tsai26note}
 & $(0,\infty)^n$ & $\Omega(n^{-\sqrt3})$\\
Ye--Liu \cite{ye2026beyond}
 & $\R^n$ & $\Omega(n^{-1.6342})$\\
Jung--Cho--Yun \cite{JCY26}
 & $(0,\infty)^n$ & $\Omega(n^{-\log_2(1+\sqrt3)})$\\
Ye--Liu \cite{ye2026almost}
 & $[0,\infty)^n$ &
 $\Omega\bigl(n^{-p}\exp(-C\sqrt{\log n\,\log\log n})\bigr)$\\
Bell \cite{bell2026predetermined}
 & $[0,\infty)^n$ & $\Omega(n^{-p}(\log n)^{-0.40998})$\\
Bell \cite{bell2026predetermined}
 & $\R^n$ & $\Omega(n^{-1.5570595})$\\
He--Cheng \cite{he2026optimal}
 & $\R^n$ &
 $\Omega\bigl(n^{-p}\exp(-C\sqrt{\log n\,\log\log n})\bigr)$\\
\Cref{thm:main} & $\R^n$ & $\Omega(n^{-p})$\\
\bottomrule
\end{tabularx}
\caption{Lower bounds for
$\sup_{f,x_\star,x_0}2(f(x_n)-f(x_\star))/\norm{x_0-x_\star}^2$ in
\eqref{eq:main} with $L=1$, valid for all sufficiently large $n$ and
all $h$ in the class indicated under `Schedules $h$'. Here $\Omega(\beta_n)$ denotes a
lower bound $c\,\beta_n$ with $c>0$ independent of $n$ and $h$. The
number $C$ is a constant independent of $n$ and $h$.}
\label{tab:lower-bounds}
\end{table}

\begin{theorem}\label{thm:sc-complexity}
Let $p=\log_2(1+\sqrt2)$ and $\kappa_0>1$. There is a number $c(\kappa_0)>0$ such that
\begin{equation}\label{eq:sc-complexity}
 n\ge c(\kappa_0)\,\kappa^{1/p}\log(1/\varepsilon)
\end{equation}
for every $\kappa\ge\kappa_0$, every $L>0$, every $0<\varepsilon<1$,
every integer $n\ge1$ and every $h\in\R^n$ that guarantees
\[
 \norm{x_n-x_\star}^2\le\varepsilon\norm{x_0-x_\star}^2
\]
for every $d\ge1$, every $L$-smooth, $L/\kappa$-strongly convex
$f:\R^d\to\R$ and every $x_0\in\R^d$.
\end{theorem}

Altschuler and Parrilo~\cite{AP-I} proved that the silver schedule gives
$\norm{x_n-x_\star}^2\le\varepsilon\norm{x_0-x_\star}^2$ after
$O(\kappa^{1/p}\log(1/\varepsilon))$ steps for all sufficiently large
$\kappa$ and $0<\varepsilon\le1/2$. Zhang and Jiang~\cite{ZJ26sc}
construct nonnegative schedules of every length for every $\kappa>1$.
For $0<\varepsilon\le1/(2\kappa)$, \cite{ZJ26sc} gives a schedule with
at most
$C\kappa^{1/p}\log(1/(\kappa\varepsilon))$ steps and
\[
 f(x_n)-f(x_\star)\le\frac{\varepsilon L}{2}\norm{x_0-x_\star}^2
\]
for every $L$-smooth, $L/\kappa$-strongly convex $f$ and every $x_0$,
where $C$ is an absolute constant. Strong convexity gives
$\norm{x_n-x_\star}^2\le2\kappa(f(x_n)-f(x_\star))/L$. Let
$\kappa>1$ and $0<\varepsilon\le1/2$. The schedule of \cite{ZJ26sc} for
$\varepsilon/\kappa$ has at most $C\kappa^{1/p}\log(1/\varepsilon)$
steps and gives
$\norm{x_n-x_\star}^2\le\varepsilon\norm{x_0-x_\star}^2$. Hence \eqref{eq:sc-complexity} is sharp for
$0<\varepsilon\le1/2$ up to a factor that depends only on $\kappa_0$.
In particular $1/p$ is the optimal exponent of $\kappa$ among all
predetermined real stepsize schedules, as conjectured by Altschuler and
Parrilo~\cite{AP-I}.

\Cref{tab:sc-bounds} lists lower bounds for the number $n$ of steps of
a schedule with the guarantee
$\norm{x_n-x_\star}^2\le\varepsilon\norm{x_0-x_\star}^2$ of
\cref{thm:sc-complexity}. Before \cite{LY26} the lower bound
$\Omega(\sqrt\kappa\log(1/\varepsilon))$ of \cite{NY83} for
sufficiently large $\kappa$ was the best lower bound for predetermined
schedules. Liu and Ye \cite{LY26} proved that for sufficiently large
$\kappa$ every nonnegative schedule needs at least
$c\kappa^{1/p}\log(1/\varepsilon)/(e^{C\sqrt{\log\kappa}}\log\kappa)$
steps, with absolute constants $c,C>0$. The lower bound of \cite{LY26}
also holds for the guarantee
$f(x_n)-f(x_\star)\le\varepsilon(f(x_0)-f(x_\star))$.
\Cref{thm:sc-complexity} holds without the factor
$e^{C\sqrt{\log\kappa}}\log\kappa$, for every $h\in\R^n$ and every
$\kappa\ge\kappa_0$. The lower bound \eqref{eq:sc-complexity} also
holds for every $h$ that guarantees
$f(x_n)-f(x_\star)\le\varepsilon(f(x_0)-f(x_\star))$ for every
$L$-smooth, $L/\kappa$-strongly convex $f$ and every $x_0$, as proved
in the paragraph after the proof of \cref{thm:sc-complexity} in
\cref{sc:sec:count}.

\begin{table}[!htbp]
\centering\small
\renewcommand{\arraystretch}{1.28}
\begin{tabularx}{\linewidth}{@{}>{\raggedright\arraybackslash}p{.29\linewidth}
 >{\centering\arraybackslash}p{.18\linewidth}X@{}}
\toprule
Reference & Schedules $h$ & Lower bound for $n$\\
\midrule
Nemirovsky--Yudin \cite{NY83}
 & $\R^n$ & $\Omega(\sqrt\kappa\log(1/\varepsilon))$\\
Liu--Ye \cite{LY26}
 & $[0,\infty)^n$ &
 $\Omega\bigl(\kappa^{1/p}\log(1/\varepsilon)\big/
 (e^{C\sqrt{\log\kappa}}\log\kappa)\bigr)$\\
\Cref{thm:sc-complexity} & $\R^n$ &
 $\Omega(\kappa^{1/p}\log(1/\varepsilon))$\\
\bottomrule
\end{tabularx}
\caption{Lower bounds for the number $n$ of steps of a schedule $h$ in
the class indicated under `Schedules $h$' that guarantees
$\norm{x_n-x_\star}^2\le\varepsilon\norm{x_0-x_\star}^2$ for all
$L$-smooth, $L/\kappa$-strongly convex functions, for
$0<\varepsilon<1$. The bounds of \cite{LY26,NY83} are for sufficiently
large $\kappa$. The constant in the row of \cref{thm:sc-complexity}
depends on $\kappa_0$. The number $C$ is an absolute constant.}
\label{tab:sc-bounds}
\end{table}

For $\kappa=1$ every $L$-smooth, $L$-strongly convex $f$ equals
$f(x_\star)+L\norm{x-x_\star}^2/2$ and one step of size $1/L$ gives
$x_1=x_\star$.

\Cref{lem:scaling} relates the trajectories of $f$ and $f/L$.

\begin{lemma}\label{lem:scaling}
Let $L>0$, $\kappa>1$ and let $f:\R^d\to\R$ be differentiable.
Then $f$ is convex and $L$-smooth if and only if $f/L$ is convex and
$1$-smooth. Also $f$ is $L/\kappa$-strongly convex if and only if
$f/L$ is $1/\kappa$-strongly convex. For every $h\in\R^n$ and
$x_0\in\R^d$, the trajectory of $f$ from $x_0$ with steps $h$ is
the trajectory of $f/L$ from $x_0$ with steps $Lh$.
\end{lemma}

\begin{proof}
Dividing by $L$ preserves convexity and divides the Lipschitz constant
of the gradient by $L$. Dividing by $L$ also takes
$f-L\norm{\cdot}^2/(2\kappa)$ to $f/L-\norm{\cdot}^2/(2\kappa)$.
Since $(Lh_i)\nabla(f/L)(x)=h_i\nabla f(x)$ for every $x$,
induction in \eqref{eq:gd} gives the same points $x_1,\ldots,x_n$
for $f$ with steps $h$ and for $f/L$ with steps $Lh$.
\end{proof}

Let $f$ be convex and $L$-smooth. Since $\argmin(f/L)=\argmin f$,
for $x_\star\in\argmin f$ and $x_0\ne x_\star$ we have
\begin{equation}\label{eq:ratio-scaling}
 \frac{2(f(x_n)-f(x_\star))}{L\norm{x_0-x_\star}^2}
 =\frac{2\big((f/L)(x_n)-(f/L)(x_\star)\big)}{\norm{x_0-x_\star}^2}.
\end{equation}
By \cref{lem:scaling} and \eqref{eq:ratio-scaling} the left side of
\eqref{eq:main} does not change when $L$ is replaced by $1$ and $h$ by
$Lh$. By
\cref{lem:scaling} a schedule $h$ guarantees
$\norm{x_n-x_\star}^2\le\varepsilon\norm{x_0-x_\star}^2$ for every
$L$-smooth, $L/\kappa$-strongly convex $f$ and every $x_0$ if and only
if $Lh$ guarantees
$\norm{x_n-x_\star}^2\le\varepsilon\norm{x_0-x_\star}^2$ for every
$1$-smooth, $1/\kappa$-strongly convex $f$ and every $x_0$. It
therefore suffices to prove \cref{thm:main,thm:sc-complexity}
for $L=1$. The constants $c$ and $c(\kappa_0)$ then apply for every $L>0$.

\subsection{Overview of the proof}\label{sec:outline}

Take $L=1$. \Cref{thm:main,thm:sc-complexity} follow from a bound on
the number of large stepsizes. \Cref{sc:thm:count} states that for
every $\kappa>1$, every integer $n\ge1$ and every $h\in\R^n$ satisfying
$\sup_{f,x_0}\norm{x_n-x_\star}^2/\norm{x_0-x_\star}^2<1$ over all
$1$-smooth, $1/\kappa$-strongly convex $f$ and $x_0\ne x_\star$,
\begin{equation}\label{eq:intro-count}
 \#\{i:h_i\ge A\}\le CnA^{-1/p}\qquad(0<A\le\kappa)
\end{equation}
with an absolute constant $C$. In \cref{sec:outline} the value of $C$
may differ from one formula to the next. On $f(x)=\lambda x^2/2$ the
step $i$ multiplies $x$ by $1-\lambda h_i$. When $\kappa\ge2$ and
$\sup_{f,x_0}\norm{x_n-x_\star}^2/\norm{x_0-x_\star}^2<1$,
\eqref{eq:intro-count} bounds the mean of
$\sum_i(\log\lvert1-\lambda h_i\rvert)^-$ over
$1/\kappa\le\lambda\le2/\kappa$ with respect to $d\lambda/\lambda$ by
$Cn\kappa^{-1/p}$. Hence
$\prod_i\lvert1-\lambda h_i\rvert\ge\exp(-Cn\kappa^{-1/p})$ for some
$1/\kappa\le\lambda\le2/\kappa$ and
$\lvert x_n\rvert\ge\exp(-Cn\kappa^{-1/p})\lvert x_0\rvert$ on
$f(x)=\lambda x^2/2$. Thus \cref{prop:sc-distance} holds for
$\kappa\ge2$. \Cref{thm:main} follows from \cref{prop:sc-distance} with
$\kappa=(n+1)^p$ and
$2(f(x_n)-f(x_\star))\ge\norm{x_n-x_\star}^2/\kappa$.

Suppose that \cref{sc:thm:count} fails. Then there are pairs
$(h,\kappa)$ with $h\in\R^n$, $\kappa>1$ and
$\sup_{f,x_0}\norm{x_n-x_\star}^2/\norm{x_0-x_\star}^2<1$ along which
$\max_{0<A\le\kappa}A^{1/p}\#\{i:h_i\ge A\}/n\to\infty$.
\Cref{sc:lem:normalization} changes each pair such that
$A^{1/p}\#\{i:h_i\ge A\}$ is largest at $A=h_n$ among $0<A\le h_n$,
$h_n\to\infty$ and $\kappa/h_n\to\infty$. Zero stepsizes are prepended
until $n-1<h_n^{1/p}\#\{i:h_i\ge h_n\}\le n$. A zero stepsize does not
change the trajectory. After the zero stepsizes are prepended, only
$o(n)$ stepsizes are positive. Limits and the notation $o(n)$ refer to
the sequence of \cref{sc:lem:normalization}.

In the convex case the silver schedule is constructed recursively by
joining two copies of a shorter schedule with a larger step between the
copies. The proof takes $h$ apart in the reverse order. For
$0\le y<\log h_n$ we cut $h$ after every step with stepsize above
$e^y$. The steps with stepsize above $e^y$ are the cuts at threshold
$e^y$. For a cut $j$ let $k$ denote the largest cut below $j$, with
$k=0$ if $j$ is the first cut. Let $I=\{k+1,\ldots,j-1\}$. Then
$h_i\le e^y<\kappa$ for every $i\in I$. Let $m_I$ be the maximum of
$-\kappa\sum_{i\in J}\log(1-h_i/\kappa)$ over the prefixes $J$ of $I$,
including $J=\varnothing$. For positive stepsizes $h_i\le\kappa/2$ we
have $\sum_{i\in I}h_i\le m_I\le\tfrac32\sum_{i\in I}h_i$. Let
$\mathcal P(y)$ be $p/(2-\sqrt2)$ times the sum of $m_I^{1/p}$ over the
cuts $j$ at threshold $e^y$. \Cref{sc:lem:initial-power} gives
$\mathcal P(y)\ge(p-o(1))n\ge n$ for $y<\log h_n$ with
$y/\log h_n\to1$. Since only $o(n)$ stepsizes are positive, we have
$\mathcal P(y)=o(n)$ for bounded $y$. It suffices to prove
$\mathcal P(y)\ge n/2$ at a bounded $y$.

Write $\rho=2^p=1+\sqrt2$. Then $\rho^2=2\rho+1$ and $\rho^{1/p}=2$.
In \cref{sec:one-example} a step with stepsize one divides an interval
with sum $2\rho+1$ into two intervals with sum $\rho$. A step with
stepsize one follows each interval. For an interval followed by a step,
the quotient is the stepsize of the step divided by the sum of the
stepsizes of the interval. The division preserves the product of the
quotients by $(1/\rho)^2=1/(2\rho+1)$ and the sum of the $1/p$ powers
of the sums by $2\rho^{1/p}=(2\rho+1)^{1/p}$.

Take $x_\star=0$. Every pair $(h,\kappa)$ of the sequence satisfies
$\norm{x_n}\le\norm{x_0}$ on every $1$-smooth, $1/\kappa$-strongly
convex function with minimizer zero. The functions of the proof are
$f(x)=\norm{x}^2/(2\kappa)+(1-1/\kappa)M_F(x)$, where
$M_F(x)=\inf_z\{F(z)+\norm{x-z}^2/2\}$ is the Moreau envelope of a
convex function $F$. For sufficiently large $y$ with $e^y<h_n$ let
$A=\min\{e^yy^3,\kappa/2\}$ and $P=e^{2y}/A$. Then
$\log A=y+O(\log y)$ and $0<P<e^y$. Let $j$ be a cut with
$\sum_{k<i\le j}\lvert h_i\rvert\le A$.
\Cref{sc:lem:ordinary} constructs points $x_k,\ldots,x_j$ in dimension
four and the gradients of $M_F$ at $x_k,\ldots,x_{j-1}$. The points are
chosen in the order $x_{j-1},\ldots,x_k$. The step $j$ gives
$\norm{x_j}\ge h_j-1$. \Cref{sc:lem:arc} bounds the norm of
the first two coordinates of the gradients of $M_F$ by $1+C/\log A$.
Combining \eqref{sc:eq:arc-bounds} for coordinates $1$ and $2$ with
\eqref{sc:eq:aux-coordinates} for coordinates $3$ and $4$ gives
$\norm{x_k}\le e^{m_I/\kappa}(1+C/\log A)(m_I+P)$. Hence
\begin{equation}\label{eq:intro-gain}
 \frac{\norm{x_j}}{\norm{x_k}}\ge\frac{h_j}{m_I+P}
                    \exp\bigl(-C/\log A-m_I/\kappa\bigr).
\end{equation}
\Cref{sc:lem:assembly} combines the points and the gradients on
consecutive intervals into one $1$-smooth, $1/\kappa$-strongly convex
function. The cuts $j$ with $\sum_{k<i\le j}\lvert h_i\rvert>A$ and the
adjacent cuts lie in sets $(a,b]$ that are unions of sets $(k,j]$. On
each set $(a,b]$, \cref{sc:lem:bad-block} bounds
$\log(\norm{x_a}/\norm{x_b})$ from above. Summing the logarithm of
\eqref{eq:intro-gain} over the cuts $j$ outside the sets $(a,b]$ and
subtracting the bounds of \cref{sc:lem:bad-block} gives a lower bound
for $\log(\norm{x_n}/\norm{x_0})\le0$. Thus
\cref{sc:prop:working-realization} bounds
$\sum_{h_j>e^y}\log[h_j/(m_I+P)]$ from above. The upper bound contains a
term $C/y$ for each cut.

For $s<t$ every cut at threshold $e^t$ is a cut at threshold $e^s$.
When $e^y$ decreases below a stepsize $h_v$, the steps with stepsize
$h_v$ become cuts. Insert the cuts with stepsize $h_v$ one at a time in
increasing order of index. The cut $v$ divides one interval $I_v$ into
two intervals $I_{vL}$ and $I_{vR}$. Write $m_v,m_{vL},m_{vR}$ for
$m_I$ with $I=I_v,I_{vL},I_{vR}$ and let $\Delta_v$ denote
$p/(2-\sqrt2)$ times $m_v^{1/p}-m_{vL}^{1/p}-m_{vR}^{1/p}$. With the
same $P$, the sum $\sum_{h_j>e^y}\log[h_j/(m_I+P)]$ changes by
\begin{equation}\label{eq:intro-split}
 \log\frac{h_v(m_v+P)}{(m_{vL}+P)(m_{vR}+P)}
\end{equation}
and $\mathcal P$ changes by $-\Delta_v$. For $0\le s<t<\log h_n$ the
difference $\mathcal P(t)-\mathcal P(s)$ is the sum of $\Delta_v$ over
the steps $v$ with $s<\log h_v\le t$. For $P=0$ the logarithm
\eqref{eq:intro-split} and $\Delta_v$ are zero at the silver
proportions $m_{vL}=m_{vR}=\rho h_v$, $m_v=\rho^2h_v$.
\Cref{lem:max-plus-excess} bounds \eqref{eq:intro-split} minus
$\Delta_v/h_v^{1/p}$ from below by
$c\min\{1,\sum_{w=vL,vR}(m_w/h_v-\rho)^2\}-C(P/h_v)^{1/p}-Ch_v/\kappa$.
Hence the integral of the upper bound of
\cref{sc:prop:working-realization} against $e^{y/p}\,dy$ over
$s\le y\le t$ bounds $\mathcal P(t)-\mathcal P(s)$. For sufficiently
large $s$ and $s<t\le4s$ with $e^tt^3\le h_n$,
\cref{sc:lem:excess-sum} gives
$\sum_{s<\log h_v\le t}\lvert\Delta_v\rvert\le Cn$.

The bound $Cn$ of \cref{sc:lem:excess-sum} comes from the term $C/y$
for each cut in \cref{sc:prop:working-realization}. The term $C/y$
comes from the factor $\exp(-C/\log A)$ in \eqref{eq:intro-gain}. The
number of ranges $(s,4s]$ between a bounded threshold and $\log h_n$ is
unbounded. Thus \cref{sc:lem:excess-sum} does not give
$\mathcal P(y)\ge n/2$ at a bounded $y$. \Cref{sc:lem:ordinary}
constructs the points for one interval at a time.
\Cref{sc:lem:four-phase} constructs the points for four consecutive
intervals at once. For each of the four intervals, the step after the
interval lies in a set $B$ and has stepsize at most $e^y$. On each of
the four intervals the gradient of $M_F$ is constant. For the points of \cref{sc:lem:four-phase},
$\log(\norm{x_j}/\norm{x_k})$ is a weighted mean of the three
logarithms on the right side of \eqref{sc:eq:four-gain} up to
$O(y^{-2})$. The
term $O(y^{-2})$ replaces $C/\log A$. Every $v\in B$ has
$\log h_v$ above an absolute constant and $h_v(\log h_v)^3\le h_n$ and
satisfies the silver proportions up to factors
$1+O((\log h_v)^{-1/4})$. \Cref{sc:lem:excess-sum} bounds the
sum of $\min\{h_v,e^t\}^{1/p}$ over the positive steps $v$ for which
$v$ or the first cut inserted in $I_{vL}$ or in $I_{vR}$ lies outside
$B$.

Let $N$ be the smallest integer with $2^N\ge(\log t)^2$.
\Cref{sc:lem:two-ranges} uses a set $S$ of cuts at threshold $e^y$. For
$j\in S$ insert cuts in $I$ until $I$ is divided into $2^N$ intervals
$I_a$. Every inserted cut lies in $B$. Let $I_\omega$ be the interval with the largest
indices among the $2^N$ intervals $I_a$. The step after $I_\omega$ is
the cut $j$ and $h_j>e^y$. For every $a\ne\omega$, \cref{sc:lem:two-ranges} uses
\cref{sc:lem:four-phase} on $I_a$ and the step after $I_a$. On
$I_\omega$ and the cut $j$, \cref{sc:lem:two-ranges} uses
\cref{sc:lem:ordinary} with $I_\omega$ divided into one, two and four
intervals. For the cut $j$ the terms
$2^NO(y^{-2})$ of \cref{sc:lem:four-phase} and the term $C/y$ of
$I_\omega$ sum to at most $C/y$. The integral over
$s\le y\le t$ then bounds a weighted integral of
$\mathcal P(t-\theta)-\mathcal P(s-\theta)$ over
$0\le\theta\le(N+2)\log\rho$ with total weight of order
$2^N\ge(\log t)^2$ by $Cn$. For sufficiently large $t$ with
$e^tt^3\le h_n$ we obtain for some $t/2\le s\le2t/3$ that the infimum
of $\mathcal P$ on $[s-\log s,s]$ is at least the infimum of
$\mathcal P$ on $[t-\log t,t]$ minus $Cn/(\log t)^2$. Applying the
inequality again with $s$ as the new $t$ and repeating reaches a
bounded threshold. The numbers $Cn/(\log t)^2$ sum to at most
$n/2$. Hence $\mathcal P(y)\ge n/2$ at a bounded $y$. The contradiction
with $\mathcal P(y)=o(n)$ proves \cref{sc:thm:count}.

\Cref{sc:sec:count} proves \cref{thm:main,thm:sc-complexity} from
\cref{sc:thm:count} and gives an example of one split.
\Cref{sec:single-split} fixes the sequence of schedules, defines $m_I$
and $\mathcal P$ and proves \cref{lem:max-plus-excess}.
\Cref{sec:backward} constructs the points for one interval and proves
\cref{sc:lem:assembly}. \Cref{sc:sec:cuts} treats all intervals at one
threshold. \Cref{sc:sec:deviations} integrates over the thresholds.
\Cref{sec:consecutive-splits} constructs the points for four
consecutive intervals and completes the proof of \cref{sc:thm:count}.

\subsection{Related work}\label{sec:previous}

Related work concerns stepsize schedules for smooth convex functions,
quadratic functions, methods with momentum and the anytime and random
stepsize settings outside \cref{thm:main,thm:sc-complexity}. The
performance estimation framework of Drori and Teboulle \cite{DT14} and
of Taylor, Hendrickx and Glineur \cite{THG17} made worst-case
guarantees for a given schedule on smooth convex functions computable
by semidefinite programming. Altschuler \cite{Alt18} determined the
optimal stepsizes for two and three steps with respect to worst-case
distance contraction in the strongly convex setting. Das Gupta, Van
Parys and Ryu \cite{DVR24} computed numerically optimized schedules for
$n\le50$. For functions with bounded sublevel set
$\{x:f(x)\le f(x_0)\}$, Grimmer \cite{Gri24} proved that certain
periodic schedules with long steps improve the leading constant of the
worst-case objective-gap guarantee of the constant stepsize $1/L$ at
the end of each period. Bok and Altschuler \cite{BA25} extended the
rates $n^{-p}$ and $\kappa^{1/p}\log(1/\varepsilon)$ to projected and
proximal gradient descent.

For every constant stepsize, one of the functions $Lx^2/2$ and
$Lx^2/(2\kappa)$ on $\R$ has
$\lvert x_i\rvert\ge\lvert x_{i-1}\rvert(\kappa-1)/(\kappa+1)$ at every
step. On $f(x)=\lambda x^2/2$ with $x_0\ne0$, the ratio $x_n/x_0$ is a
polynomial in $\lambda$ of degree at most $n$ with value one at
$\lambda=0$. Among the polynomials of degree at most $n$ with value one
at $0$, a rescaled Chebyshev polynomial has the least maximum of the
absolute value on $[L/\kappa,L]$. The least maximum is at least
$[(\sqrt\kappa-1)/(\sqrt\kappa+1)]^n$. Hence for every $h\in\R^n$,
\begin{equation}\label{eq:chebyshev}
 \max_{L/\kappa\le\lambda\le L}\,\prod_{i=1}^n\lvert1-\lambda h_i\rvert
 \ge\left(\frac{\sqrt\kappa-1}{\sqrt\kappa+1}\right)^n .
\end{equation}
By \eqref{eq:chebyshev}, for sufficiently large $\kappa$ and
$0<\varepsilon\le1/2$, no $h\in\R^n$ improves the order
$\sqrt\kappa\log(1/\varepsilon)$ of the stepsizes of \cite{You53} on
quadratic functions. The order in which the stepsizes of \cite{You53}
are used affects the numerical stability \cite{LF71}. Agarwal, Goel and
Zhang \cite{AGZ21} showed that a fractal permutation of the stepsizes
of \cite{You53} is stable. Goujaud, Scieur, Dieuleveut, Taylor and
Pedregosa \cite{GSDTP22} analyzed momentum with cyclical stepsizes for
Hessians with a spectral gap. Kelner, Marsden, Sharan, Sidford, Valiant
and Yuan \cite{KMSSVY22} gave gradient methods for sums of
non-interacting smooth strongly convex functions with different
scales.

For sufficiently large $\kappa$ and $0<\varepsilon\le1/2$, the order
$\sqrt\kappa\log(1/\varepsilon)$ is a lower bound for all first-order
methods \cite{NY83}. Momentum methods go back to Polyak \cite{Pol64}
and Nesterov \cite{Nes83}. For smooth convex functions and final
objective gap, the method of Kim and Fessler \cite{KF16} attains the
lower bound of Drori \cite{Dro17} for first-order methods. For smooth
strongly convex functions and the distance to $x_\star$, the lower
bound of Drori and Taylor \cite{DT22} equals the guarantee of the
method of Taylor and Drori \cite{TD23}. Thus a lower bound with an
exponent of $\kappa$ larger than $1/2$ requires functions that are not
quadratic and fails for first-order methods with momentum.

In the anytime setting the steps $h_1,\ldots,h_n$ are the first $n$
terms of one infinite schedule for every $n$. Kornowski and Shamir
\cite{KS24} asked whether an infinite schedule can accelerate the rate
$n^{-1}$ for every $n$. Zhang, Lee, Du and Chen \cite{ZLDC25} obtained
the objective-gap rate $n^{-2p/(1+p)}$. For strongly convex functions
\cite{ZLDC25} gives an anytime objective-gap bound of the form
$\exp(-cn/\kappa^{(1+p)/(2p)})$. Tsai, Fatkhullin, Zhang and He
\cite{TFZH26} proved lower bounds in the anytime setting for positive
stepsizes. Ye and Liu \cite{ye2026almost} proved a lower bound that
matches the exponent $2p/(1+p)$ for nonnegative stepsizes. He and Cheng
\cite{he2026optimal} extended the lower bound of \cite{ye2026almost} to
real stepsizes. Random stepsizes form another setting. For separable
$L$-smooth, $L/\kappa$-strongly convex functions with $x_0\ne x_\star$,
independent stepsizes whose reciprocals have the arcsine distribution
on $[L/\kappa,L]$ give
$(\norm{x_n-x_\star}/\norm{x_0-x_\star})^{1/n}\to(\sqrt\kappa-1)/(\sqrt\kappa+1)$
almost surely \cite{AP24}. The present paper concerns one deterministic schedule and
a fixed $n$.

\section{Step counts and convergence lower bounds}\label{sc:sec:count}

Take $L=1$ by \cref{lem:scaling}. In this section we reduce
\cref{thm:main,thm:sc-complexity} to \cref{sc:thm:count}, a bound for
the number of stepsizes above each threshold.
\Cref{thm:main,thm:sc-complexity} follow from \cref{prop:sc-distance},
an estimate on strongly convex functions. Under the strict contraction
condition $\sup_{f,x_0}\norm{x_n-x_\star}^2/\norm{x_0-x_\star}^2<1$ on
the full $1$-smooth, $1/\kappa$-strongly convex class,
\cref{sc:thm:count} bounds the step counts and testing scalar quadratic
functions then proves \cref{prop:sc-distance}.
\Crefrange{sec:single-split}{sec:consecutive-splits} prove
\cref{sc:thm:count}.

\begin{proposition}\label{prop:sc-distance}
There is an absolute constant $C>0$ such that for every $\kappa>1$,
every integer $n\ge1$ and every $h\in\R^n$,
\begin{equation}\label{eq:sc-distance}
 \sup_{f,x_0}\frac{\norm{x_n-x_\star}^2}{\norm{x_0-x_\star}^2}
 \ge (1-\kappa^{-1/p})^{Cn},
\end{equation}
where the supremum is over all integers $d\ge1$, all $1$-smooth,
$1/\kappa$-strongly convex $f:\R^d\to\R$ and all
$x_0\in\R^d\setminus\{x_\star\}$.
\end{proposition}

\subsection{Proofs of the two main theorems}

\begin{proof}[Proof of \cref{thm:main}]
Fix $n$ and choose $\kappa=(n+1)^p$. Then
$(1-\kappa^{-1/p})^{Cn}=(1+1/n)^{-Cn}\ge e^{-C}$.
Every $1$-smooth, $1/\kappa$-strongly convex function $f:\R^d\to\R$ is
convex and $1$-smooth. Strong convexity and $\nabla f(x_\star)=0$ give
\[
 2\bigl(f(x_n)-f(x_\star)\bigr)\ge
 \norm{x_n-x_\star}^2/\kappa .
\]
Taking suprema and applying \cref{prop:sc-distance} gives a lower
bound of $e^{-C}(n+1)^{-p}\ge 2^{-p}e^{-C}n^{-p}$ in \eqref{eq:main}.
\end{proof}

\begin{proof}[Proof of \cref{thm:sc-complexity}]
Let $h$ guarantee
$\norm{x_n-x_\star}^2\le\varepsilon\norm{x_0-x_\star}^2$ for every
$d\ge1$, every $1$-smooth, $1/\kappa$-strongly convex $f:\R^d\to\R$ and
every $x_0\in\R^d$. Then \cref{prop:sc-distance} gives
$(1-\kappa^{-1/p})^{Cn}\le\varepsilon$.
Since $-\log(1-\kappa^{-1/p})\le\kappa^{-1/p}/(1-\kappa^{-1/p})$ and $\kappa\ge\kappa_0>1$,
\[
 \log(1/\varepsilon)\le -Cn\log(1-\kappa^{-1/p})
 \le \frac{Cn\kappa^{-1/p}}{1-\kappa_0^{-1/p}} .
\]
Rearranging gives \eqref{eq:sc-complexity}.
\end{proof}

Suppose instead that $h$ guarantees
$f(x_n)-f(x_\star)\le\varepsilon(f(x_0)-f(x_\star))$ for every $d\ge1$,
every $1$-smooth, $1/\kappa$-strongly convex $f:\R^d\to\R$ and every
$x_0\in\R^d$. Repeating $h$ for $r$ rounds gives
$f(x_{nr})-f(x_\star)\le\varepsilon^r(f(x_0)-f(x_\star))$
\cite{LY26}. Every $1$-smooth, $1/\kappa$-strongly convex $f$ satisfies
\[
 \norm{x-x_\star}^2/\kappa\le 2(f(x)-f(x_\star))\le\norm{x-x_\star}^2 .
\]
Hence
$\norm{x_{nr}-x_\star}^2\le\kappa\varepsilon^r\norm{x_0-x_\star}^2$.
\Cref{prop:sc-distance} for the repeated schedule of length $nr$ gives
$(1-\kappa^{-1/p})^{Cnr}\le\kappa\varepsilon^r$. Taking $r$th roots
and letting $r\to\infty$ gives $(1-\kappa^{-1/p})^{Cn}\le\varepsilon$.
The proof of \cref{thm:sc-complexity} then gives
\eqref{eq:sc-complexity}.

\subsection{From step counts to distance ratios}\label{sc:sec:step-counts}

Let $\kappa>1$. For $h\in\R^n$ we use the condition
\begin{equation}\label{eq:guarantee}
 \sup_{f,x_0}\frac{\norm{x_n-x_\star}^2}{\norm{x_0-x_\star}^2}<1
\end{equation}
where the supremum is over all integers $d\ge1$, all $1$-smooth,
$1/\kappa$-strongly convex $f:\R^d\to\R$ and all
$x_0\in\R^d\setminus\{x_\star\}$.

\begin{theorem}\label{sc:thm:count}
There is an absolute constant $C>0$ such that for every $\kappa>1$,
every integer $n\ge1$ and every $h\in\R^n$ with \eqref{eq:guarantee},
\begin{equation}\label{sc:eq:count-goal}
 \#\{i:h_i\ge A\}\le C n A^{-1/p}\qquad(0<A\le\kappa).
\end{equation}
\end{theorem}

From here on $c$ and $C$ denote positive absolute constants unless a
dependence on fixed auxiliary constants is stated explicitly. The
values of $c$ and $C$ may differ from one formula to the next. For
$f(x)=\lambda x^2/2$ on $\R$ with $x_0=1$,
\begin{equation}\label{sc:eq:quadratic-iterate}
 x_n=\prod_{i=1}^n(1-\lambda h_i).
\end{equation}

Suppose that $h$ satisfies \eqref{eq:guarantee}. Let
$1/\kappa\le\lambda\le1$ with $\lambda h_i\ne1$ for every $i$. By
\eqref{sc:eq:quadratic-iterate},
$\sum_i\log\lvert1-\lambda h_i\rvert\le0$.
In $\sum_i\log[(1+\lambda\lvert h_i\rvert)/\lvert1-\lambda h_i\rvert]$
all terms with $h_i\le0$ vanish. Adding $\sum_i\log\lvert1-\lambda h_i\rvert\le0$
shows that $\sum_i\log(1+\lambda\lvert h_i\rvert)$ is at most
$\sum_{h_i>0}\log[(1+\lambda h_i)/\lvert1-\lambda h_i\rvert]$. Because
$1+\lambda h_i\ge\max\{1,\lvert1-\lambda h_i\rvert\}$ for $h_i>0$ and
$\lvert1-\lambda h_i\rvert\ge1$ for $h_i\le0$, every term satisfies
$(\log\lvert1-\lambda h_i\rvert)^-\le\log[(1+\lambda\lvert h_i\rvert)/\lvert1-\lambda h_i\rvert]$.
Hence each of $\sum_i\log(1+\lambda\lvert h_i\rvert)$ and
$\sum_i(\log\lvert1-\lambda h_i\rvert)^-$ is at most
$\sum_{h_i>0}\log[(1+\lambda h_i)/\lvert1-\lambda h_i\rvert]$.
For $\tau_0>0$,
\[
 \int_0^{\tau_0}\log\frac{1+\tau}{|1-\tau|}\,\frac{d\tau}{\tau}
 \le C\min\{\tau_0,1\}.
\]
Indeed, $\log((1+\tau)/\lvert1-\tau\rvert)/\tau$ is bounded near zero
and the singularity at one is integrable. The substitution
$\tau\mapsto1/\tau$ shows that the integrals of
$\tau^{-1}\log[(1+\tau)/\lvert1-\tau\rvert]$ over $(0,1)$ and over
$(1,\infty)$ are equal. At most $n$ numbers $\lambda$ have
$\lambda h_i=1$ for some $i$. For $1\le A\le\kappa$, integrating
$\sum_{h_i>0}\log[(1+\lambda h_i)/\lvert1-\lambda h_i\rvert]\,d\lambda/\lambda$
with $\tau=\lambda h_i$ over $[1/\kappa,1/A]\subset[1/\kappa,1]$ gives
\begin{equation}\label{sc:eq:quadratic-average}
 \int_{1/\kappa}^{1/A}\sum_i\left(\log(1+\lambda|h_i|)
           +(\log|1-\lambda h_i|)^-\right)\frac{d\lambda}{\lambda}
 \le \frac CA\sum_i\min\{h_i^+,A\}.
\end{equation}
The proof of \cref{prop:sc-distance} uses $(\log|1-\lambda h_i|)^-$ in
\eqref{sc:eq:quadratic-average} and \cref{sc:sec:setup} uses both
$\log(1+\lambda|h_i|)$ and $(\log|1-\lambda h_i|)^-$ in
\eqref{sc:eq:quadratic-average} to bound the stepsize tail sums.

\begin{proof}[Proof of \cref{prop:sc-distance}]
We may assume that
$\sup_{f,x_0}\norm{x_n-x_\star}^2/\norm{x_0-x_\star}^2<1$. Hence $h$
satisfies \eqref{eq:guarantee}. By \eqref{sc:eq:count-goal},
\[
  \sum_i\min\{h_i^+,A\}=\int_0^A\#\{i:h_i\ge\tau\}\,d\tau
 \le CnA^{1-1/p}\qquad(0<A\le\kappa).
\]
Let $\kappa\ge2$ and take $A=\kappa/2$ in
\eqref{sc:eq:quadratic-average}. Since
$\int_{1/\kappa}^{2/\kappa}d\lambda/\lambda=\log2$, some
$1/\kappa\le\lambda\le2/\kappa$ satisfies
$\sum_i(\log\lvert1-\lambda h_i\rvert)^-\le Cn\kappa^{-1/p}$. For
$f(x)=\lambda x^2/2$ with $x_0=1$, \eqref{sc:eq:quadratic-iterate}
gives $-\log\lvert x_n\rvert\le Cn\kappa^{-1/p}$. Squaring and using
$\log(1-\kappa^{-1/p})\le-\kappa^{-1/p}$ proves \eqref{eq:sc-distance}.

Let $1<\kappa<2$. By \eqref{eq:chebyshev} with $L=1$ and
\eqref{sc:eq:quadratic-iterate}, some $1/\kappa\le\lambda\le1$ has
$\lvert x_n\rvert\ge[(\sqrt\kappa-1)/(\sqrt\kappa+1)]^n$. Since
$1-\kappa^{-1/p}<1/2$,
\[
 \frac{\sqrt\kappa-1}{\sqrt\kappa+1}
 =\frac{1-1/\kappa}{(1+1/\sqrt\kappa)^2}
 \ge\frac{1-\kappa^{-1/p}}4\ge(1-\kappa^{-1/p})^3 .
\]
Squaring gives \eqref{eq:sc-distance}.
\end{proof}

\subsection{Distance ratios at a split}\label{sec:one-example}

Write $\rho=2^p=1+\sqrt2$. Then $\rho^2=2\rho+1$ and $\rho^{1/p}=2$.
Take positive stepsizes below one with sum $\rho$, then a stepsize one,
then positive stepsizes below one with sum $\rho$, then a stepsize one.
For an interval followed by a step, divide the stepsize of the step by
the sum of the stepsizes of the interval.
Cut the schedule after each step with stepsize one. Then there are two
intervals before the cuts with sum $\rho$ and quotients $1/\rho$ and
$1/\rho$. Cut the schedule after the last step only. Then there is one
interval before the last step with sum $\rho+1+\rho$ and quotient
$1/(2\rho+1)$. We have
\begin{gather}
 \frac1\rho\cdot\frac1\rho=\frac1{2\rho+1},\label{eq:outline-split}\\
 \rho^{1/p}+\rho^{1/p}=(2\rho+1)^{1/p}.\label{eq:outline-powers}
\end{gather}
By \eqref{eq:outline-split} the product of the quotients is the same
for one cut and for two cuts. By \eqref{eq:outline-powers} the sum of
the $1/p$ powers of the sums of the intervals is the same for one cut
and for two cuts.

For an interval $I$ with $h_i<\kappa$ for every $i\in I$, the number
$m_I\ge0$ of \eqref{sc:eq:mass-def} replaces the sum of the stepsizes
of $I$. A negative stepsize increases the distance on $x^2/(2\kappa)$.
For positive stepsizes $h_i\le A\le\kappa/2$ we have
$\sum_{i\in I}h_i\le m_I\le(1+A/\kappa)\sum_{i\in I}h_i$. Let $A$ be
sufficiently large and let $I$ be the interval between consecutive cuts
$k<j$ with $\sum_{k<i\le j}\lvert h_i\rvert\le A\le\kappa/2$. Take
$A^{19/20}\le P\le h_j$. \Cref{sc:lem:ordinary} constructs a trajectory
with $\norm{x_j}\ge h_j-1$ and
$\norm{x_k}\le e^{m_I/\kappa}(1+C/\log A)(m_I+P)$. Hence
$\norm{x_j}/\norm{x_k}$ is at least $h_j/(m_I+P)$ up to the factor
$\exp(-C/\log A-m_I/\kappa)$.

Inequality \eqref{sc:eq:count-goal} is sharp for the silver schedules.
Let $n+1\ge4$ be a power of two. The silver schedule of length $n$ contains the stepsize
$1+\rho^{i-2}$ for $(n+1)2^{-i}$ indices, for every
$1\le i\le\log_2(n+1)$ \cite{altschuler2025silver}. Hence
$(n+1)2^{1-i}-1$ stepsizes are at least $A=1+\rho^{i-2}$. Since
$\rho^{1/p}=2$, we have
$cnA^{-1/p}\le(n+1)2^{1-i}-1\le CnA^{-1/p}$. Write
$\kappa=1+\rho^{-2}(n+1)^p$ for the largest stepsize of the silver
schedule of length $n$. Then
$A\le\kappa\le2\rho^{-2}(n+1)^p<((n+1)/2)^p\le n^p$. On $1$-smooth
convex functions the silver schedule satisfies
$f(x_n)-f(x_\star)\le\norm{x_0-x_\star}^2/(2n^p)$
\cite{altschuler2025silver}. Every $1$-smooth, $1/\kappa$-strongly
convex $f$ satisfies
$\norm{x_n-x_\star}^2\le2\kappa(f(x_n)-f(x_\star))$. Hence
$\norm{x_n-x_\star}^2\le\kappa n^{-p}\norm{x_0-x_\star}^2$ and the
silver schedule satisfies \eqref{eq:guarantee} for $\kappa$.

\section{Step counts and estimates for a single split}\label{sec:single-split}

\Cref{sc:sec:count} reduces \cref{thm:main,thm:sc-complexity} to
\cref{sc:thm:count}. In \crefrange{sec:single-split}{sec:consecutive-splits}
we suppose that \cref{sc:thm:count} fails and derive a contradiction.
In this section we fix a sequence of schedules by
\cref{sc:lem:normalization}, define $m_I$ and $\mathcal P(y)$ and
prove \cref{lem:max-plus-excess}. Since only $o(n)$
stepsizes are positive, \eqref{sc:eq:power-small-threshold} gives
$\mathcal P(y)=o(n)$ for bounded $y$. The contradiction is
$\mathcal P(y)\ge n/2$ at a bounded $y$.
\subsection{Choosing the threshold and the condition number}\label{sc:sec:setup}

Let $h\in\R^n$ with $h_n>0$. Suppose that
$A^{1/p}\#\{i:h_i\ge A\}$ is largest at $A=h_n$ among $0<A\le h_n$ and
that
\begin{equation}\label{sc:eq:normalized-length}
 n-1<h_n^{1/p}\,\#\{i:h_i\ge h_n\}\le n.
\end{equation}
For $0<A\le h_n$ we then have
\begin{align}
 \#\{i:h_i\ge A\}&\le nA^{-1/p},\label{sc:eq:count-bound}\\
 \sum_{0<h_i\le A}h_i\le\sum_i\min\{h_i^+,A\}
  &=\int_0^A\#\{i:h_i\ge\tau\}\,d\tau
 \le\frac{nA^{1-1/p}}{1-1/p}.\label{sc:eq:positive-tail}
\end{align}

If \cref{sc:thm:count} fails, \cref{sc:lem:normalization} gives a
sequence of schedules with \eqref{sc:eq:normalized-length}. Conditions
\eqref{sc:eq:standing-log} and \eqref{sc:eq:top-quadratic} are used in
the proof of \cref{sc:prop:working-realization}.

\begin{lemma}\label{sc:lem:normalization}
Suppose that \cref{sc:thm:count} fails. Then there is a sequence of
pairs $(h,\kappa)$ of a schedule $h\in\R^n$ with $h_n>0$ and a number
$\kappa>1$, where $n$ depends on the pair, such that:
\begin{enumerate}
\item\label{it:norm-guar} every pair $(h,\kappa)$ satisfies
\eqref{eq:guarantee} for $\kappa$, the product
$A^{1/p}\#\{i:h_i\ge A\}$ is largest at $A=h_n$ among $0<A\le h_n$ and
\eqref{sc:eq:normalized-length} holds.
\item\label{it:norm-scale} along the sequence,
\begin{gather}
 \#\{i:h_i>0\}=o(n),\label{sc:eq:standing-sparse}\\
 \kappa/h_n\to\infty,\label{sc:eq:standing}\\
 \kappa\le h_n\log h_n.\label{sc:eq:standing-log}
\end{gather}
\item\label{it:norm-top} $h_i\ne\kappa$ for every $i$ and
\begin{equation}\label{sc:eq:top-quadratic}
 \sum_i\left(\log\!\left(1+\frac{|h_i|}{\sqrt{\kappa h_n}}\right)
       +\bigl(\log|1-h_i/\kappa|\bigr)^-\right)
 \le\frac{Cn h_n^{-1/p}}{\log(\kappa/h_n)}
\end{equation}
with an absolute constant $C$.
\end{enumerate}
\end{lemma}
\begin{proof}
Since \cref{sc:thm:count} fails, there are numbers $\kappa>1$,
schedules $h\in\R^n$ with \eqref{eq:guarantee} for $\kappa$ and numbers
$0<A\le\kappa$ such that $R=A^{1/p}\#\{i:h_i\ge A\}/n\to\infty$.
Pass to a tail of the sequence where $R>1$. The count
$\#\{i:h_i\ge A\}$ does not decrease when $A$ decreases. Replace $A$
by $A/R^{p/2}\le A$. Then
$A^{1/p}\#\{i:h_i\ge A\}/n\ge\sqrt R\to\infty$ and $\kappa/A\to\infty$.

We change $h$ three times. After each change $h$ and $n$ denote the
changed schedule and the length of the changed schedule.

First replace $h$ by a schedule that consists of copies of $h$ and
satisfies
$\sup_{f,x_0}\norm{x_n-x_\star}/\norm{x_0-x_\star}<1/(1+A)$. The length
$n$ and every count $\#\{i:h_i\ge\tau\}$ with $\tau>0$ are multiplied
by the number of copies.

Second append one step with a stepsize $h_{\text{new}}\in(0,A]$ that
maximizes $\tau^{1/p}\#\{i:h_i\ge\tau\}$ over $0<\tau\le A$. A step of size at
most $A$ multiplies every distance to the minimizer by at most $1+A$.
Hence $h$ satisfies \eqref{eq:guarantee} and $h_n=h_{\text{new}}\le A$. For $0<\tau\le h_n$, appending $h_n$
increases $\tau^{1/p}\#\{i:h_i\ge\tau\}$ by
$\tau^{1/p}\le h_n^{1/p}$. Therefore
$\tau^{1/p}\#\{i:h_i\ge\tau\}$ is largest at $\tau=h_n$ among
$0<\tau\le h_n$. Appending one step at most doubles $n$. Thus
\[
 h_n^{1/p}\#\{i:h_i\ge h_n\}/n\to\infty .
\]

Third pass to a tail of the sequence where
$h_n^{1/p}\#\{i:h_i\ge h_n\}\ge n$ and prepend zero stepsizes until $n$
is the smallest integer with $n\ge h_n^{1/p}\#\{i:h_i\ge h_n\}$. A zero
stepsize leaves the iterate unchanged. Hence \eqref{eq:guarantee} still
holds. By the choice of $n$, \eqref{sc:eq:normalized-length} holds. The
number of positive stepsizes is at most the length before the zero
stepsizes were prepended. Since $h_n^{1/p}\#\{i:h_i\ge h_n\}$ divided by
the length before the zero stepsizes were prepended tends to infinity,
\eqref{sc:eq:standing-sparse} holds. Since $h_n\le A$, we have
$\kappa/h_n\to\infty$. By \eqref{sc:eq:normalized-length} and
\eqref{sc:eq:standing-sparse}, $h_n\to\infty$.

Pass to a tail of the sequence where $h_n>e$ and $\kappa>h_n$. Replace
$\kappa$ by $\min\{\kappa,h_n\log h_n\}>h_n$. The class of
$1$-smooth, $1/\kappa$-strongly convex functions becomes smaller.
Hence \eqref{eq:guarantee} still holds. Since $h_n\to\infty$, also
$\kappa/h_n\to\infty$ still holds. For $\lambda>0$ outside the finite
set $\{1/h_i:h_i>0\}$ define
\[
 G(\lambda)=\sum_i\left(\log\bigl(1+\lvert h_i\rvert\sqrt{\lambda/h_n}\bigr)
       +\bigl(\log\lvert1-\lambda h_i\rvert\bigr)^-\right).
\]
For $1/\kappa\le\lambda\le1/\sqrt{\kappa h_n}$ we have
$\sqrt{\lambda/h_n},\lambda\in[1/\kappa,1/h_n]$. Since
$d\log\sqrt{\lambda/h_n}=\tfrac12\,d\log\lambda$,
\eqref{sc:eq:quadratic-average} at $A=h_n\ge1$ and
\eqref{sc:eq:positive-tail} give
$\int_{1/\kappa}^{1/\sqrt{\kappa h_n}}G(\lambda)\,d\lambda/\lambda\le Cn h_n^{-1/p}$.
Also
\[
 \int_{1/\kappa}^{1/\sqrt{\kappa h_n}}\frac{d\lambda}\lambda
 =\tfrac12\log(\kappa/h_n)>0 .
\]
Hence some $\lambda\in[1/\kappa,1/\sqrt{\kappa h_n}]$ outside
$\{1/h_i:h_i>0\}$ satisfies $G(\lambda)\le Cnh_n^{-1/p}/\log(\kappa/h_n)$.
Replace $\kappa$ by $1/\lambda\in[\sqrt{\kappa h_n},\kappa]$. Since
$1/\lambda\le\kappa$ does not increase $\log(\kappa/h_n)$ and
$1/\lambda\ge\sqrt{\kappa h_n}>h_n$, we obtain $h_i\ne\kappa$ for every
$i$, \eqref{sc:eq:top-quadratic}, \eqref{sc:eq:standing} and
\eqref{sc:eq:standing-log}. Also \eqref{eq:guarantee} still holds.
\end{proof}

From here on $h\in\R^n$ and $\kappa$ denote a schedule and a number
of the sequence in \cref{sc:lem:normalization}. Unless another limit is
named, the notation $o(\cdot)$ and limits refer to the sequence in
\cref{sc:lem:normalization}. In particular $h_n\to\infty$. By
\eqref{sc:eq:standing} we may assume $h_n\le\kappa/4$. By
\cref{sc:lem:normalization}\ref{it:norm-top}, $h_i\ne\kappa$ for every
$i$.

For $1\le A\le h_n$ with $2A\le\kappa$, \eqref{sc:eq:positive-tail}
gives $(C/A)\sum_i\min\{h_i^+,A\}\le CnA^{-1/p}$. Since
$\log(1+|h_i|/A)\le2\log(1+\lambda|h_i|)$ for $1/(2A)\le\lambda\le1/A$,
restricting \eqref{sc:eq:quadratic-average} to $1/(2A)\le\lambda\le1/A$
gives
\begin{equation}\label{sc:eq:quadratic-tails}
 \sum_i\log(1+|h_i|/A)\le CnA^{-1/p}.
\end{equation}

\subsection{The intervals at each threshold}\label{sc:sec:interval-tree}\label{sc:sec:curvature}\label{sc:sec:gaps}

Call a set of consecutive indices an \emph{interval}. For integers
$a<b$ write $(a,b]=\{i\in\mathbb Z:a<i\le b\}$. For real numbers $s<t$
the set $(s,t]$ is a real interval. Call $\varnothing$ and every set
$\{i\in I:i\le u\}$ with $u\in I$ a \emph{prefix} of the interval
$I$. Call an index
$1\le i\le n$ a \emph{step} and an index $v$ with $h_v>0$ a
\emph{positive step}. In every sum or condition involving $\log h_v$ or
$\log h_b$, the index runs over positive steps.

Starting with $\{1,\ldots,n\}$, remove the positive steps one at a time
in decreasing order of stepsize with ties broken in increasing order of
index. For a positive step $v$ let $I_v$ be the largest interval that
contains $v$ and no step removed before $v$. Define
\begin{equation}\label{eq:tree-split}
 I_{vL}=\{i\in I_v:i<v\},\qquad I_{vR}=\{i\in I_v:i>v\}.
\end{equation}
If $I_{vL}$ contains a positive step, let $vL$ also denote the first
step removed from $I_{vL}$. Then the interval of the step $vL$ is
$I_{vL}$. If $I_{vR}$ contains a positive step, let $vR$ also denote
the first step removed from $I_{vR}$. Then the interval of the step
$vR$ is $I_{vR}$. When present, call $vL$ and $vR$ the \emph{children}
of $v$. Call $v$ the \emph{parent} of $I_{vL}$, of $I_{vR}$ and of each
child of $v$.

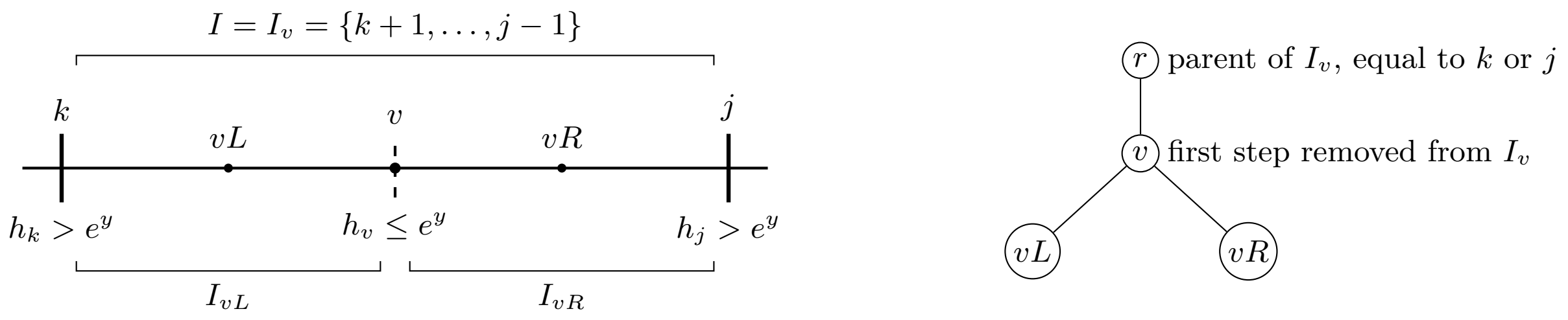
\begin{figure}[htbp]
\centering
\begin{tikzpicture}[font=\footnotesize,x=1cm,y=1cm]
 \draw[thick] (0,1.4)--(7.6,1.4);
 \foreach \x/\lab in {0.4/k, 7.2/j} {
  \draw[very thick] (\x,1.05)--(\x,1.75);
  \node[above] at (\x,1.75) {$\lab$};
  \node[below] at (\x,1.05) {$h_\lab>e^y$};
 }
 \draw[thick,dashed] (3.8,1.1)--(3.8,1.7);
 \fill (3.8,1.4) circle (1.6pt);
 \node[above] at (3.8,1.72) {$v$};
 \node[below] at (3.8,1.08) {$h_v\le e^y$};
 \fill (2.1,1.4) circle (1.3pt) node[above=2pt] {$vL$};
 \fill (5.5,1.4) circle (1.3pt) node[above=2pt] {$vR$};
 \draw (0.55,0.45)--(0.55,0.35)--(3.65,0.35)--(3.65,0.45);
 \node[below] at (2.1,0.35) {$I_{vL}$};
 \draw (3.95,0.45)--(3.95,0.35)--(7.05,0.35)--(7.05,0.45);
 \node[below] at (5.5,0.35) {$I_{vR}$};
 \draw (0.55,2.45)--(0.55,2.55)--(7.05,2.55)--(7.05,2.45);
 \node[above] at (3.8,2.55) {$I=I_v=\{k+1,\ldots,j-1\}$};
 \begin{scope}[shift={(9.6,0.2)}]
  \node[draw,circle,inner sep=1.6pt] (p) at (1.8,2.3) {$r$};
  \node[draw,circle,inner sep=1.6pt] (v) at (1.8,1.35) {$v$};
  \node[draw,circle,inner sep=1.2pt] (vl) at (0.7,0.35) {$vL$};
  \node[draw,circle,inner sep=1.2pt] (vr) at (2.9,0.35) {$vR$};
  \draw (p)--(v);
  \draw (v)--(vl);
  \draw (v)--(vr);
  \node[right=4pt] at (p) {parent of $I_v$, equal to $k$ or $j$};
  \node[right=4pt] at (v) {first step removed from $I_v$};
 \end{scope}
\end{tikzpicture}
\caption{Left: an interval $I=I_v$ between consecutive cuts $k$ and $j$
at threshold $e^y$, shown in the case $k>0$ with both child steps
present. Inserting $v$ as a cut in the removal order divides $I_v$ into
$I_{vL}$ and $I_{vR}$. Right: the corresponding fragment of the binary
removal tree showing the parent $r$, the step $v$ and the children $vL$
and $vR$.}
\label{fig:tree-split}
\end{figure}

Let $0\le y<\log h_n$. Call the steps $j$ with $h_j>e^y$ the
\emph{cuts} at threshold $e^y$. Since $h_n>e^y$, the step $n$ is a cut.
For a cut $j$ let $k$ denote the largest cut below $j$, with $k=0$ if
$j$ is the first cut. Let $I=\{k+1,\ldots,j-1\}$. Thus the schedule is
cut after every step with stepsize above $e^y$ and $I$ is the interval
between the consecutive cuts $k$ and $j$, as shown in
\cref{fig:tree-split}. In every formula that contains a cut $j$, the
letters $k$ and $I$ denote the largest cut below $j$ and the set
$\{k+1,\ldots,j-1\}$ unless $k$ is defined otherwise.

The positive steps of $I$ have stepsizes at most $e^y$ and are
removed after $j$ and also after $k$ when $k>0$. Let $k>0$ and $h_k\ge h_j$. Then $k$
is removed before $j$ and $I=I_{jL}$. Let $k>0$ and $h_k<h_j$. Then $j$
is removed before $k$ and $I=I_{kR}$. For $k=0$ we have $I=I_{jL}$.
The parent of $I$ denotes $k$ for $k>0$ and $h_k<h_j$ and denotes $j$
for $k=0$ or $h_k\ge h_j$. For $I=\varnothing$ the parent of $I$ is determined by the
cut $j$ and not by the set $I$. Hence the parent of $I$ is a cut. If $I$ contains a
positive step, then $I=I_v$ for the first step $v$ removed from $I$.

Let $I$ be an interval with $h_i<\kappa$ for every $i\in I$. On
$x^2/(2\kappa)$ the step $i$ multiplies $x$ by $1-h_i/\kappa>0$. A
negative stepsize increases $\lvert x\rvert$. Define
\begin{equation}\label{sc:eq:mass-def}
 m_I=\max_{J}\left[-\kappa\sum_{i\in J}\log(1-h_i/\kappa)\right],
\end{equation}
where the maximum is over the prefixes $J$ of $I$. Thus $e^{-m_I/\kappa}$
is the smallest product of the factors $1-h_i/\kappa$ over a prefix of
$I$. The empty prefix
gives $m_I\ge0$ and $m_\varnothing=0$. For a positive step $v$ with
$h_v<\kappa$ write $m_v,m_{vL},m_{vR}$ for $m_I$ with
$I=I_v,I_{vL},I_{vR}$. For $0<h_i\le\kappa/2$,
\begin{equation}\label{sc:eq:phi-defect}
 0\le -\kappa\log(1-h_i/\kappa)-h_i
  =\int_0^{h_i}\frac{\tau}{\kappa-\tau}\,d\tau\le h_i^2/\kappa .
\end{equation}
Let $h_i\le A\le\kappa/2$ for every $i\in I$. In
\eqref{sc:eq:mass-def} the terms with $h_i\le0$ are nonpositive. By
\eqref{sc:eq:phi-defect} the terms with $h_i>0$ are at most
$(1+A/\kappa)h_i$. Hence
\begin{equation}\label{sc:eq:mass-product}
 m_I\le (1+A/\kappa)\sum_{i\in I}h_i^+ .
\end{equation}

For $0\le y<\log h_n$ define
\begin{equation}\label{sc:eq:threshold-set}
 \mathcal P(y)=\frac{p}{2-\sqrt2}\sum_{h_j>e^y}m_I^{1/p},
\end{equation}
where the sum is over the cuts $j$ at threshold $e^y$. By
\eqref{sc:eq:mass-product} with $A=e^y$,
$m_I\le 2e^y\#\{i\in I:h_i>0\}$. Since $1/p<1$, the $1/p$ power of a
nonnegative integer is at most the integer. Summing over the cuts
gives
\begin{equation}\label{sc:eq:power-small-threshold}
 \mathcal P(y)\le C e^{y/p}\#\{i:h_i>0\}.
\end{equation}
\Cref{sc:lem:initial-power} gives $\mathcal P(y)\ge(p-o(1))n\ge n$ for
$y<\log h_n$ with $y/\log h_n\to1$. It suffices to prove
$\mathcal P(y)\ge n/2$ at some $y$ that stays bounded along the
sequence. Then \eqref{sc:eq:power-small-threshold} gives
$\#\{i:h_i>0\}\ge cn$, contrary to \eqref{sc:eq:standing-sparse}.

\subsection{Changes at a split}\label{sec:split}\label{sec:one-split}

Let $v$ be a positive step with $h_v<\kappa$. A prefix of $I_v$ is a
prefix of $I_{vL}$ or the union of $I_{vL}\cup\{v\}$ and a prefix of
$I_{vR}$. Thus \eqref{sc:eq:mass-def} gives
\begin{equation}\label{sc:eq:split-m}
 m_v=\max\left\{m_{vL},\;
 m_{vR}-\kappa\sum_{i\in I_{vL}\cup\{v\}}\log(1-h_i/\kappa)\right\}.
\end{equation}
Since $-\kappa\sum_{i\in I_{vL}}\log(1-h_i/\kappa)\le m_{vL}$,
\begin{equation}\label{sc:eq:split-size}
 m_{vL}\le m_v\le m_{vL}+m_{vR}-\kappa\log(1-h_v/\kappa).
\end{equation}

Let $j$ be a cut with $I=I_v$ and let $P>0$. \Cref{sc:lem:ordinary}
bounds $\norm{x_j}/\norm{x_k}$ from below in terms of $h_j/(m_v+P)$.
Since $m_v\ge0$, we have $m_v+P>0$. Insert $v$ as an additional cut
with the same $P$. Then the intervals before $v$ and before $j$ are
$I_{vL}$ and $I_{vR}$. The product $[h_v/(m_{vL}+P)][h_j/(m_{vR}+P)]$
replaces $h_j/(m_v+P)$ and the sum of $\log[h_j/(m_I+P)]$ over the
cuts changes by
\begin{equation}\label{sc:eq:Lambda-Delta}
 \log\frac{h_v(m_v+P)}{(m_{vL}+P)(m_{vR}+P)}.
\end{equation}
In the sum in \eqref{sc:eq:threshold-set},
$m_{vL}^{1/p}+m_{vR}^{1/p}$ replaces $m_v^{1/p}$. Define
\begin{equation}\label{sc:eq:Delta}
 \Delta_v=\frac p{2-\sqrt2}\left(m_v^{1/p}-m_{vL}^{1/p}-m_{vR}^{1/p}\right).
\end{equation}
Let $0\le s<t<\log h_n$. The cuts at threshold $e^s$ that are not cuts
at threshold $e^t$ are the steps $v$ with $e^s<h_v\le e^t$. Take the
steps $v$ with $e^s<h_v\le e^t$ as cuts one at a time in the order of
removal. Each $v$ replaces $m_v^{1/p}$ by $m_{vL}^{1/p}+m_{vR}^{1/p}$
in the sum in \eqref{sc:eq:threshold-set}. Hence
\begin{equation}\label{sc:eq:band-sums}
 \mathcal P(t)-\mathcal P(s)=\sum_{s<\log h_v\le t}\Delta_v .
\end{equation}

For $P=0$ the logarithm \eqref{sc:eq:Lambda-Delta} and $\Delta_v$ are
zero when
\begin{equation}\label{sc:eq:silver-proportions}
 m_{vL}=m_{vR}=\rho h_v,\qquad m_v=\rho^2h_v .
\end{equation}
We call \eqref{sc:eq:silver-proportions} the \emph{silver proportions}.
The constant $p/(2-\sqrt2)$ in \eqref{sc:eq:threshold-set} and
\eqref{sc:eq:Delta} is chosen such that for $P=0$ and
$m_v=m_{vL}+m_{vR}+h_v$ the difference of \eqref{sc:eq:Lambda-Delta}
and $\Delta_v/h_v^{1/p}$ is stationary in $(m_{vL},m_{vR})$ at
\eqref{sc:eq:silver-proportions}. \Cref{lem:max-plus-excess} bounds
\eqref{sc:eq:Lambda-Delta} minus $\Delta_v/h_v^{1/p}$ from below.

\begin{lemma}\label{lem:max-plus-excess}
Let $h_v>0$ and $P\ge0$.
\begin{enumerate}
\item\label{it:split-real} Let $m_v,m_{vL},m_{vR}\ge0$ be real numbers
with $m_{vL}+P>0$, $m_{vR}+P>0$ and
\begin{equation}\label{sc:eq:split-hypothesis}
 m_{vL}\le m_v\le m_{vL}+m_{vR}+h_v .
\end{equation}
With $\Delta_v$ defined by \eqref{sc:eq:Delta},
\begin{multline}\label{sc:eq:split-stability}
 \log\frac{h_v(m_v+P)}{(m_{vL}+P)(m_{vR}+P)}
   -\frac{\Delta_v}{h_v^{1/p}}\\
 \ge c\min\!\left\{1,\sum_{w=vL,vR}\left(\frac{m_w}{h_v}-\rho\right)^2\right\}
 -C\left(\frac P{h_v}\right)^{1/p}.
\end{multline}
\item\label{it:split-mass} Let $v$ be a positive step with $h_v\le h_n$
and let $m_v,m_{vL},m_{vR}$ be the numbers of \eqref{sc:eq:mass-def}
with $m_{vL}+P>0$ and $m_{vR}+P>0$. Then
\begin{multline}\label{sc:eq:split-stability-mass}
 \log\frac{h_v(m_v+P)}{(m_{vL}+P)(m_{vR}+P)}
   -\frac{\Delta_v}{h_v^{1/p}}\\
 \ge c\min\!\left\{1,\sum_{w=vL,vR}\left(\frac{m_w}{h_v}-\rho\right)^2\right\}
 -C\left(\frac P{h_v}\right)^{1/p}-\frac{Ch_v}\kappa .
\end{multline}
\end{enumerate}
\end{lemma}
\begin{proof}
In \ref{it:split-real} both sides of \eqref{sc:eq:split-stability} and
condition \eqref{sc:eq:split-hypothesis} are unchanged when
$m_v,m_{vL},m_{vR},P$ and $h_v$ are multiplied by one positive number.
We may therefore take $h_v=1$ in \ref{it:split-real}. Define
\[
 \Phi(m_v,m_{vL},m_{vR},P)
 =\log\frac{m_v+P}{(m_{vL}+P)(m_{vR}+P)}-\Delta_v .
\]

First let $P=0$. As a function of $m_v^{1/p}$,
$\Phi(m_v,m_{vL},m_{vR},0)$ is concave. On the range
\eqref{sc:eq:split-hypothesis}, $\Phi(m_v,m_{vL},m_{vR},0)$ is
therefore smallest at $m_v=m_{vL}$ or at $m_v=m_{vL}+m_{vR}+1$. At
$m_v=m_{vL}$ we have
$\Phi(m_{vL},m_{vL},m_{vR},0)=p(2-\sqrt2)^{-1}m_{vR}^{1/p}-\log m_{vR}\ge p$
by $\log \tau\le p(\tau^{1/p}-1)$ with $\tau=m_{vR}$.

Let $m_v=m_{vL}+m_{vR}+1$ and define
\[
 \varphi(\tau)=\frac{p}{2-\sqrt2}\tau^{1/p}-\log \tau .
\]
Then $\Phi(m_v,m_{vL},m_{vR},0)=\varphi(m_{vL})+\varphi(m_{vR})-\varphi(m_v)$.
The derivative of $\varphi$ is
\[
 \varphi'(\tau)=\frac{\tau^{1/p}/(2-\sqrt2)-1}{\tau}.
\]
The second derivative $\varphi''(\tau)$ has the sign of
$1-(1-1/p)\tau^{1/p}/(2-\sqrt2)$. Thus $\varphi'(\tau)$ increases and then
decreases as $\tau$ increases. Also $\varphi'(\rho)=\varphi'(\rho^2)=1$
and the unique maximizer of $\varphi'$ lies in $(\rho,\rho^2)$.

Let $w\in\{vL,vR\}$ have $m_w\ge\rho^2$. Since $m_v>m_w\ge\rho^2$ and
$\varphi'$ decreases on $[\rho^2,\infty)$, we have
$\varphi'(m_v)<\varphi'(m_w)$ and decreasing $m_w$ strictly decreases
$\Phi(m_v,m_{vL},m_{vR},0)$. Let $\max\{m_{vL},m_{vR}\}\ge\rho^2$. Replacing $m_{vL}$ and $m_{vR}$
by $\min\{m_{vL},\rho^2\}$ and $\min\{m_{vR},\rho^2\}$ therefore does
not increase $\Phi(m_{vL}+m_{vR}+1,m_{vL},m_{vR},0)$. Since
$\rho^2-\rho>1$, we have $\min\{1,\sum_{w=vL,vR}(m_w-\rho)^2\}=1$
before and after the replacement.
Hence it suffices to prove \eqref{sc:eq:split-stability} for
$0<m_{vL},m_{vR}\le\rho^2$. We have
$\Phi(m_{vL}+m_{vR}+1,m_{vL},m_{vR},0)\to\infty$ as $m_{vL}\to0$ or
$m_{vR}\to0$. Because decreasing $m_w$ from $\rho^2$ strictly decreases
$\Phi(m_{vL}+m_{vR}+1,m_{vL},m_{vR},0)$, any minimum on $(0,\rho^2]^2$
lies in $(0,\rho^2)^2$. At a minimum in $(0,\rho^2)^2$, the derivatives
satisfy
$\varphi'(m_{vL})=\varphi'(m_{vR})=\varphi'(m_v)$. Since $m_v>m_{vL}$
and $m_v>m_{vR}$, both $m_{vL}$ and $m_{vR}$ lie where $\varphi'$
increases. Hence $m_{vL}=m_{vR}$. For $m_{vL}=m_{vR}=\tau$ and
$m_v=2\tau+1$, $\varphi'(\tau)$ increases with $\tau$ and $\varphi'(2\tau+1)$
decreases with $\tau$ wherever $\varphi'(\tau)=\varphi'(2\tau+1)$ can hold.
Thus $\Phi(m_{vL}+m_{vR}+1,m_{vL},m_{vR},0)$ is smallest only at
$m_{vL}=m_{vR}=\rho$ and $m_v=2\rho+1=\rho^2$, where
$\Phi(\rho^2,\rho,\rho,0)=0$.

Near $\tau=\rho$ we have $\varphi''(\tau)>0$. Near $\tau=\rho^2$ we
have $\varphi''(\tau)<0$. Thus the second-order term
$\tfrac12\varphi''(\rho)[(m_{vL}-\rho)^2+(m_{vR}-\rho)^2]-\tfrac12\varphi''(\rho^2)(m_{vL}+m_{vR}-2\rho)^2$
is positive definite and Taylor's theorem gives
$\Phi(m_{vL}+m_{vR}+1,m_{vL},m_{vR},0)\ge c\sum_{w=vL,vR}(m_w-\rho)^2$
for $(m_{vL},m_{vR})$ near $(\rho,\rho)$. For all other
$m_{vL},m_{vR}\le\rho^2$, compactness bounds
$\Phi(m_{vL}+m_{vR}+1,m_{vL},m_{vR},0)$ from below by a positive
constant. Thus \eqref{sc:eq:split-stability} holds for $P=0$.

Let $P>0$. Since $m_{vL}\le m_v\le m_{vL}+m_{vR}+1$, we have
$m_{vL}+P\le m_v+P\le(m_{vL}+P)+(m_{vR}+P)+1$. Thus $m_v+P$,
$m_{vL}+P$, $m_{vR}+P$ satisfy \eqref{sc:eq:split-hypothesis} with
$h_v=1$. Apply
\eqref{sc:eq:split-stability} with $P=0$ to $m_v+P,m_{vL}+P,m_{vR}+P$.
Replacing $m_v+P,m_{vL}+P,m_{vR}+P$ by $m_v,m_{vL},m_{vR}$ changes each
of the three $1/p$ powers in $\Delta_v$ by at most $P^{1/p}$ and
changes $\min\{1,\sum_{w=vL,vR}(m_w-\rho)^2\}$ by at most
$C\min\{P,1\}\le CP^{1/p}$. Hence \eqref{sc:eq:split-stability} holds.

Let $m_v,m_{vL},m_{vR}$ be the numbers of \eqref{sc:eq:mass-def} with
$h_v\le h_n$ and set $\tilde m_w=m_w/h_v$ for $w\in\{v,vL,vR\}$ and
$\tilde P=P/h_v$. If $\tilde m_v\le\tilde m_{vL}+\tilde m_{vR}+1$, then
$\tilde m_v,\tilde m_{vL},\tilde m_{vR}$ satisfy
\eqref{sc:eq:split-hypothesis} with $h_v=1$ and \ref{it:split-real}
applied with
$(h_v,P,m_v,m_{vL},m_{vR})=(1,\tilde P,\tilde m_v,\tilde m_{vL},\tilde m_{vR})$
gives \eqref{sc:eq:split-stability} for $h_v$, $P$, $m_v$, $m_{vL}$,
$m_{vR}$. Let
$\tilde m_v>\tilde m_{vL}+\tilde m_{vR}+1$. By \eqref{sc:eq:split-size}
and \eqref{sc:eq:phi-defect}, $\tilde m_v$ exceeds
$\tilde m_{vL}+\tilde m_{vR}+1$ by at most $h_v/\kappa$. For
$\tilde m_v\ge1$,
$\lvert\partial\Phi(\tilde m_v,\tilde m_{vL},\tilde m_{vR},\tilde P)/\partial\tilde m_v\rvert
=\lvert(\tilde m_v+\tilde P)^{-1}-(2-\sqrt2)^{-1}\tilde m_v^{1/p-1}\rvert\le C$.
Decreasing $\tilde m_v$ to $\tilde m_{vL}+\tilde m_{vR}+1$ therefore
changes $\Phi(\tilde m_v,\tilde m_{vL},\tilde m_{vR},\tilde P)$ by at
most $Ch_v/\kappa$. Also \eqref{sc:eq:split-hypothesis} with $h_v=1$ holds for
$(m_v,m_{vL},m_{vR})=(\tilde m_{vL}+\tilde m_{vR}+1,\tilde m_{vL},\tilde m_{vR})$.
Hence \ref{it:split-real} gives \eqref{sc:eq:split-stability-mass}.
\end{proof}

\section{Backward construction of a gradient descent trajectory}\label{sec:backward}
In this section we construct a trajectory of a $1$-smooth,
$1/\kappa$-strongly convex function with minimizer zero on one interval
$I=\{k+1,\ldots,j-1\}$ followed by the step $j$. The points are chosen
in the order $x_{j-1},x_{j-2},\ldots,x_k$. \Cref{sc:lem:assembly}
reduces the construction of one function for consecutive intervals to
\eqref{eq:tail-planes} for adjacent intervals. \Cref{sc:lem:ordinary}
gives $\norm{x_j}/\norm{x_k}\ge h_j/(m_I+P)$ up to the factor
$\exp(-C/\log A-m_I/\kappa)$. The proof of \cref{sc:lem:ordinary} uses
a two-dimensional solution from \cref{sc:lem:arc} and a one-dimensional solution
from \cref{sc:lem:ramp}. For a vector $x\ne0$ write
$\widehat x=x/\norm x$.
\subsection{Prescribing gradients at earlier points}\label{sc:sec:trajectories}

\subsubsection{One function for consecutive intervals}\label{sc:sec:objectives}

For integers $k<j$ write $(k,j]=\{k+1,\ldots,j\}$. All vectors lie in
one Euclidean space. Let $x_k,\ldots,x_j$ be points. For $k\le i<j$
let $z_i,g_i$ be vectors with $x_i=z_i+g_i$ and let $q_i$ be a real
number. We use the conditions
\begin{align}
 x_{i+1}&=x_i-h_{i+1}\bigl(x_i/\kappa+(1-1/\kappa)g_i\bigr),
 \label{sc:eq:tail-update}\\
 q_i+\ip{g_i}{z_u-z_i}&\le q_u .
 \label{eq:tail-planes}
\end{align}
In \eqref{sc:eq:tail-update} the gradient at $x_i$ of the function
$f$ of \cref{sc:lem:assembly} is $x_i/\kappa+(1-1/\kappa)g_i$. Condition \eqref{eq:tail-planes} states
that the affine function $z\mapsto q_i+\ip{g_i}{z-z_i}$ is at most
$q_u$ at $z_u$. Hence there is a convex function $F$ with
$F(z_i)=q_i$ and $g_i\in\partial F(z_i)$. The Moreau envelope
$M_F(x)=\inf_z\{F(z)+\norm{x-z}^2/2\}$ is convex and $1$-smooth. Since
$g_i\in\partial F(z_i)$, the infimum at $x_i=z_i+g_i$ is attained at
$z_i$ and $\nabla M_F(x_i)=x_i-z_i=g_i$. The function of
\cref{sc:lem:assembly} is
$f(x)=(1-1/\kappa)M_F(x)+\norm x^2/(2\kappa)$.

We call points $x_k,\ldots,x_j$ and triples $(z_i,g_i,q_i)$ for
$k\le i<j$ with $x_i=z_i+g_i$ a \emph{solution} on $(k,j]$ if
\eqref{sc:eq:tail-update} holds for $k\le i<j$ and
\eqref{eq:tail-planes} holds for every pair of triples among $(0,0,0)$
and $(z_i,g_i,q_i)$ for $k\le i<j$. With $(0,0,0)$ as one of the two
triples, \eqref{eq:tail-planes} gives
\[
 0\le q_i\le\ip{g_i}{z_i}.
\]
Let $\mathcal S$ be a solution on $(k,j]$ with points $x_i$ and triples
$(z_i,g_i,q_i)$. We call $(k,j]$ the \emph{interval} of $\mathcal S$,
$x_k$ the \emph{first point} of $\mathcal S$, $x_j$ the \emph{last
point} of $\mathcal S$ and $z_i,g_i$ for $k\le i<j$ the \emph{vectors}
of $\mathcal S$. For $r\in\{1,2,4\}$ we call $\mathcal S$
$r$-dimensional if the points and the vectors of $\mathcal S$ lie in
$\R^r$. Applying one rotation to every $x_i,z_i,g_i$ gives a
solution on $(k,j]$. For $\sigma>0$ the map
\begin{equation}\label{eq:scaling}
 (x_i,z_i,g_i,q_i)\mapsto(\sigma x_i,\sigma z_i,\sigma g_i,\sigma^2q_i)
\end{equation}
gives a solution on $(k,j]$. We call a solution obtained from
$\mathcal S$ by one rotation and by \eqref{eq:scaling} with one
$\sigma>0$ an \emph{image} of $\mathcal S$ under a rotation and
\eqref{eq:scaling}. Let $\tilde{\mathcal S}$ be a second
solution on $(k,j]$ with points $\tilde x_i$ and triples
$(\tilde z_i,\tilde g_i,\tilde q_i)$. If the vectors of
$\tilde{\mathcal S}$ lie in a subspace orthogonal to the vectors of
$\mathcal S$, then $\ip{g_i}{\tilde z_u}=\ip{\tilde g_i}{z_u}=0$. Hence the points
$x_i+\tilde x_i$ and triples
\begin{equation}\label{sc:eq:piece-sum}
 (z_i+\tilde z_i,\;g_i+\tilde g_i,\;q_i+\tilde q_i)
\end{equation}
form a solution on $(k,j]$.

\begin{lemma}\label{sc:lem:assembly}
Let finitely many solutions on consecutive intervals be given in one
Euclidean space. Assume that:
\begin{enumerate}
\item\label{it:asm-shared} the last point of each solution is the first
point of the next solution and the points shared by two solutions are
nonzero and pairwise orthogonal.
\item\label{it:asm-span} for each solution $\mathcal S$ the vectors
of $\mathcal S$ lie in the sum of a subspace
$U_{\mathcal S}$ and the span of the shared points of $\mathcal S$.
\item\label{it:asm-perp} the subspaces $U_{\mathcal S}$ of distinct
solutions are orthogonal to each other and to every shared point.
\end{enumerate}
If \eqref{eq:tail-planes} holds for every two triples of two solutions
with a shared point, then there is a $1$-smooth, $1/\kappa$-strongly
convex function with minimizer zero whose trajectory consists of the
points of the given solutions.
\end{lemma}
\begin{proof}
Let the triples at $i$ and at $j$ belong to two solutions without a
shared point. By \ref{it:asm-shared}--\ref{it:asm-perp},
$\ip{g_i}{z_j}=0$ and hence
\[
 q_i+\ip{g_i}{z_j-z_i}\le0\le q_j .
\]
Thus \eqref{eq:tail-planes} holds for every two triples of the given
solutions.

Let $F$ be the maximum of the zero function and the affine functions
$z\mapsto q_i+\ip{g_i}{z-z_i}$ of all triples. By \eqref{eq:tail-planes}, $g_i$ is
a subgradient of $F$ at $z_i$ and zero is a subgradient of $F$ at the
origin. The Moreau envelope $M_F(x)=\inf_z\{F(z)+\norm{x-z}^2/2\}$ is
convex and $1$-smooth \cite{THG17}. At $x_i=z_i+g_i$ the convex
function $z\mapsto F(z)+\norm{x_i-z}^2/2$ has zero subgradient at $z_i$. Hence
$\nabla M_F(x_i)=x_i-z_i=g_i$. Since zero is a subgradient of $F$ at
the origin, $\nabla M_F(0)=0$. Define
$f(x)=(1-1/\kappa)M_F(x)+\norm x^2/(2\kappa)$. Then $f$ is $1$-smooth
and $1/\kappa$-strongly convex with minimizer zero and gradient
$x_i/\kappa+(1-1/\kappa)g_i$ at $x_i$. Hence \eqref{sc:eq:tail-update}
is the gradient step with stepsize $h_{i+1}$.
\end{proof}

Let \cref{sc:lem:assembly}\ref{it:asm-shared}--\ref{it:asm-perp} hold
and let $x$ be the point shared by two solutions $\mathcal S$ and
$\mathcal S'$. Then every triple $(z,g,q)$ of $\mathcal S$ and every
triple $(z',g',q')$ of $\mathcal S'$ satisfy
\begin{equation}\label{eq:adjacent-support}
 \ip g{z'}=\ip g{\widehat x}\ip{\widehat x}{z'},\qquad
 \ip{g'}z=\ip{g'}{\widehat x}\ip{\widehat x}z .
\end{equation}
Let $\mathcal S_1,\ldots,\mathcal S_r$ be solutions on consecutive
intervals. Suppose that the last point of $\mathcal S_1$ is nonzero,
that the first point of $\mathcal S_r$ is nonzero and that the first
and last points of $\mathcal S_i$ are nonzero and orthogonal for
$1<i<r$. Then rotations and \eqref{eq:scaling} in a Euclidean space of
sufficiently large dimension give images of
$\mathcal S_1,\ldots,\mathcal S_r$ that satisfy
\cref{sc:lem:assembly}\ref{it:asm-shared}--\ref{it:asm-perp}.

\subsubsection{A two-dimensional trajectory}\label{sc:sec:planar}

Let $K$ be a compact convex set containing zero and let $\Pi_K$ be the
Euclidean projection onto $K$. At a point $x$ take the triple $(z,g,q)$ with
\[
 g=\Pi_K(x),\qquad z=x-g,\qquad q=\ip gz .
\]
The projection satisfies
\begin{equation}\label{sc:eq:planar-subgradient}
 \ip\xi z\le\ip gz\qquad(\xi\in K).
\end{equation}
Let $(z_i,g_i,q_i)$ and $(z_u,g_u,q_u)$ be the triples at two points.
By \eqref{sc:eq:planar-subgradient} for $z_u$ with $\xi=g_i$,
\[
 q_i+\ip{g_i}{z_u-z_i}=\ip{g_i}{z_u}\le\ip{g_u}{z_u}=q_u .
\]
Since $0\in K$, \eqref{sc:eq:planar-subgradient} gives $q_u\ge0$. Hence
\eqref{eq:tail-planes} holds for every pair of triples among $(0,0,0)$
and the triples $(z,g,q)$ at any finite set of points $x$.

A point of the plane is written $x=(x_1,x_2)$. For $g=\Pi_K(x)$ and
$z=x-g$ the coordinates are $g_1,g_2$ and $z_1,z_2$. The coordinates of
a point $x_i$ are $(x_i)_1$ and $(x_i)_2$. \Cref{sc:lem:arc} starts from $x_{j-1}=(1/A,1)\in K$ and chooses
$x_{j-2},x_{j-3},\ldots,x_k$ one after another. In the proof of
\cref{sc:lem:ordinary} the step $j$ gives $x_j=(1-h_j)(1/A,1)$ in the
first two coordinates.

\begin{lemma}\label{sc:lem:arc}
Let $A<\kappa$ be sufficiently large. Let $0\le k<j\le n$ and
$h_i<\kappa$ for $k<i<j$. There are a compact convex set
$K\subset[0,\infty)^2$ and points $x_k,\ldots,x_{j-1}\in(0,\infty)^2$
with
\[
 x_{j-1}=(1/A,1)\in K,\qquad x_i=x_{i-1}-h_i
 [x_{i-1}/\kappa+(1-1/\kappa)\Pi_K(x_{i-1})]\quad(k<i<j).
\]
\begin{enumerate}
\item\label{sc:it:arc-size} $K$ contains the quarter unit disk
$\{\xi\in[0,\infty)^2:\norm\xi\le1\}$ and $\norm\xi\le1+C/\log A$ for
every $\xi\in K$.
\item\label{sc:it:arc-signs} For every $x\in[0,\infty)^2$,
$g=\Pi_K(x)$, $z=x-g$ and $q=\ip gz$ satisfy $q\ge\norm z$ and
$g,z\in[0,\infty)^2$.
\item\label{sc:it:arc-bounds} For every $k\le i<j$ the interval
$J=\{i+1,\ldots,j-1\}$ satisfies
\begin{align}
 \norm{x_i}&\le e^{m_J/\kappa}\big[(1+C/\log A)m_J+\sqrt A\big],
 \label{sc:eq:arc-bounds}\\
 (x_i)_2&\le e^{m_J/\kappa}\sqrt A .\label{sc:eq:arc-height}
\end{align}
\end{enumerate}
\end{lemma}

For nonnegative steps, Ye and Liu \cite{ye2026almost} construct
backward trajectories by projection onto the convex hull of a circular
arc and the origin. For real steps, He and Cheng \cite{he2026optimal}
obtain a predecessor by a fixed-point argument on a closed convex set.
We use a concave curve through $x_{j-1}$ to bound the second
coordinate.

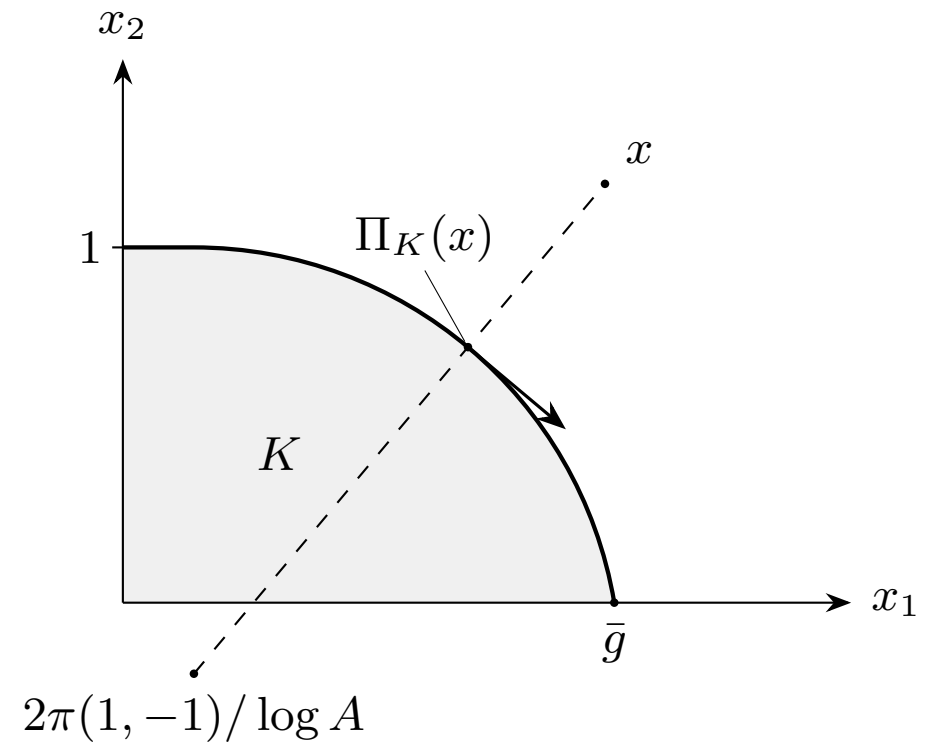
\begin{figure}[htbp]
\centering
\begin{tikzpicture}[x=2.35cm,y=2.35cm,>=Stealth,font=\footnotesize]
 \pgfmathsetmacro{\lastangle}{asin(1/6)}
 \coordinate (c) at (.2,-.2);
 \coordinate (v) at ({.2+1.2*cos(50)},{-.2+1.2*sin(50)});
 \coordinate (X) at ({.2+1.8*cos(50)},{-.2+1.8*sin(50)});
 \fill[black!6] (0,0)--(0,1)--(.2,1)
       arc[start angle=90,end angle=\lastangle,radius=1.2]--cycle;
 \draw[->] (0,0)--(2.05,0) node[right] {$x_1$};
 \draw[->] (0,0)--(0,1.53) node[above] {$x_2$};
 \draw (-.03,1)--(0,1) node[left] {$1$};
 \draw[thick] (0,1)--(.2,1)
       arc[start angle=90,end angle=\lastangle,radius=1.2];
 \node at (.44,.42) {$K$};
 \fill ({.2+sqrt(1.4)},0) circle (.8pt) node[below] {$\bar g$};
 \draw[dashed] (c)--(X);
 \draw[->,semithick] (v)--+({.36*sin(50)},{-.36*cos(50)});
 \fill (c) circle (.8pt) node[below] {$2\pi(1,-1)/\log A$};
 \fill (v) circle (.8pt);
 \node[anchor=south,inner sep=1.5pt] (lab) at (.85,.935) {$\Pi_K(x)$};
 \draw[very thin,shorten >=1.2pt] (lab.south)--(v);
 \fill (X) circle (.8pt) node[above right] {$x$};
\end{tikzpicture}
\caption{The set $K$ in the proof of \cref{sc:lem:arc}, drawn with
$2\pi/\log A=1/5$. The arc of the boundary of $K$ has center
$2\pi(1,-1)/\log A$ and ends at $\bar g$. The arrow is the clockwise
unit tangent of the arc at $\Pi_K(x)$.}
\label{sc:fig:arc}
\end{figure}

\begin{proof}
Let $K\subset[0,\infty)^2$ be the set of $x=(x_1,x_2)\in[0,\infty)^2$
satisfying
\begin{equation}\label{sc:eq:arc-set}
 \left[\left(x_1-\frac{2\pi}{\log A}\right)^+\right]^2
 +\left(x_2+\frac{2\pi}{\log A}\right)^2
 \le\left(1+\frac{2\pi}{\log A}\right)^2 .
\end{equation}
The upper boundary of $K$ consists of a horizontal segment from $(0,1)$
to $(2\pi/\log A,1)$ and of an arc from $(2\pi/\log A,1)$ to a point of
the horizontal axis, as shown in \cref{sc:fig:arc}. The arc lies on the
circle with center $2\pi(1,-1)/\log A$ and radius $1+2\pi/\log A$.

Every point $\xi=(\xi_1,\xi_2)\in[0,\infty)^2$ with $\norm\xi\le1$
satisfies \eqref{sc:eq:arc-set} since $\xi_1^2+\xi_2^2\le1$ and
$\xi_2\le1$. Let $x\in K$. By \eqref{sc:eq:arc-set} the vector
$((x_1-2\pi/\log A)^+,x_2)$ has norm at most $1+2\pi/\log A$. Hence
\[
 \norm x\le1+4\pi/\log A .
\]
Thus \ref{sc:it:arc-size} holds. The function
$x\mapsto[(x_1-2\pi/\log A)^+]^2+(x_2+2\pi/\log A)^2$ is convex and
nondecreasing in each coordinate on $[0,\infty)^2$. Therefore $K$ is
convex and $0\le g\le x$ coordinatewise for every $x\in[0,\infty)^2$.
For $z\ne0$ we have $z/\norm z\in K$. By
\eqref{sc:eq:planar-subgradient} with $\xi=z/\norm z$,
\[
 q=\ip gz\ge\norm z .
\]
Thus \ref{sc:it:arc-signs} holds.

We construct $x_{i-1}$ from $x_i$. Let $k<i<j$ and let
$x_i\in(0,\infty)^2$. Let $g$ maximize
\[
 \ip{x_i}g+(h_i-1)\norm g^2/2
\]
over $g\in K$. Define
\[
 z=\frac{x_i+(h_i-1)g}{1-h_i/\kappa},\qquad x_{i-1}=g+z .
\]
At the maximizer $g$ the derivative in the direction $\xi-g$ is
nonpositive for every $\xi\in K$. Hence
\[
 \ip{\xi-g}z\le0\qquad(\xi\in K)
\]
and therefore $\Pi_K(x_{i-1})=g$. By the definition of $z$,
\[
 x_{i-1}-h_i\bigl[x_{i-1}/\kappa+(1-1/\kappa)g\bigr]
 =(1-h_i/\kappa)z+(1-h_i)g=x_i .
\]
Let $r\in\{1,2\}$. If $g_r>0$, a decrease of $g_r$ gives a point of $K$. Hence $z_r\ge0$. If $g_r=0$, then
$z_r=(x_i)_r/(1-h_i/\kappa)>0$. Thus $x_{i-1}\in(0,\infty)^2$. Starting
from $x_{j-1}=(1/A,1)$ we obtain $x_{j-2},\ldots,x_k\in(0,\infty)^2$.

To bound $x_i$ we use the differential equation
\begin{equation}\label{sc:eq:arc-ode}
 \frac{dx_2}{dx_1}=\frac{x_2+(\kappa-1)g_2}{x_1+(\kappa-1)g_1},
 \qquad g=\Pi_K(x),
\end{equation}
for a positive function $x_2$ of $x_1>0$. For $x\in(0,\infty)^2$ the
quotient $[x_2+(\kappa-1)g_2]/[x_1+(\kappa-1)g_1]$ in
\eqref{sc:eq:arc-ode} is locally Lipschitz and lies between
$x_2/(\kappa x_1)$ and $\kappa x_2/x_1$. Hence every positive function
$x_2$ that satisfies \eqref{sc:eq:arc-ode} on an interval of $x_1$
extends to a positive function that satisfies \eqref{sc:eq:arc-ode} for
all $x_1>0$. In this proof a graph is the graph of a positive function
$x_2$ of $x_1>0$ that satisfies \eqref{sc:eq:arc-ode}. The graphs are integral curves of the
vector field $x\mapsto x/\kappa+(1-1/\kappa)\Pi_K(x)$. The graphs are
increasing and pairwise disjoint. Every point $x\in(0,\infty)^2$ lies
on exactly one graph, the graph through $x$.

We show that every graph is concave. Let $\bar g$ denote the endpoint of
the arc on the horizontal axis. Inside $K$ the graphs are rays from the
origin. Let $x$ move along a graph outside $K$ with increasing $x_1$.
Then $g=\Pi_K(x)$ moves to the right along the horizontal segment,
moves clockwise along the arc or stays at $\bar g$. Indeed, above the
horizontal segment $g=(x_1,1)$. Let $g$ lie on the arc with $g\ne\bar g$.
Both $g-2\pi(1,-1)/\log A$ and $x-g$ are orthogonal to the arc at $g$.
Hence the component of $x/\kappa+(1-1/\kappa)g$ along the clockwise
unit tangent equals the component of the center $2\pi(1,-1)/\log A$
along the clockwise unit tangent. Since the clockwise unit tangent has
nonnegative first coordinate and nonpositive second coordinate, the
component of the center is at least $2\pi/\log A$. Let $g=\bar g$.
Because the circle center $2\pi(1,-1)/\log A$ lies below the horizontal
axis, every $\xi\in K$ satisfies $\xi_1\le\bar g_1$. Hence the normal
cone of $K$ at $\bar g$ contains $\bar g$ and $x-\bar g$.
Therefore $\Pi_K=\bar g$ on the half-line from $x$ in the direction
$x/\kappa+(1-1/\kappa)\bar g$. Thus $g_1$ is nondecreasing and $g_2$ is nonincreasing along a graph
outside $K$. Along a graph outside $K$, using \eqref{sc:eq:arc-ode}
to cancel the terms without derivatives of $g$, we obtain almost
everywhere
\[
 \frac{d}{dx_1}\,\frac{x_2+(\kappa-1)g_2}{x_1+(\kappa-1)g_1}
 =(\kappa-1)\,
 \frac{[x_1+(\kappa-1)g_1]\,\dfrac{dg_2}{dx_1}
      -[x_2+(\kappa-1)g_2]\,\dfrac{dg_1}{dx_1}}
      {[x_1+(\kappa-1)g_1]^2}\le0 .
\]
Therefore every graph is concave.

By \eqref{sc:eq:arc-ode} and the recursion for $x_i$ in
\cref{sc:lem:arc},
$x_i$ lies on the tangent line at $x_{i-1}$ of the graph through
$x_{i-1}$. By concavity $x_i$ lies on or above the graph through
$x_{i-1}$. Since the graphs are disjoint, $x_{i-1}$ lies on or below
the graph through $x_i$. Hence every $x_i$ with $k\le i<j$ lies on or
below the graph through $x_{j-1}=(1/A,1)$.

We bound the graph through $x_{j-1}$. Above the horizontal segment
$g=(x_1,1)$ and \eqref{sc:eq:arc-ode} becomes
$dx_2/dx_1=(x_2+\kappa-1)/(\kappa x_1)$ with $x_2(1/A)=1$. Hence
\[
 x_2=1+\kappa\bigl[(Ax_1)^{1/\kappa}-1\bigr]
 \qquad(1/A\le x_1\le2\pi/\log A).
\]
Since $A<\kappa$, at $x_1=2\pi/\log A$ we get $x_2=O(\log A)$. Let $R$
denote the distance from $x$ to the center $2\pi(1,-1)/\log A$ and let
$\theta$ denote the clockwise polar angle of $x$ around the center
$2\pi(1,-1)/\log A$. Let $x\notin K$, $R<\sqrt A$ and let
$g=\Pi_K(x)$ lie on the arc with $g\ne\bar g$. Then
\ref{sc:it:arc-size} and $A<\kappa$ give
\[
 \norm{x/\kappa+(1-1/\kappa)g}\le1+O(1/\log A).
\]
Also $x-g\ne0$. Along the graph through $x$, the derivative
$d\log R/d\theta$ equals the component of $x/\kappa+(1-1/\kappa)g$
along the unit vector $\widehat{x-g}$ divided by the component of
$x/\kappa+(1-1/\kappa)g$ along the clockwise unit tangent of the arc
at $g$. Hence
\[
 \frac{d\log R}{d\theta}\le\frac{1+O(1/\log A)}{2\pi/\log A}
 =\frac{\log A+O(1)}{2\pi} .
\]
Here $\theta$ increases by at most $\pi/2$ while $g$ lies on the
arc. We have $R=O(\log A)$ at $x_1=2\pi/\log A$ and
$C(\log A)A^{1/4}<\sqrt A$. Hence the graph through $x_{j-1}$ does not
reach $R=\sqrt A$ before $g=\bar g$ and satisfies
$R\le C(\log A)A^{1/4}<\sqrt A$ while $g\ne\bar g$.
Because $(0,-1)$ lies in the normal cone of $K$ at $\bar g$, the set
$\{x:\Pi_K(x)=\bar g\}$ is closed under downward displacement. Let $x\in(0,\infty)^2$ lie on or below the graph through
$x_{j-1}$ with $\Pi_K(x)\ne\bar g$. Then the point on the graph through
$x_{j-1}$ directly above $x$ with first coordinate $x_1$ also has
projection different from $\bar g$. Therefore
\[
 \norm x\le x_1+x_2\le C(\log A)A^{1/4}<\sqrt A/2 .
\]

Fix $k\le i<j$ and put $J=\{i+1,\ldots,j-1\}$. Let $u\ge i$ be the
smallest index with $\Pi_K(x_u)\ne \bar g$. Then $\norm{x_u}\le\sqrt A/2$.
For $i<r\le u$ we have $\Pi_K(x_{r-1})=\bar g$ and the recursion for
$x_r$ in \cref{sc:lem:arc} gives
\[
 x_r+(\kappa-1)\bar g
 =(1-h_r/\kappa)\bigl(x_{r-1}+(\kappa-1)\bar g\bigr).
\]
Put $E=\prod_{i<r\le u}(1-h_r/\kappa)^{-1}$. Then
$x_i=Ex_u+(E-1)(\kappa-1)\bar g$. By \eqref{sc:eq:mass-def},
$E\le e^{m_J/\kappa}$ and $e^{m_J/\kappa}\ge1$. Also $\bar g$ lies
on the horizontal axis with $\norm{\bar g}\le1+C/\log A$. Since
$x_i,x_u,\bar g\in[0,\infty)^2$, we have
$x_i\le Ex_u+(E-1)^+(\kappa-1)\bar g$ coordinatewise and dividing by
$e^{m_J/\kappa}\ge\max\{E,1\}$ gives
\[
 e^{-m_J/\kappa}(x_i)_2\le\sqrt A/2,\qquad
 e^{-m_J/\kappa}\norm{x_i}\le\sqrt A/2+
 (1-e^{-m_J/\kappa})(\kappa-1)(1+C/\log A).
\]
Now $1-e^{-m_J/\kappa}\le m_J/\kappa$ gives \eqref{sc:eq:arc-bounds}
and \eqref{sc:eq:arc-height}.
\end{proof}

\subsubsection{A scalar trajectory}\label{sc:sec:ramp}

Let $(z_i,g_i,q_i)$ be a triple on $(k,j]$ and let $(z,g,q)$ be a
triple with $0\le q\le\ip gz$. Then $q+\ip g{z_i-z}\le q_i$ and
$q_i+\ip{g_i}{z-z_i}\le q$ follow from $\ip g{z_i}\le q_i$ and
$\ip{g_i}z\le\ip{g_i}{z_i}-q_i$. \Cref{sc:lem:ramp} constructs a
one-dimensional solution on $(k,j]$ with last point zero and a lower
bound for $\min\{q_i,g_iz_i-q_i\}$. For integers $k<j$ and a number $A$
we use the condition
\begin{equation}\label{sc:eq:attach-setting}
 2<h_j,\qquad \sum_{k<i\le j}|h_i|\le A\le\kappa/2.
\end{equation}
\Cref{sc:lem:ramp}\ref{it:ramp-good} is used in the proof of
\cref{sc:lem:ordinary} and \cref{sc:lem:ramp}\ref{it:ramp-bad} in the
proof of \cref{sc:lem:bad-block}.

\begin{lemma}\label{sc:lem:ramp}
Let $0\le k<j\le n$ and $I=\{k+1,\ldots,j-1\}$.
\begin{enumerate}
\item\label{it:ramp-good} Assume \eqref{sc:eq:attach-setting}. Then there is a one-dimensional solution
on $(k,j]$ with $x_j=0$ such that for $k\le i<j$,
\begin{align}
 ch_j\le x_i&\le CA,\qquad h_j\le x_{j-1}\le 2h_j,\label{sc:eq:ramp-step}\\
 0<g_i&\le g_{j-1}=1,\label{sc:eq:ramp-slope}\\
 \min\{q_i,g_iz_i-q_i\}&\ge ch_j^2/A .\label{sc:eq:ramp-height}
\end{align}
\item\label{it:ramp-bad} Assume $2<h_j\le A\le\kappa/2$ and
$h_u\ne\kappa$ for every $u\in I$. Then there is a one-dimensional
solution on $(k,j]$ with $x_j=0$ such that for $k\le i<j$,
\begin{align}
 \min\{q_i,g_iz_i-q_i\}&\ge c/A,\label{sc:eq:ramp-value}\\
 \log\frac{|x_i|+|z_i|+A|g_i|}{A/h_j}
 &\le C+2\sum_{u\in I}\bigl[\log(1+|h_u|/A)+(\log|1-h_u/\kappa|)^-\bigr].
 \label{sc:eq:ramp-cost}
\end{align}
\end{enumerate}
\end{lemma}
\begin{proof}
Let $\varphi$ be the even convex function on $\R$ with
\[
 \varphi(z)=\frac{z^2}{4A}\quad(\lvert z\rvert\le h_j/2),\qquad
 \varphi(z)=\frac{h_j^2}{16A}+\lvert z\rvert-\frac{h_j}2
 \quad(\lvert z\rvert\ge h_j/2).
\]
For every $x\in\R$, let $z$ be the unique minimizer of
$z'\mapsto\varphi(z')+(x-z')^2/2$ and set $g=x-z$. Both $z$ and $g$
depend continuously on $x$ with $g\in\partial\varphi(z)$. At $x$ take
the triple $(z,g,\varphi(z))$. Since $\varphi$ is convex with
$\varphi\ge\varphi(0)=0$, \eqref{eq:tail-planes} holds for every two
triples among $(0,0,0)$ and the triples $(z,g,\varphi(z))$ at any
finite set of points $x$. We have
\[
 gx\ge0,\qquad\lvert g\rvert\le1,\qquad
 \lvert z\rvert\le\lvert x\rvert\le2\lvert z\rvert .
\]
For $\lvert z\rvert<h_j/2$ we have
$\varphi(z)=gz-\varphi(z)=z^2/(4A)$. For $\lvert z\rvert\ge h_j/2$ we
have $\varphi(z)\ge h_j^2/(16A)$ and $gz-\varphi(z)\ge h_j^2/(16A)$.
Hence
\begin{equation}\label{sc:eq:ramp-min}
 \min\{\varphi(z),gz-\varphi(z)\}\ge(c/A)\min\{\lvert x\rvert^2,h_j^2\}.
\end{equation}
Take
\[
 x_{j-1}=\frac{h_j(\kappa-1)}{\kappa-h_j}\in[h_j,2h_j].
\]
Then $z_{j-1}=x_{j-1}-1>h_j/2$ and $g_{j-1}=1$. Hence
\eqref{sc:eq:tail-update} gives $x_j=0$.

Let $k\le i<j-1$ and $h_{i+1}\ne\kappa$. Define
$\Psi(x_i)=x_i-h_{i+1}(x_i/\kappa+(1-1/\kappa)g_i)$ for $x_i\in\R$.
Since $g_i$ depends continuously on $x_i$ with $\lvert g_i\rvert\le1$
and $1-h_{i+1}/\kappa\ne0$, $\Psi$ is continuous on $\R$ with opposite
infinite limits as $x_i\to\pm\infty$. Hence for every $x_{i+1}$ there is $x_i$ with $\Psi(x_i)=x_{i+1}$. Thus
\eqref{sc:eq:tail-update} holds. For
$h_{i+1}<\kappa$ we have $\Psi(0)=0$ and $\Psi(x_i)\to+\infty$ as
$x_i\to+\infty$. Hence for $h_{i+1}<\kappa$ and $x_{i+1}>0$ we choose
$x_i>0$. Since $\lvert g_i\rvert\le1$, \eqref{sc:eq:tail-update} gives
\begin{equation}\label{sc:eq:ramp-back}
 \lvert x_i\rvert+A\le\frac{1+\lvert h_{i+1}\rvert/A}
 {\min\{1,\lvert1-h_{i+1}/\kappa\rvert\}}(\lvert x_{i+1}\rvert+A).
\end{equation}
For $2\lvert x_i\rvert\le h_j$ we have $g_i=x_i/(2A+1)$. Thus for
every $x_i$,
\begin{equation}\label{sc:eq:ramp-forward}
 \min\{h_j,2\lvert x_{i+1}\rvert\}
 \le(1+\lvert h_{i+1}\rvert/A)\min\{h_j,2\lvert x_i\rvert\}.
\end{equation}

Assume \eqref{sc:eq:attach-setting}. Then
\[
 \sum_{u\in I}\bigl(\lvert h_u\rvert/A+2h_u^+/\kappa\bigr)\le2 .
\]
Hence multiplying \eqref{sc:eq:ramp-back} evaluated at $r$ over
$i\le r\le j-2$ gives $x_i\le CA$ and multiplying
\eqref{sc:eq:ramp-forward} evaluated at $r$ over $i\le r\le j-2$ gives
$x_i\ge ch_j$. Thus \eqref{sc:eq:ramp-step} holds.
Since $x_i>0$ for $k\le i<j$, \eqref{sc:eq:ramp-slope} holds. By
\eqref{sc:eq:ramp-min} and \eqref{sc:eq:ramp-step} we obtain
\eqref{sc:eq:ramp-height}.

Assume $2<h_j\le A\le\kappa/2$ and $h_u\ne\kappa$ for $u\in I$.
Choose $x_{j-2},\ldots,x_k$ with \eqref{sc:eq:tail-update}.
Then apply \eqref{eq:scaling} with
\[
 \sigma=h_j^{-1}\prod_{u\in I}(1+\lvert h_u\rvert/A).
\]
Before \eqref{eq:scaling} is applied, \eqref{sc:eq:ramp-forward} gives
\[
 \min\{h_j,2\lvert x_i\rvert\}\prod_{u\in I}(1+\lvert h_u\rvert/A)\ge h_j .
\]
By \eqref{sc:eq:ramp-min}, \eqref{sc:eq:ramp-value} holds after
\eqref{eq:scaling} is applied. Before \eqref{eq:scaling} is applied,
\eqref{sc:eq:ramp-back} gives
\[
 \lvert x_i\rvert+A\le CA\prod_{u\in I}
 \frac{1+\lvert h_u\rvert/A}{\min\{1,\lvert1-h_u/\kappa\rvert\}} .
\]
Since $\lvert z_i\rvert\le\lvert x_i\rvert$ and
$\lvert g_i\rvert\le1$, \eqref{sc:eq:ramp-cost} holds after
\eqref{eq:scaling} is applied.
\end{proof}

\subsection{Adding one interval of steps}\label{sc:sec:attach}\label{sc:sec:attachment}

Let $0\le k<j\le n$ and $I=\{k+1,\ldots,j-1\}$. For a vector $x$ in
dimension four, $(x)_r$ is coordinate $r$ and $(x)_{1:2}$, $(x)_{1:3}$,
$(x)_{3:4}$ are the vectors of the indicated coordinates. Let
$C_0>0$. We say that a solution $\mathcal S$ with nonzero first point
$x'$ satisfies \eqref{sc:eq:start-projections} with constant $C_0$ if
\begin{equation}\label{sc:eq:start-projections}
 \ip{\widehat{x'}}z\ge0,\qquad
 \norm z+A\norm g\le C_0(A/P)^2\norm{x'}
\end{equation}
for every triple $(z,g,q)$ of $\mathcal S$.
In every statement that assumes \eqref{sc:eq:start-projections} with
constant $C_0$, the constants $C$ and $c$ and the lower bound for $A$
or for $y$ depend on $C_0$. A rotation and \eqref{eq:scaling} preserve
\eqref{sc:eq:start-projections} with constant $C_0$.

\Cref{sc:lem:ordinary} constructs a four-dimensional solution on
$(k,j]$ with $\norm{x_j}/\norm{x_k}$ at least $h_j/(m_I+P)$ up to the
factor $\exp(-C/\log A-m_I/\kappa)$. Coordinates $1$ and $2$ come from
\cref{sc:lem:arc} and give $m_I$. Coordinate $3$ comes from
\cref{sc:lem:ramp}\ref{it:ramp-good} and gives the lower bound for
$\min\{q_i,\ip{g_i}{z_i}-q_i\}$ used in
\cref{sc:lem:ordinary}\ref{it:ord-join}. Coordinate $4$ makes $x_k$
orthogonal to $x_j$. \Cref{sc:lem:ordinary}\ref{it:ord-boundary} is
used in the proof of \cref{sc:lem:two-ranges}.
\Cref{sc:lem:ordinary}\ref{it:ord-buffer} is used in the proof of
\cref{sc:lem:bad-block}.

\begin{lemma}\label{sc:lem:ordinary}
Let $0\le k<j\le n$ and $I=\{k+1,\ldots,j-1\}$. Assume
\eqref{sc:eq:attach-setting} and $A^{19/20}\le P\le h_j$. There is an
absolute constant $C_{\mathrm{ord}}$ such that for sufficiently large
$A$ there are a point $x_j\in\R^4$ with $\norm{x_j}\ge ch_j$ and a
four-dimensional solution $\mathcal S$ on $(k,j]$ with last point $x_j$
and nonzero first point $x_k\perp x_j$ such that:
\begin{enumerate}
\item\label{it:ord-gain} $\mathcal S$ satisfies
\eqref{sc:eq:start-projections} with constant $C_{\mathrm{ord}}$ and
\begin{equation}\label{sc:eq:ordinary-gain}
 \frac{\norm{x_j}}{\norm{x_k}}\ge\frac{h_j}{m_I+P}
                    \exp\bigl(-C/\log A-m_I/\kappa\bigr).
\end{equation}
\item\label{it:ord-join} Let $C_0>0$ and let $A$ be sufficiently large
depending on $C_0$. Let $\tilde{\mathcal S}$ be the image of
$\mathcal S$ under a rotation and \eqref{eq:scaling}. Let $j<b\le n$
and let $\mathcal S'$ be a solution on $(j,b]$ such that the first
point of $\mathcal S'$ is the last point of $\tilde{\mathcal S}$ and
$\mathcal S'$ satisfies \eqref{sc:eq:start-projections} with constant
$C_0$. If $\tilde{\mathcal S}$ and $\mathcal S'$ satisfy
\cref{sc:lem:assembly}\ref{it:asm-shared}--\ref{it:asm-perp}, then
\eqref{eq:tail-planes} holds for every triple of $\tilde{\mathcal S}$
together with every triple of $\mathcal S'$.
\item\label{it:ord-boundary} The first point $x_k$ and the triples
$(z_i,g_i,q_i)$ of $\mathcal S$ satisfy
\begin{gather}
 \norm{x_k}\le e^{m_I/\kappa}(1+C/\log A)(m_I+P),\label{sc:eq:ord-norm}\\
 \norm{g_i}\le C,\qquad\ip{x_k}{g_i}\ge0,\qquad
 0\le\ip{x_k}{z_i}\le\norm{x_k}q_i+P^2\qquad(k\le i<j),
 \label{sc:eq:ord-boundary}\\
 \norm{x_k}\ge h_w/2-C\sum_{k<u<w}h_u^-\qquad(w\in I,\ h_w>1).
 \label{sc:eq:ord-lower}
\end{gather}
\item\label{it:ord-buffer} For every $R\ge P$ there is a
four-dimensional solution $\mathcal S_R$ on $(k,j]$ with last point
$x_j$ and nonzero first point $x_k'\perp x_j$ such that
$\norm{x_k'}\le C(A/P)^2R$ and $\mathcal S_R$ satisfies
\eqref{sc:eq:start-projections} with constant $C_{\mathrm{ord}}$. Let
$(z,g,q)$ be a triple with $0\le q\le\ip gz$ and
$\norm z+A\norm g\le R$ such that every triple $(z_i,g_i,q_i)$ of
$\mathcal S_R$ satisfies
\begin{equation}\label{sc:eq:buffer-support}
 \ip{g_i}z=\ip{g_i}{\widehat{x_j}}\ip{\widehat{x_j}}z,\qquad
 \ip g{z_i}=\ip g{\widehat{x_j}}\ip{\widehat{x_j}}{z_i}.
\end{equation}
Then \eqref{eq:tail-planes} holds for every triple of $\mathcal S_R$
together with $(z,g,q)$.
\end{enumerate}
\end{lemma}

From here on a solution given by \cref{sc:lem:ordinary} on $(k,j]$
denotes a four-dimensional solution $\mathcal S$ on $(k,j]$ with last
point $x_j$, $\norm{x_j}\ge ch_j$, nonzero first point $x_k\perp x_j$
and the properties \ref{it:ord-gain}--\ref{it:ord-buffer} of
\cref{sc:lem:ordinary}.

\begin{proof}
Let $K$ and $x_k,\ldots,x_{j-1}$ be given by \cref{sc:lem:arc}. At each
$x=x_k,\ldots,x_{j-1}$ take the triple $(z,g,q)$ with $g=\Pi_K(x)$,
$z=x-g$ and $q=\ip gz$. Since $x_{j-1}=(1/A,1)\in K$,
\eqref{sc:eq:tail-update} gives $x_j=(1-h_j)(1/A,1)$. The points
$x_k,\ldots,x_j$ and the triples at $x_k,\ldots,x_{j-1}$ form a
two-dimensional solution $\mathcal S_1$ on $(k,j]$. Every interval
$J=\{i+1,\ldots,j-1\}$ with $k\le i<j$ satisfies
$m_J\le2\sum_{u\in I}h_u^+\le2A$ by \eqref{sc:eq:mass-product}. Hence
\eqref{sc:eq:arc-bounds} and \eqref{sc:eq:arc-height} give
$\norm{x_i}\le CA$ and $(x_i)_2\le C\sqrt A$ for $k\le i<j$.

Let $\sigma\ge c_1P/A$ where $c_1>0$ is a fixed sufficiently small
absolute constant. Let $\mathcal S_2$ be the image under
\eqref{eq:scaling} with factor $\sigma$ of a solution on $(k,j]$ given
by \cref{sc:lem:ramp}\ref{it:ramp-good}. Let $\mathcal S_3$ be a
one-dimensional solution on $(k,j]$ with points
$x_{i+1}=(1-h_{i+1}/\kappa)x_i$ and triples $(x_i,0,0)$. Place
$\mathcal S_1$ in coordinates $1$ and $2$, $\mathcal S_2$ in coordinate
$3$ and $\mathcal S_3$ in coordinate $4$. Two applications of
\eqref{sc:eq:piece-sum} give a four-dimensional solution
$\mathcal S(\sigma)$ on $(k,j]$. Write $x_i$ and $(z_i,g_i,q_i)$ for
the points and triples of $\mathcal S(\sigma)$. Then
\begin{gather}
 (g_i)_{1:2}=\Pi_K((x_i)_{1:2}),\qquad (g_i)_4=0,\qquad z_i=x_i-g_i ,
 \label{sc:eq:attach-data}\\
 \norm{(x_i)_{1:2}}\le CA,\qquad 0\le(x_i)_2\le C\sqrt A\qquad(k\le i<j),
 \label{sc:eq:geom-sizes}\\
 (x_j)_{1:3}=(1-h_j)(1/A,1,0).\label{sc:eq:attach-Y}
\end{gather}
By \eqref{sc:eq:attach-setting}, every product of consecutive factors
$1-h_i/\kappa$ with $k<i\le j$ lies between two positive absolute
constants. Choose $(x_k)_4\ge0$ such that
\[
 (x_k)_4(x_j)_4=(h_j-1)\big((x_k)_1/A+(x_k)_2\big).
\]
Then $x_k\perp x_j$ by \eqref{sc:eq:attach-Y}. The points of
$\mathcal S_1$ and of $\mathcal S_3$ do not depend on $\sigma$ and
$(x_j)_3=0$. Hence $x_j$ does not depend on $\sigma$. By
\eqref{sc:eq:geom-sizes},
\begin{equation*}
 (x_i)_4^2\le Ch_j\sqrt A\qquad(k\le i\le j).
\end{equation*}
By \eqref{sc:eq:attach-Y},
$h_j-1\le\norm{(x_j)_{1:3}}\le h_j$. Therefore
\[
 h_j-1\le\norm{x_j}
 \le\norm{(x_j)_{1:3}}+\frac{(x_j)_4^2}{2\norm{(x_j)_{1:3}}}
 \le h_j+C\sqrt A .
\]

By \eqref{sc:eq:piece-sum}, \cref{sc:lem:arc}\ref{sc:it:arc-signs},
\eqref{sc:eq:ramp-step} and \eqref{sc:eq:ramp-height},
\begin{equation}\label{sc:eq:geom-scalar}
 q_i\ge\norm{(z_i)_{1:2}}+c(x_{j-1})_3^2/A,\qquad
 q_i-\ip{g_i}{z_i}\le-c(x_{j-1})_3^2/A .
\end{equation}
All coordinates of $x_k,z_i,g_i$ are nonnegative. By
\eqref{sc:eq:ramp-step} and \eqref{sc:eq:ramp-slope},
\[
 \norm{x_k}\ge ch_j(g_{j-1})_3,\qquad
 \norm{z_i}+A\norm{g_i}\le CA\bigl(1+(g_{j-1})_3\bigr).
\]
Since $(g_{j-1})_3=\sigma\ge c_1P/A$ and $h_j\ge P$, we have $x_k\ne0$
and $\mathcal S(\sigma)$ satisfies \eqref{sc:eq:start-projections} with
an absolute constant $C_{\mathrm{ord}}>0$ for every $\sigma\ge c_1P/A$.

Let $(z,g,q)$ be a triple with $0\le q\le\ip gz$ such that
\eqref{sc:eq:buffer-support} holds for $k\le i<j$. Since
$(x_j)_3=0$, \eqref{sc:eq:attach-data} and \eqref{sc:eq:attach-Y} give
$\ip{x_j}{g_i}=(1-h_j)((g_i)_1/A+(g_i)_2)$ and
$\ip{x_j}{z_i}=(1-h_j)((z_i)_1/A+(z_i)_2)+(x_j)_4(x_i)_4$. Dividing by
$\norm{x_j}\ge h_j-1$ and applying \eqref{sc:eq:geom-sizes},
\cref{sc:lem:arc}\ref{sc:it:arc-size} and $(x_i)_4^2\le Ch_j\sqrt A$
yields
\[
 -2\le\ip{\widehat{x_j}}{g_i}\le0,\qquad
 \lvert\ip{\widehat{x_j}}{z_i}\rvert\le C\sqrt A .
\]
Hence by \eqref{sc:eq:buffer-support} and \eqref{sc:eq:geom-scalar},
\eqref{eq:tail-planes} holds for the triple at $i$ together with
$(z,g,q)$ if
\begin{equation}\label{sc:eq:ordinary-join}
 2\bigl(\ip{-\widehat{x_j}}z\bigr)^+
   + C\sqrt A\,|\ip{\widehat{x_j}}g|
 \le c(x_{j-1})_3^2/A .
\end{equation}

Take $\mathcal S=\mathcal S(\sigma)$ with $\sigma=c_1P/A$. By
\eqref{sc:eq:ramp-step}, $\lvert(x_i)_3\rvert\le Cc_1P$. Since
$(x_i)_4^2\le Ch_j\sqrt A\le CA^{3/2}$, $A^{3/4}/P\le A^{-1/5}$,
$c_1$ is sufficiently small and $A$ is sufficiently large,
\begin{equation}\label{sc:eq:aux-coordinates}
 \norm{(z_i)_{3:4}}\le\norm{(x_i)_{3:4}}\le P/2\qquad(k\le i<j).
\end{equation}
By \eqref{sc:eq:arc-bounds} for $(x_k)_{1:2}$ and
\eqref{sc:eq:aux-coordinates}, \eqref{sc:eq:ord-norm} holds. Since
$1/h_j\le A^{-19/20}$ and $\norm{x_j}\ge h_j-1$,
\eqref{sc:eq:ordinary-gain} holds. Since $h_j>2$, we have
$\norm{x_j}\ge h_j/2$. Thus \ref{it:ord-gain} holds.

Let $C_0>0$ and let $A$ be sufficiently large depending on $C_0$. Let
$\tilde{\mathcal S}$ and $\mathcal S'$ satisfy the hypotheses of
\ref{it:ord-join}. One rotation and \eqref{eq:scaling} with one
$\sigma$ applied to $\tilde{\mathcal S}$ and to $\mathcal S'$ preserve
\eqref{sc:eq:start-projections} with constant $C_0$,
\cref{sc:lem:assembly}\ref{it:asm-shared}--\ref{it:asm-perp} and
\eqref{eq:tail-planes}. Hence we may assume
$\tilde{\mathcal S}=\mathcal S$. Let $(z,g,q)$ be a triple of
$\mathcal S'$. By \eqref{eq:adjacent-support} for $\mathcal S$ and
$\mathcal S'$, \eqref{sc:eq:buffer-support} holds for $k\le i<j$.
By $\norm{x_j}\le CA$ and \eqref{sc:eq:start-projections} for
$\mathcal S'$,
\[
 \bigl(\ip{-\widehat{x_j}}z\bigr)^+=0,\qquad
 C\sqrt A\,\lvert\ip{\widehat{x_j}}g\rvert\le CA^{5/2}/P^2 .
\]
Since $(x_{j-1})_3\ge cP^2/A$ and $P^6/A^{11/2}\ge A^{1/5}$, we have
$c(x_{j-1})_3^2/A\ge cP^4/A^3>CA^{5/2}/P^2$ for sufficiently large $A$
depending on $C_0$. Hence \eqref{sc:eq:ordinary-join} holds and \eqref{eq:tail-planes}
holds for every triple of $\mathcal S$ together with $(z,g,q)$. Thus
\ref{it:ord-join} holds.

By \eqref{sc:eq:attach-data},
\cref{sc:lem:arc}\ref{sc:it:arc-size} and \eqref{sc:eq:ramp-slope},
$\norm{g_i}\le C$. Because all coordinates of $x_k,z_i,g_i$ are
nonnegative, $\ip{x_k}{g_i}\ge0$ and $\ip{x_k}{z_i}\ge0$. By
\eqref{sc:eq:aux-coordinates} and \eqref{sc:eq:geom-scalar},
\[
 \ip{x_k}{z_i}=\ip{(x_k)_{1:2}}{(z_i)_{1:2}}+\ip{(x_k)_{3:4}}{(z_i)_{3:4}}
 \le\norm{x_k}q_i+P^2 .
\]
Thus \eqref{sc:eq:ord-boundary} holds. Let $w\in I$ have $h_w>1$. Since
$(x_w)_{1:2}\in(0,\infty)^2$ by \cref{sc:lem:arc} and
$(x_{w-1})_{1:2}\in K$ would give
$(x_w)_{1:2}=(1-h_w)(x_{w-1})_{1:2}\notin(0,\infty)^2$ for $h_w>1$, we
have $(x_{w-1})_{1:2}\notin K$. Because $K$ contains
$\{\xi\in[0,\infty)^2:\norm\xi\le1\}$, we have
$\norm{(g_{w-1})_{1:2}}\ge1$. By \eqref{sc:eq:tail-update} and
nonnegativity of the first two coordinates of $x_w$ and $x_{w-1}$,
$\norm{(x_{w-1})_{1:2}}\ge(1-1/\kappa)h_w\ge h_w/2$. For $k<u\le w-1$
with $h_u\ge0$, $(x_u)_{1:2}\le(x_{u-1})_{1:2}$ coordinatewise. For
$k<u\le w-1$ with $h_u<0$, \eqref{sc:eq:tail-update},
\eqref{sc:eq:geom-sizes}, $A\le\kappa/2$ and
\cref{sc:lem:arc}\ref{sc:it:arc-size} give
$\norm{(x_u)_{1:2}}-\norm{(x_{u-1})_{1:2}}\le Ch_u^-$. Summing over
$k<u\le w-1$ gives \eqref{sc:eq:ord-lower}. Thus \ref{it:ord-boundary} holds.

Let $R\ge P$ and choose $\sigma$ such that
$(x_{j-1})_3=C_2(A/P)R$ for a sufficiently large absolute constant
$C_2$. By \eqref{sc:eq:ramp-step} and $h_j\le A$, we have
$\sigma\ge(x_{j-1})_3/(2h_j)\ge c_1P/A$. Take $\mathcal S_R=\mathcal S(\sigma)$ and write $x_k'$ for the
first point of $\mathcal S_R$. By \eqref{sc:eq:geom-sizes},
\eqref{sc:eq:ramp-step}, $h_j\ge P$ and $R\ge P$,
\[
 \norm{x_k'}\le CA\bigl(1+(x_{j-1})_3/h_j\bigr)\le C(A/P)^2 R .
\]
Let $(z,g,q)$ be a triple with $0\le q\le\ip gz$,
$\norm z+A\norm g\le R$ and \eqref{sc:eq:buffer-support} for
$k\le i<j$. Then
$2(\ip{-\widehat{x_j}}z)^++C\sqrt A\,\lvert\ip{\widehat{x_j}}g\rvert\le CR$
and $c(x_{j-1})_3^2/A=cC_2^2(AR/P^2)R\ge cC_2^2 R$. For sufficiently
large $C_2$, \eqref{sc:eq:ordinary-join} holds. Thus
\ref{it:ord-buffer} holds.
\end{proof}

\section{Trajectory bounds for intervals between large steps}\label{sc:sec:cuts}

In this section we combine solutions on all intervals at one threshold
$e^y$ by \cref{sc:lem:assembly}. On the function of
\cref{sc:lem:assembly}, \eqref{eq:guarantee} gives
$\norm{x_n}\le\norm{x_0}$. With \eqref{sc:eq:ordinary-gain} and
\cref{sc:lem:bad-block} we obtain \eqref{sc:eq:working-quantitative}, a
lower bound
for $\sum_{h_j>e^y}\log[(m_I+P)/h_j]$. \Cref{sc:sec:deviations}
integrates \eqref{sc:eq:working-quantitative} over $y$.

Let $y$ be sufficiently large with $e^y<h_n$ and define
\begin{equation}\label{eq:cut-scale}
 A=\min\{e^y y^3,\kappa/2\},\qquad P=e^{2y}/A .
\end{equation}
Since $h_n\le\kappa/4$, we have $2e^y\le A\le e^y y^3$. Hence
$A/y^6\le P\le e^y/2$. Since $y\le\log A$,
\[
 \log A=y+O(\log y),\qquad A^{19/20}\le P<e^y .
\]
Every cut $j$ at threshold $e^y$ has $h_j>e^y>P>2$. Thus the hypothesis
$A^{19/20}\le P\le h_j$ of \cref{sc:lem:ordinary} holds for every cut
$j$ at threshold $e^y$. Let $D_0$ be the set of the cuts $j$ at
threshold $e^y$ with
\[
 \sum_{k<i\le j}\lvert h_i\rvert>A .
\]
For a cut $j\notin D_0$ the hypothesis \eqref{sc:eq:attach-setting} of
\cref{sc:lem:ordinary} holds. For a cut $j$ at threshold $e^y$ we use
the condition
\begin{equation}\label{eq:three-cuts}
 \text{$j\notin D_0$, the cut before $j$ is not in $D_0$ and the cut
 after $j$ is not in $D_0$.}
\end{equation}
For the smallest cut $j$ the condition on the cut before $j$ is
omitted.
For $j=n$ the condition on the cut after $j$ is omitted. The cuts
without \eqref{eq:three-cuts} are treated with \cref{sc:lem:bad-block}.
The conditions on the cut before $j$ and on the cut after $j$ in
\eqref{eq:three-cuts} give hypotheses \ref{it:bad-hyp-a} and
\ref{it:bad-hyp-b} of \cref{sc:lem:bad-block}.

\begin{lemma}\label{sc:lem:bad-block}
Let $C_0>0$, let $y$ be sufficiently large depending on $C_0$ with
$e^y<h_n$ and let $A$ and $P$ be given by \eqref{eq:cut-scale}. Let
$0\le a<b\le n$. Assume that:
\begin{enumerate}
\item\label{it:bad-hyp-a} if $a>0$, then some $a<j<b$ satisfies $h_j>P$
and $\sum_{a<i\le j}\lvert h_i\rvert\le A$.
\item\label{it:bad-hyp-b} if $b<n$, then $P<h_b\le A$.
\end{enumerate}
Then there is a solution $\mathcal S$ on $(a,b]$ with nonzero first
point $x_a$ and nonzero last point $x_b$, with $x_a\perp x_b$ when
$a>0$ and $b<n$, such that:
\begin{enumerate}[label=\textup{(\roman*)}]
\item\label{it:bad-cost} with a constant $C$ depending on $C_0$,
\begin{equation}\label{eq:run-cost}
 \log\frac{\norm{x_a}}{\norm{x_b}}\le C\log(A/P)
 +C\sum_{a<i\le b}\bigl[\log(1+\lvert h_i\rvert/A)
                         +(\log\lvert1-h_i/\kappa\rvert)^-\bigr].
\end{equation}
\item\label{it:bad-left} Let $a>0$ and let $\tilde{\mathcal S}$ be the
image of $\mathcal S$ under a rotation and \eqref{eq:scaling}. Let
$0\le k'<a$ and let $\mathcal S'$ be the image under a rotation and
\eqref{eq:scaling} of a solution given by \cref{sc:lem:ordinary} on
$(k',a]$ such that the last point of $\mathcal S'$ is the first point of
$\tilde{\mathcal S}$. If $\mathcal S'$ and $\tilde{\mathcal S}$ satisfy
\cref{sc:lem:assembly}\ref{it:asm-shared}--\ref{it:asm-perp}, then
\eqref{eq:tail-planes} holds for every triple of $\mathcal S'$ together
with every triple of $\tilde{\mathcal S}$.
\item\label{it:bad-right} Let $b<n$ and let $\tilde{\mathcal S}$ be the
image of $\mathcal S$ under a rotation and \eqref{eq:scaling}. Let
$\mathcal S'$ be a solution such that the first point of $\mathcal S'$
is the last point of $\tilde{\mathcal S}$ and $\mathcal S'$ satisfies
\eqref{sc:eq:start-projections} with constant $C_0$. If
$\tilde{\mathcal S}$ and $\mathcal S'$ satisfy
\cref{sc:lem:assembly}\ref{it:asm-shared}--\ref{it:asm-perp}, then
\eqref{eq:tail-planes} holds for every triple of $\tilde{\mathcal S}$
together with every triple of $\mathcal S'$.
\end{enumerate}
\end{lemma}
\begin{proof}
One rotation and \eqref{eq:scaling} with one $\sigma$ applied to
$\tilde{\mathcal S}$ and to $\mathcal S'$ preserve the hypotheses and
the conclusions of \ref{it:bad-left} and \ref{it:bad-right}. Hence it
suffices to prove \ref{it:bad-left} and \ref{it:bad-right} for
$\tilde{\mathcal S}=\mathcal S$. By
\cref{sc:lem:normalization}\ref{it:norm-top}, $h_i\ne\kappa$ for every
$i$. Set $j=0$ if $a=0$. If $a>0$, let $j$ be an index with
$a<j<b$, $h_j>P$ and $\sum_{a<i\le j}\lvert h_i\rvert\le A$. We
construct a solution $\mathcal S_0$ on $(j,b]$ with $\norm{x_b}=1$.
Take the first coordinate along $x_b$ and define $x_i$ for
$j\le i\le b$ on the line through $x_b$ by
\[
 (x_i)_1=(1-h_i/\kappa)(x_{i-1})_1\qquad(j<i\le b).
\]
The points $x_j,\ldots,x_b$ and the triples $(x_i,0,0)$ for
$j\le i<b$ form a one-dimensional solution $\mathcal S_1$ on $(j,b]$
with
\begin{equation}\label{eq:run-line}
 \log\max_{j\le i\le b}\lvert(x_i)_1\rvert
 \le\sum_{j<u\le b}(\log\lvert1-h_u/\kappa\rvert)^- .
\end{equation}
If $b=n$, let $\mathcal S_0=\mathcal S_1$. If $b<n$, take a line $U$
orthogonal to $x_b$. Let $\mathcal S_2$ be the image under
\eqref{eq:scaling} with
\[
 \sigma=C_1(A/P)\sqrt{\max_{j\le i\le b}\lvert(x_i)_1\rvert}
\]
of a solution on $(j,b]$ in $U$ given by
\cref{sc:lem:ramp}\ref{it:ramp-bad}, where $C_1$ is a sufficiently
large constant depending on $C_0$. Let $\mathcal S_0$ be the sum
\eqref{sc:eq:piece-sum} of $\mathcal S_1$ and $\mathcal S_2$. From here
on $x_i$ and $(z_i,g_i,q_i)$ denote the points and the triples of
$\mathcal S_0$. For $b<n$ we have $h_b>P$. Hence
\eqref{sc:eq:ramp-cost} and \eqref{eq:run-line} give
\begin{equation}\label{eq:run-tail-size}
 \log(\norm{x_i}+\norm{z_i}+A\norm{g_i})
 \le C\log(A/P)
 +C\sum_{j<u\le b}\bigl[\log(1+\lvert h_u\rvert/A)
                         +(\log\lvert1-h_u/\kappa\rvert)^-\bigr]
\end{equation}
for every $j\le i<b$. For $b=n$ we have $g_i=0$ and \eqref{eq:run-line}
gives \eqref{eq:run-tail-size}.

Let $b<n$. Let $(z,g,q)$ be a triple with $0\le q\le\ip gz$ and
$\norm g\le C_0(A/P)^2/A$ such that $z$ and $g$ are orthogonal to $U$
and $\ip g{z_i}=\ip g{\widehat{x_b}}\ip{\widehat{x_b}}{z_i}$ for
$j\le i<b$. Since $U\perp x_b$, we have
$\ip{\widehat{x_b}}{z_i}=(x_i)_1$. By $q\le\ip gz$,
\eqref{sc:eq:ramp-value} and the choice of $C_1$,
\begin{equation}\label{eq:S0-right-in}
 q+\ip g{z_i-z}\le\ip g{\widehat{x_b}}(x_i)_1
 \le C_0(A/P)^2\max_{j\le u\le b}\lvert(x_u)_1\rvert/A\le q_i
 \qquad(j\le i<b).
\end{equation}
Since $g_i\in U$ and $z\perp U$, we have $\ip{g_i}z=0$. Since
$\ip{g_i}{z_i}-q_i>0$ by \eqref{sc:eq:ramp-value},
\begin{equation}\label{eq:S0-right-out}
 q_i+\ip{g_i}{z-z_i}=q_i-\ip{g_i}{z_i}<0\le q\qquad(j\le i<b).
\end{equation}
Under \eqref{eq:scaling} both sides of \eqref{eq:tail-planes} are
multiplied by $\sigma^2$. Hence $q+\ip g{z_i-z}\le q_i$ and
$q_i+\ip{g_i}{z-z_i}\le q$ also hold after a rotation and
\eqref{eq:scaling} of $\mathcal S_0$ when
$\norm g\le C_0(A/P)^2\norm{x_b}/A$.

Let $a=0$ and $\mathcal S=\mathcal S_0$. Then \eqref{eq:run-tail-size}
at $i=0$ gives \eqref{eq:run-cost}. Let $b<n$ and let $\mathcal S'$
satisfy the hypotheses of \ref{it:bad-right} with
$\tilde{\mathcal S}=\mathcal S$. Let $(z,g,q)$ be a triple of
$\mathcal S'$. Since $\norm{x_b}=1$,
\eqref{sc:eq:start-projections} with constant $C_0$ gives
$\norm g\le C_0(A/P)^2/A$. By \eqref{sc:eq:ramp-value}, $g_{b-1}$ is a
nonzero vector of $U$. Hence $U$ lies in the span of the vectors of
$\mathcal S$ and $U\perp x_b$. By
\cref{sc:lem:assembly}\ref{it:asm-shared}--\ref{it:asm-perp}, $U$ is
orthogonal to $z$ and to $g$. By \eqref{eq:adjacent-support} for
$\mathcal S$ and $\mathcal S'$ with shared point $x_b$,
$\ip g{z_i}=\ip g{\widehat{x_b}}\ip{\widehat{x_b}}{z_i}$. Thus
\eqref{eq:S0-right-in} and \eqref{eq:S0-right-out} give
\eqref{eq:tail-planes} for every triple of $\mathcal S$ together with
$(z,g,q)$ and \ref{it:bad-right} holds.

Let $a>0$. Since $h_j>P>2$ and $A\le\kappa/2$,
\eqref{sc:eq:attach-setting} holds with $k=a$. Let $\mathcal S^\ast$
be a solution given by \cref{sc:lem:ordinary} on $(a,j]$, let $x_j$ be
the last point of $\mathcal S^\ast$ and let $W$ be the
four-dimensional space of $\mathcal S^\ast$. Apply a rotation and
\eqref{eq:scaling} to $\mathcal S_0$ such that the first point of
$\mathcal S_0$ is $x_j$ and every vector of $\mathcal S_0$ lies in
$\operatorname{span}(x_j)\oplus W^\perp$. Define
\[
 R=\max\Bigl\{P,\;\max_{j\le i<b}(\norm{z_i}+A\norm{g_i}),\;
              C_0(A/P)^2\norm{x_b}\Bigr\},
\]
where $C_0(A/P)^2\norm{x_b}$ is omitted for $b=n$. Let $\mathcal S_R$
be the solution on $(a,j]$ in $W$ given by
\cref{sc:lem:ordinary}\ref{it:ord-buffer} for $\mathcal S^\ast$ and
$R$. Let $\mathcal S$ consist of the points and triples of
$\mathcal S_R$ on $(a,j]$ and of $\mathcal S_0$ on $(j,b]$. Every
triple $(z,g,q)$ at $j\le i<b$ has
$z,g\in\operatorname{span}(x_j)\oplus W^\perp$, $0\le q\le\ip gz$ and
$\norm z+A\norm g\le R$. The vectors of $\mathcal S_R$ lie in $W$.
Hence every triple $(z_i,g_i,q_i)$ at $a\le i<j$ and every triple
$(z,g,q)$ at an index in $[j,b)$ satisfy
$\ip{g_i}z=\ip{g_i}{\widehat{x_j}}\ip{\widehat{x_j}}z$ and
$\ip g{z_i}=\ip g{\widehat{x_j}}\ip{\widehat{x_j}}{z_i}$. By
\cref{sc:lem:ordinary}\ref{it:ord-buffer}, \eqref{eq:tail-planes}
holds for every triple at an index in $[a,j)$ together with every
triple at an index in $[j,b)$. Thus $\mathcal S$ is a solution on
$(a,b]$. The first point $x_a$ of $\mathcal S_R$ lies in $W$ and is
orthogonal to $x_j$. Every vector of $\mathcal S_0$ lies in
$\operatorname{span}(x_j)\oplus W^\perp$. Hence $x_a\perp x_b$ and
$z_i,g_i\perp x_a$ for $j\le i<b$.

By \cref{sc:lem:ordinary}\ref{it:ord-buffer},
$\norm{x_a}\le C(A/P)^2R$. Since $\norm{x_j}\ge ch_j>cP$ and
$x_j=z_j+g_j$, we have $P\le C(\norm{z_j}+A\norm{g_j})$. Before the
rotation and \eqref{eq:scaling} of $\mathcal S_0$ we had
$\norm{x_b}=1$. Hence after the rotation and \eqref{eq:scaling},
\eqref{eq:run-tail-size} holds for
$\log[(\norm{x_i}+\norm{z_i}+A\norm{g_i})/\norm{x_b}]$. Dividing
$\norm{x_a}\le C(A/P)^2R$ by $\norm{x_b}$ gives \eqref{eq:run-cost}.

Let $b<n$ and let $\mathcal S'$ satisfy the hypotheses of
\ref{it:bad-right} with $\tilde{\mathcal S}=\mathcal S$. Let $(z,g,q)$
be a triple of $\mathcal S'$. By
\cref{sc:lem:assembly}\ref{it:asm-shared}--\ref{it:asm-perp}, we can
write $z=\alpha x_b+\zeta$ and $g=\beta x_b+\zeta'$, where $\zeta$ and
$\zeta'$ are orthogonal to $x_b$ and to every vector of
$\mathcal S$. Since $j<b$, also $x_j=z_j+g_j$ is orthogonal to
$\zeta$ and to $\zeta'$.
Let $a\le i<j$. Since $g_i\in W$ and
$x_b\in\operatorname{span}(x_j)\oplus W^\perp$, we have
$\ip{g_i}{x_b}=\ip{g_i}{\widehat{x_j}}\ip{\widehat{x_j}}{x_b}$ and
therefore
\[
 \ip{g_i}z=\alpha\ip{g_i}{x_b}
 =\alpha\ip{g_i}{\widehat{x_j}}\ip{\widehat{x_j}}{x_b}
  =\ip{g_i}{\widehat{x_j}}\ip{\widehat{x_j}}{\alpha x_b+\zeta}
 =\ip{g_i}{\widehat{x_j}}\ip{\widehat{x_j}}z .
\]
Since $z_i\in W$,
\[
 \ip g{z_i}=\beta\ip{x_b}{z_i}
 =\beta\ip{x_b}{\widehat{x_j}}\ip{\widehat{x_j}}{z_i}
 =\ip g{\widehat{x_j}}\ip{\widehat{x_j}}{z_i} .
\]
Also $0\le q\le\ip gz$ and
$\norm z+A\norm g\le C_0(A/P)^2\norm{x_b}\le R$ by
\eqref{sc:eq:start-projections} with constant $C_0$. By
\cref{sc:lem:ordinary}\ref{it:ord-buffer}, \eqref{eq:tail-planes}
holds for every triple at an index in $[a,j)$ together with $(z,g,q)$.
By \eqref{sc:eq:ramp-value}, $g_{b-1}$ is a nonzero vector of the
rotated line $U$. Hence $\zeta\perp U$ and $\zeta'\perp U$. Since $U\perp x_b$,
also $z\perp U$ and $g\perp U$. By
\eqref{eq:adjacent-support} for $\mathcal S$ and $\mathcal S'$,
$\ip g{z_i}=\ip g{\widehat{x_b}}\ip{\widehat{x_b}}{z_i}$ for
$j\le i<b$. Thus \eqref{eq:S0-right-in} and \eqref{eq:S0-right-out}
give \eqref{eq:tail-planes} for every triple at an index in $[j,b)$
together with $(z,g,q)$. Hence \ref{it:bad-right} holds.

Let $\mathcal S'$ satisfy the hypotheses of \ref{it:bad-left} with
$\tilde{\mathcal S}=\mathcal S$ and let $(z',g',q')$ be a triple of
$\mathcal S'$. Let $j\le i<b$. By \eqref{eq:adjacent-support} for
$\mathcal S'$ and $\mathcal S$ and by $z_i,g_i\perp x_a$, we have
$\ip{g'}{z_i}=\ip{g_i}{z'}=0$. Since $q'-\ip{g'}{z'}\le0\le q_i$ and
$q_i-\ip{g_i}{z_i}\le0\le q'$, \eqref{eq:tail-planes} holds for
$(z',g',q')$ together with the triple at $i$. By
\cref{sc:lem:ordinary}\ref{it:ord-buffer}, $\mathcal S_R$ has first
point $x_a$ and satisfies \eqref{sc:eq:start-projections} with constant
$C_{\mathrm{ord}}$. Hence
\cref{sc:lem:ordinary}\ref{it:ord-join} applied to $\mathcal S'$ and
$\mathcal S_R$ with $C_0=C_{\mathrm{ord}}$ gives
\eqref{eq:tail-planes} for $(z',g',q')$ together with every triple at
an index in $[a,j)$. Thus \ref{it:bad-left} holds.
\end{proof}

\Cref{sc:prop:working-realization} is used in three forms. In
\eqref{sc:eq:initial-log-floor}, $S=V=\varnothing$. In
\cref{sc:lem:excess-sum}, $S=\varnothing$ and the steps of a set $V$
are used as cuts in addition to the cuts at threshold $e^y$. In
\cref{sc:lem:two-ranges} solutions given by \cref{sc:lem:four-phase}
are prescribed for the cuts of a set $S$.

\begin{proposition}\label{sc:prop:working-realization}
Let $C_0>0$, let $y$ be sufficiently large depending on $C_0$ with
$e^y<h_n$ and let $A$ and $P$ be given by \eqref{eq:cut-scale}. Let $S$
be a set of cuts at threshold $e^y$ with \eqref{eq:three-cuts} and let
$V$ be a set of positive steps $v$ with $P<h_v\le e^y$. Assume that:
\begin{enumerate}
\item\label{it:work-S} for every $j\in S$ solutions
$\mathcal S_{j,1},\ldots,\mathcal S_{j,r_j}$ on consecutive intervals
with union $(k,j]$ are given. The first and last points of each of
$\mathcal S_{j,1},\ldots,\mathcal S_{j,r_j}$ are nonzero and
orthogonal. \Cref{sc:lem:assembly}\ref{it:asm-shared}--\ref{it:asm-perp} hold for
$\mathcal S_{j,1},\ldots,\mathcal S_{j,r_j}$ and
\eqref{eq:tail-planes} holds for every two triples of
$\mathcal S_{j,1},\ldots,\mathcal S_{j,r_j}$. Also $\mathcal S_{j,1}$
satisfies \eqref{sc:eq:start-projections} with constant $C_0$ and
$\mathcal S_{j,r_j}$ is the image under a rotation and
\eqref{eq:scaling} of a solution given by \cref{sc:lem:ordinary} on
the interval of $\mathcal S_{j,r_j}$. Write
$x_k$ for the first point of $\mathcal S_{j,1}$ and $x_j$ for the last
point of $\mathcal S_{j,r_j}$.
\item\label{it:work-V} for every cut $j$ at threshold $e^y$ the set
$V\cap I$ is an initial segment of the positive steps of $I$ in the
order of removal and $V\cap I=\varnothing$ for $j\in S$.
\end{enumerate}
Then
\begin{multline}\label{sc:eq:working-quantitative}
 \sum_{j\in S}
 \log\frac{(m_I+P)\norm{x_j}}{h_j\norm{x_k}}
 +\sum_{v\in V}\log\frac{h_v(m_v+P)}{(m_{vL}+P)(m_{vR}+P)}
 -\frac{C}{y}\bigl(\#\{j\notin S:h_j>e^y\}+\#V\bigr)
 \\ \le \sum_{h_j>e^y}\log\frac{m_I+P}{h_j}
 +\frac{Cn}{e^{y/p}}\left(\frac1y+
 \frac{(e^y/h_n)^{1/(2p)}}{\log(\kappa/h_n)}\right).
\end{multline}
The constant $C$ in \eqref{sc:eq:working-quantitative} depends on
$C_0$. The number of cuts at threshold $e^y$ without
\eqref{eq:three-cuts} is at most
$Cne^{-y/p}[1/y+(e^y/h_n)^{1/(2p)}/\log(\kappa/h_n)]$ with an absolute
constant $C$.
\end{proposition}

\begin{proof}
Let $T$ be the union of $V$ and the set of the cuts at threshold $e^y$.
For $j\in T$ let $k_T$ denote the largest element of $T$ below $j$,
with $k_T=0$ for the first element of $T$. Let
$I_T=\{k_T+1,\ldots,j-1\}$. By \ref{it:work-V}, $k_T=k$ for $j\in S$.
Let $D_T$ be the set of the $j\in T$ with
\[
 \sum_{k_T<i\le j}\lvert h_i\rvert>A .
\]
For $j\in T$ let $j'$ be the smallest cut with $j'\ge j$ and let $k'$
be the largest cut below $j'$, with $k'=0$ for the first cut. Then
$(k_T,j]\subset(k',j']$. Hence $j\in D_T$ implies $j'\in D_0$. Let
$T^\ast$ be the set of the $j\in T$ such that $j$, the element of $T$
before $j$ or the element of $T$ after $j$ lies in $D_T$. Then every
element of $T^\ast$ lies in the set $(k',j']$ of a cut $j'$ without
\eqref{eq:three-cuts}. In particular $S\cap T^\ast=\varnothing$. We use
sets $(a,b]$ with
\begin{equation}\label{eq:run-def}
 \text{$(a,b]$ is a maximal union of consecutive sets $(k_T,j]$ with
 $j\in T^\ast$.}
\end{equation}
Let $(a,b]$ satisfy \eqref{eq:run-def}. Let $a>0$ and let $j$ be the
first element of $T$ in $(a,b]$. Then $a\in T\setminus T^\ast$. Hence
$j\notin D_T$ and $\sum_{a<i\le j}\lvert h_i\rvert\le A$. Since
$j\in T$, we have $h_j>P$. Since $j\in T^\ast$, $a\notin D_T$ and
$j\notin D_T$, the element of $T$ after $j$ lies in $D_T$. Hence $j<b$.
Let $b<n$. The element of $T$ after $b$ is not in $T^\ast$. Hence
$b\notin D_T$ and $P<h_b\le A$. Thus $(a,b]$ satisfies
\cref{sc:lem:bad-block}\ref{it:bad-hyp-a} and \ref{it:bad-hyp-b}.

Set $C_\ast=\max\{C_0,C_{\mathrm{ord}}\}$. Let $\Gamma$ be the family
of the sets $(k_T,j]$ with $j\in T\setminus(S\cup T^\ast)$, the
intervals of $\mathcal S_{j,i}$ with $j\in S$ and
$1\le i\le r_j$ and the sets $(a,b]$ with \eqref{eq:run-def}. Then
$\Gamma$ is a partition of $(0,n]$ into sets of consecutive indices.
For each $J\in\Gamma$ we define a solution $\mathcal S_J$ on $J$. Let
$J=(k_T,j]$ with $j\in T\setminus(S\cup T^\ast)$. Then $j\notin D_T$
and $h_j>P>2$. Hence \eqref{sc:eq:attach-setting} holds with $k=k_T$.
Let $\mathcal S_J$ be a solution given by \cref{sc:lem:ordinary} on
$(k_T,j]$. Let $J$ be the interval of $\mathcal S_{j,i}$ with $j\in S$
and $1\le i\le r_j$. Let $\mathcal S_J=\mathcal S_{j,i}$. Let
$J=(a,b]$ satisfy \eqref{eq:run-def}. Let $\mathcal S_J$ be a solution
given by \cref{sc:lem:bad-block} applied with $C_0=C_\ast$.

Let $(a,b]\in\Gamma$ satisfy \eqref{eq:run-def}. Let $a>0$ and let
$J\in\Gamma$ be the set with last index $a$. Then $J=(k_T,a]$ with
$a\in T\setminus(S\cup T^\ast)$ or $J$ is the interval of
$\mathcal S_{a,r_a}$ with $a\in S$. Hence $\mathcal S_J$ is a solution
given by \cref{sc:lem:ordinary} on $J$ or the image under a rotation
and \eqref{eq:scaling} of a solution given by \cref{sc:lem:ordinary}
on $J$. Let $b<n$ and let $J\in\Gamma$ be the set with first index
$b+1$. Then $J=(b,j]$ with $j\in T\setminus(S\cup T^\ast)$ or $J$ is
the interval of $\mathcal S_{j,1}$ with $j\in S$. Hence $\mathcal S_J$
satisfies \eqref{sc:eq:start-projections} with constant $C_\ast$.

For every $J\in\Gamma$ the first and last points of $\mathcal S_J$ are
nonzero. For every $J\in\Gamma$ with $1\notin J$ and $n\notin J$ the
first and last points of $\mathcal S_J$ are orthogonal. Hence
rotations and \eqref{eq:scaling} in a Euclidean space of sufficiently
large dimension give images $\tilde{\mathcal S}_J$ of $\mathcal S_J$
for $J\in\Gamma$ that satisfy
\cref{sc:lem:assembly}\ref{it:asm-shared}--\ref{it:asm-perp}. Let
$J_-,J_+\in\Gamma$ be consecutive with $J_-$ before $J_+$ and write
$\mathcal S_-=\tilde{\mathcal S}_{J_-}$ and
$\mathcal S_+=\tilde{\mathcal S}_{J_+}$. We show that
\eqref{eq:tail-planes} holds for every triple of $\mathcal S_-$
together with every triple of $\mathcal S_+$.

Let $J_-$ and $J_+$ be the intervals of $\mathcal S_{j,i}$ and
$\mathcal S_{j,i+1}$ for one $j\in S$ and one $1\le i<r_j$. The last
point of $\mathcal S_{j,i}$ is the first point of $\mathcal S_{j,i+1}$
and the last point of $\mathcal S_-$ is the first point of
$\mathcal S_+$. Hence the maps from $\mathcal S_{j,i}$ to
$\mathcal S_-$ and from $\mathcal S_{j,i+1}$ to $\mathcal S_+$ have
one factor $\sigma$ in \eqref{eq:scaling}. Let $(z,g,q)$ be a triple of
$\mathcal S_{j,i}$ and $(z',g',q')$ a triple of $\mathcal S_{j,i+1}$.
By \eqref{eq:adjacent-support} for $\mathcal S_{j,i}$ and
$\mathcal S_{j,i+1}$ and for $\mathcal S_-$ and $\mathcal S_+$, the
two maps multiply $\ip g{z'}$ and $\ip{g'}z$ by $\sigma^2$. The two
maps also multiply $q$, $q'$, $\ip gz$ and $\ip{g'}{z'}$ by
$\sigma^2$. By \ref{it:work-S}, \eqref{eq:tail-planes}
holds for $(z,g,q)$ together with $(z',g',q')$. Thus
\eqref{eq:tail-planes} holds for every triple of $\mathcal S_-$
together with every triple of $\mathcal S_+$.

Let $J_+=(a,b]$ satisfy \eqref{eq:run-def}. Then
\cref{sc:lem:bad-block}\ref{it:bad-left} gives \eqref{eq:tail-planes}
for every triple of $\mathcal S_-$ together with every triple of
$\mathcal S_+$. Let $J_-=(a,b]$ satisfy \eqref{eq:run-def}. Then
\cref{sc:lem:bad-block}\ref{it:bad-right} applied with $C_0=C_\ast$
gives \eqref{eq:tail-planes} for every triple of $\mathcal S_-$
together with every triple of $\mathcal S_+$.

Let neither $J_-$ nor $J_+$ satisfy \eqref{eq:run-def} and let $J_-$
and $J_+$ not be the intervals of $\mathcal S_{j,i}$ and
$\mathcal S_{j,i+1}$ for one $j\in S$ and one $1\le i<r_j$. Then
$J_-=(k_T,j]$ with $j\in T\setminus(S\cup T^\ast)$ or $J_-$ is the
interval of $\mathcal S_{j,r_j}$ with $j\in S$. Hence $\mathcal S_-$ is
the image under a rotation and \eqref{eq:scaling} of a solution given
by \cref{sc:lem:ordinary} on $J_-$. Also $J_+=(k_T,j]$ with
$j\in T\setminus(S\cup T^\ast)$ or $J_+$ is the interval of
$\mathcal S_{j,1}$ with $j\in S$. Hence $\mathcal S_+$ satisfies
\eqref{sc:eq:start-projections} with constant $C_\ast$. Thus
\cref{sc:lem:ordinary}\ref{it:ord-join} applied with $C_0=C_\ast$ gives
\eqref{eq:tail-planes} for every triple of $\mathcal S_-$ together with
every triple of $\mathcal S_+$.

Let $x_0,\ldots,x_n$ be the points of $\tilde{\mathcal S}_J$ with
$J\in\Gamma$. By \cref{sc:lem:assembly}
there is a $1$-smooth, $1/\kappa$-strongly convex function with
minimizer zero and trajectory $x_0,\ldots,x_n$.
By \eqref{eq:guarantee},
\[
 \log\frac{\norm{x_n}}{\norm{x_0}}\le0 .
\]

Let $1/A\le\lambda\le1$. We have $\log A\ge y$ and
$m_{I_T}/\kappa\le\sum_{i\in I_T}(\log\lvert1-h_i/\kappa\rvert)^-$.
Hence \eqref{sc:eq:ordinary-gain} gives for
$j\in T\setminus(S\cup T^\ast)$
\[
 \log\frac{\norm{x_j}}{\norm{x_{k_T}}}\ge\log\frac{h_j}{m_{I_T}+P}
 -\frac Cy-\sum_{i\in I_T}(\log\lvert1-h_i/\kappa\rvert)^- .
\]
For $(a,b]$ with \eqref{eq:run-def} we have \eqref{eq:run-cost} and
$\log(1+\lvert h_i\rvert/A)\le\log(1+\lambda\lvert h_i\rvert)$. For
$j\in T$ with $a<j\le b$,
\[
 \log\frac{h_j}{m_{I_T}+P}\le\log(A/P)+\log(1+\lambda h_j) .
\]
By \eqref{sc:eq:Lambda-Delta} and \ref{it:work-V}, the sum of
$\log[h_j/(m_{I_T}+P)]$ over $j\in T\setminus S$ equals the sum of
$\log[h_j/(m_I+P)]$ over the cuts $j\notin S$ at threshold $e^y$ plus
the sum over $v\in V$ in \eqref{sc:eq:working-quantitative}. We deduce
that the left side of \eqref{sc:eq:working-quantitative} is at most
$\sum_{h_j>e^y}\log[(m_I+P)/h_j]$ plus $C$ times
\begin{equation}\label{sc:eq:working-error}
 \sum_i\bigl[\log(1+\lambda\lvert h_i\rvert)
 +(\log\lvert1-h_i/\kappa\rvert)^-\bigr]
 +\log(A/P)\,\#T^\ast .
\end{equation}

Each $j\in D_T$ satisfies
\[
 \sum_{k_T<i\le j}\log(1+\lambda\lvert h_i\rvert)\ge\log(1+\lambda A).
\]
The sets $(k_T,j]$ with $j\in D_T$ are disjoint and each $j\in D_T$
gives at most three elements of $T^\ast$. Hence
$\#T^\ast\le3\#D_T\le3\sum_i\log(1+\lambda\lvert h_i\rvert)/\log(1+\lambda A)$.
The sets $(k,j]$ with $j\in D_0$ are disjoint and each $j\in D_0$ gives
at most three cuts without \eqref{eq:three-cuts}. Hence the number of
cuts without \eqref{eq:three-cuts} is at most
$3\#D_0\le3\sum_i\log(1+\lambda\lvert h_i\rvert)/\log(1+\lambda A)$.

Let $A=e^y y^3$ and $\lambda=1/(e^y y^2)$. Then $e^y y^2\le h_n$.
Otherwise $\log h_n<2y$, contrary to
$2e^y y^3\le\kappa\le h_n\log h_n$. By \eqref{sc:eq:quadratic-tails}
and \eqref{sc:eq:top-quadratic},
\[
 \sum_i\bigl[\log(1+\lambda\lvert h_i\rvert)
 +(\log\lvert1-h_i/\kappa\rvert)^-\bigr]\le Cn(e^y y^2)^{-1/p}.
\]
Here $\log(1+\lambda A)=\log(1+y)$ and $\log(A/P)=6\log y$. Since
$p<2$, both \eqref{sc:eq:working-error} and the number of cuts without
\eqref{eq:three-cuts} are at most $Cne^{-y/p}/y$.

Let $A=\kappa/2$ and $\lambda=1/\sqrt{\kappa h_n}$. Since
$h_n\le\kappa/4$, we have $1/A\le\lambda\le1$. By
\eqref{sc:eq:top-quadratic},
\[
 \sum_i\bigl[\log(1+\lambda\lvert h_i\rvert)
 +(\log\lvert1-h_i/\kappa\rvert)^-\bigr]
 \le\frac{Cnh_n^{-1/p}}{\log(\kappa/h_n)} .
\]
Here $\log(1+\lambda A)\ge c\log(\kappa/h_n)$. Hence both $\#T^\ast$
and the number of cuts without \eqref{eq:three-cuts} are at most
$C\sum_i\log(1+\lambda\lvert h_i\rvert)/\log(\kappa/h_n)$. Also
\[
 \log(A/P)\le2\log(\kappa/h_n)+2\log(h_n/e^y),\qquad
 (e^y/h_n)^{1/(2p)}\bigl(1+\log(h_n/e^y)\bigr)\le C .
\]
Hence \eqref{sc:eq:working-error} and the number of cuts without
\eqref{eq:three-cuts} are at most
$Cn(e^y h_n)^{-1/(2p)}/\log(\kappa/h_n)$. For $A=e^yy^3$ and for
$A=\kappa/2$, \eqref{sc:eq:working-quantitative} and the bound for the
number of cuts without \eqref{eq:three-cuts} hold.
\end{proof}

\section{Deviations from the silver proportions}\label{sc:sec:deviations}

\Cref{sc:prop:working-realization} bounds
$\sum_{h_j>e^y}\log[(m_I+P(y))/h_j]$ from below at one threshold $e^y$.
In this section we integrate \eqref{sc:eq:working-quantitative} against
$e^{y/p}\,dy$. When $y$ increases past $\log h_v$, the sum
$\sum_{h_j>e^y}\log[(m_I+P(y))/h_j]$ changes by
\eqref{sc:eq:Lambda-Delta} and $\mathcal P$ changes by $\Delta_v$.
\Cref{lem:max-plus-excess} bounds \eqref{sc:eq:Lambda-Delta} minus
$\Delta_v/h_v^{1/p}$ from below. \Cref{sc:lem:initial-power} gives
$\mathcal P(y)\ge(p-o(1))n$ for $y<\log h_n$ with $y/\log h_n\to1$. For
$s<t\le4s$ with $e^t t^3\le h_n$, \cref{sc:lem:excess-sum} gives
$\sum_{s<\log h_v\le t}\lvert\Delta_v\rvert\le Cn$ and bounds the sum
of $\min\{h_v,e^t\}^{1/p}$ over the steps $v$ that are not near the
silver proportions \eqref{sc:eq:silver-proportions}.

By \eqref{eq:cut-scale}, $P=P(y)=\max\{e^y/y^3,2e^{2y}/\kappa\}$. The
function $P$ is locally absolutely continuous and increases for $y>3$.
Let $v$ be a positive step with $3<\log h_v<\log h_n$. The step
$v$ is a cut at threshold $e^y$ for $y<\log h_v$ and is not a cut for
$y\ge\log h_v$. Inserting the cut $v$ changes
$\sum_{h_j>e^y}\log[h_j/(m_I+P)]$ by \eqref{sc:eq:Lambda-Delta}. Hence
removing the cut $v$ changes $\sum_{h_j>e^y}\log[(m_I+P)/h_j]$ by
\eqref{sc:eq:Lambda-Delta}. Thus the value of
$\sum_{h_j>e^y}\log[(m_I+P(y))/h_j]$ at $y=\log h_v$ minus the limit
from the left equals the sum of \eqref{sc:eq:Lambda-Delta} over the
positive steps with stepsize $h_v$, with $P=P(\log h_v)$. On an
interval of $y>3$ without a number $\log h_v$ for a positive step $v$,
the set of cuts is constant and $y\mapsto\log[(m_I+P(y))/h_j]$ is
locally absolutely continuous with nonnegative derivative
$P'(y)/(m_I+P(y))$ almost everywhere. For sufficiently large
$s<t<\log h_n$, integration by parts and dropping the nonnegative
integral $\int_s^t e^{y/p}\sum_{h_j>e^y}P'(y)/(m_I+P(y))\,dy$ therefore
give
\begin{multline}\label{eq:integrated-split}
 \sum_{s<\log h_v\le t}h_v^{1/p}
 \log\frac{h_v(m_v+P(\log h_v))}{(m_{vL}+P(\log h_v))(m_{vR}+P(\log h_v))}
 +\frac1p\int_s^t e^{y/p}\sum_{h_j>e^y}\log\frac{m_I+P(y)}{h_j}\,dy\\
 \le \left[e^{y/p}\sum_{h_j>e^y}\log\frac{m_I+P(y)}{h_j}\right]_{y=s}^{y=t},
\end{multline}
where $\sum_{s<\log h_v\le t}$ is the sum over the positive steps $v$
with $s<\log h_v\le t$.
Let $y<\log h_n$ be sufficiently large. In
\cref{sc:prop:working-realization} take $C_0=1$ and $S=V=\varnothing$.
By \eqref{sc:eq:count-bound} there are at most $ne^{-y/p}$ cuts. Hence
\begin{equation}\label{sc:eq:initial-log-floor}
 e^{y/p}\sum_{h_j>e^y}\log\frac{m_I+P(y)}{h_j}
 \ge -Cn\left(\frac1y+\frac{(e^y/h_n)^{1/(2p)}}{\log(\kappa/h_n)}\right).
\end{equation}

\subsection{Thresholds near the last stepsize}
\label{sc:sec:initial-power}

\begin{lemma}\label{sc:lem:initial-power}
Along the sequence of \cref{sc:lem:normalization}, let $s<\log h_n$
with $s/\log h_n\to1$. Then $\mathcal P(s)\ge (p-o(1))n$.
\end{lemma}
\begin{proof}
For a fixed schedule, as $y\uparrow\log h_n$ the cuts at threshold
$e^y$ are the indices $j$ with $h_j\ge h_n$ and
$P(y)\to2h_n^2/\kappa$ by \eqref{sc:eq:standing-log}. Along the
sequence of \cref{sc:lem:normalization},
$2h_n^2/\kappa=o(h_n)$ by \eqref{sc:eq:standing}. Define
\[
  G(y)=\mathcal P(y)-e^{y/p}\sum_{h_j>e^y}\log\frac{m_I+P(y)}{h_j} .
\]
By \eqref{sc:eq:threshold-set}, $G(y)$ is the sum over the cuts $j$ of
\[
 \frac p{2-\sqrt2}\,m_I^{1/p}-e^{y/p}\log\frac{m_I+P(y)}{h_j} .
\]
We have $\log\tau\le p(\tau^{1/p}-1)$ for
$\tau=(m_I+P(y))/h_j$, $(m_I+P(y))^{1/p}\le m_I^{1/p}+P(y)^{1/p}$,
$h_j\ge h_n$ and $2-\sqrt2<1$. Hence the summand of each cut $j$ has
limit at least $p(h_n^{1/p}-(2h_n^2/\kappa)^{1/p})$ as
$y\uparrow\log h_n$. By \eqref{sc:eq:normalized-length},
$\liminf_{y\uparrow\log h_n}G(y)\ge(p-o(1))n$.

Apply \eqref{eq:integrated-split} with $s<t<\log h_n$. We have
$\log[(\log h_n)/s]\to0$ and $\kappa/h_n\to\infty$. By
\eqref{sc:eq:initial-log-floor}, both
$e^{s/p}\sum_{h_j>e^s}\log[(m_I+P(s))/h_j]$ and
$\frac1p\int_s^t e^{y/p}\sum_{h_j>e^y}\log[(m_I+P(y))/h_j]\,dy$ are at
least $-o(n)$. For $s<\log h_v\le t<\log h_n$, at $y=\log h_v$ we have
$P(y)/h_v\le y^{-3}+2h_v/\kappa$. By \eqref{sc:eq:count-bound},
$\sum h_v^{1/p}\le Cn$ over the steps $v$ with $\log h_v$ in an
interval of length one inside $(s,\log h_n)$. Hence the terms
$C(P(\log h_v)/h_v)^{1/p}$ and $Ch_v/\kappa$ of
\cref{lem:max-plus-excess}\ref{it:split-mass},
multiplied by $h_v^{1/p}$ and summed over $s<\log h_v\le t$, are at
most
\[
 Cn\left(\int_{s-1}^{\infty}y^{-3/p}\,dy+(h_n/\kappa)^{1/p}\right)=o(n).
\]
By \cref{lem:max-plus-excess}\ref{it:split-mass} and
\eqref{sc:eq:band-sums},
\[
 \sum_{s<\log h_v\le t}h_v^{1/p}
 \log\frac{h_v(m_v+P(\log h_v))}{(m_{vL}+P(\log h_v))(m_{vR}+P(\log h_v))}
 \ge\mathcal P(t)-\mathcal P(s)-o(n).
\]
Thus \eqref{eq:integrated-split} gives $\mathcal P(s)\ge G(t)-o(n)$.
Letting $t\uparrow\log h_n$ proves \cref{sc:lem:initial-power}.
\end{proof}

\subsection{Splits near the silver proportions}\label{sc:sec:bands}\label{sc:sec:balanced}\label{sc:sec:balanced-splits}

In the rest of \cref{sc:sec:deviations} and in
\cref{sec:consecutive-splits} we assume that $y$ is sufficiently large
and that
\begin{equation}\label{sc:eq:band-B}
 e^y y^3\le h_n .
\end{equation}
Then \eqref{eq:cut-scale} gives $A=e^yy^3$, $P=P(y)=e^y/y^3$ and
$e^y/\kappa\le e^y/h_n\le y^{-3}$.

Let $v$ be a positive step with $e^y/y\le h_v\le e^y$. Then
$P(y)/h_v\le y^{-2}$. Hence \cref{lem:max-plus-excess}\ref{it:split-mass}
with $P=P(y)$ gives
\begin{equation}\label{sc:eq:bulk-split}
 \log\frac{h_v(m_v+P(y))}{(m_{vL}+P(y))(m_{vR}+P(y))}
 \ge \frac{\Delta_v}{h_v^{1/p}}-Cy^{-2/p}.
\end{equation}
Fix a sufficiently large absolute constant $y_0>3$. Let $B$ be the set
of the positive steps $v$ such that $\log h_v\ge y_0$,
$h_v(\log h_v)^3\le h_n$ and
\begin{equation}\label{eq:B-def}
 \frac{m_v}{h_v}\ge\rho^2-(\log h_v)^{-1/4},\qquad
 \frac1{h_v}\sum_{i\in I_w}h_i^+\le\rho+(\log h_v)^{-1/4}
 \quad(w=vL,vR).
\end{equation}
Let $B'$ be the set of the positive steps $v$ such that $v\notin B$ or
a child of $v$ is not in $B$. By \cref{sc:lem:balance}, every $v\in B$
satisfies \eqref{sc:eq:silver-proportions} up to factors
$1+O((\log h_v)^{-1/4})$.

\begin{lemma}\label{sc:lem:balance}
Let $v\in B$ and $y=\log h_v$. Then:
\begin{enumerate}
\item\label{it:bal-prop} $m_v=(\rho^2+O(y^{-1/4}))h_v$ and
$m_w=(\rho+O(y^{-1/4}))h_v$ for $w=vL,vR$. In particular $I_{vL}$ and
$I_{vR}$ each contain a positive step.
\item\label{it:bal-parent} if the parent $r$ of $v$ lies in $B$, then
$\lvert\log(h_r/(\rho h_v))\rvert\le Cy^{-1/4}$.
\item\label{it:bal-debt} the negative steps in $I_{vL}$ satisfy
\begin{equation}\label{sc:eq:balanced-debt}
 \sum_{i\in I_{vL}}h_i^-\le Cy^{-1/4}h_v .
\end{equation}
\end{enumerate}
\end{lemma}
\begin{proof}
By \eqref{eq:B-def}, \eqref{sc:eq:mass-product} and
$h_v/\kappa\le y^{-3}$,
\[
 m_w\le(1+h_v/\kappa)\sum_{i\in I_w}h_i^+\le(\rho+Cy^{-1/4})h_v
 \qquad(w=vL,vR).
\]
By \eqref{sc:eq:split-size} and \eqref{sc:eq:phi-defect},
\[
 m_v\le m_{vL}+m_{vR}+(1+y^{-3})h_v .
\]
Since $m_v\ge(\rho^2-y^{-1/4})h_v$ and $1+2\rho=\rho^2$, we deduce
that
\[
 m_w\ge(\rho-Cy^{-1/4})h_v\quad(w=vL,vR),\qquad
 m_v\le(\rho^2+Cy^{-1/4})h_v .
\]
Thus \ref{it:bal-prop} holds. Let the parent $r$ of $v$ lie in $B$.
Then
\[
 (\rho+O(y^{-1/4}))h_r=m_v=(\rho^2+O(y^{-1/4}))h_v
\]
and therefore $\lvert\log(h_r/(\rho h_v))\rvert\le Cy^{-1/4}$. Thus
\ref{it:bal-parent} holds.

Since $m_v>m_{vL}$, a maximizing prefix $J$ of $I_v$ with
$-\kappa\sum_{i\in J}\log(1-h_i/\kappa)=m_v$ contains $v$ and therefore
contains $I_{vL}$. By \eqref{sc:eq:phi-defect}
and \eqref{eq:B-def},
\[
 \kappa\sum_{i\in I_{vL}}\log(1+h_i^-/\kappa)
 \le (1+h_v/\kappa)\sum_{i\in I_v}h_i^+-m_v\le Cy^{-1/4}h_v .
\]
Since $h_v/\kappa\le y^{-3}$, we have
$\sum_{i\in I_{vL}}\log(1+h_i^-/\kappa)\le Cy^{-13/4}$ and therefore
\[
 \kappa\log(1+h_i^-/\kappa)\ge h_i^-/2\qquad(i\in I_{vL}).
\]
Thus \eqref{sc:eq:balanced-debt} holds.
\end{proof}

\Cref{sc:lem:deviation-integrals} integrates
\eqref{sc:eq:working-quantitative} over $s\le y\le t$ with $t\le4s$.

\begin{lemma}\label{sc:lem:deviation-integrals}
Let $s$ be sufficiently large and $s<t\le4s$ with $e^t t^3\le h_n$.
Then:
\begin{enumerate}
\item\label{it:dev-P} $\mathcal P(y)\le Cn$ for $s\le y\le t$ and
\begin{equation}\label{sc:eq:first-log-integral}
 \int_s^t e^{y/p}\sum_{h_j>e^y}\log\frac{P(y)+m_I}{h_j}\,dy\le Cn .
\end{equation}
\item\label{it:dev-quad}
\begin{equation}\label{sc:eq:deviation-sum}
 \sum_{s<\log h_v\le t}h_v^{1/p}\min\Bigl\{1,\sum_{w=vL,vR}
 \Bigl(\frac1{h_v}\sum_{i\in I_w}h_i^+-\rho\Bigr)^2\Bigr\}\le Cn .
\end{equation}
\item\label{it:dev-gap}
\begin{equation}\label{sc:eq:second-log-integral}
 \int_s^t e^{y/p}\sum_{h_j>e^y}
 \Bigl(\log\frac{P(y)+\sum_{i\in I}h_i^+}{P(y)+m_I}\Bigr)^+dy\le Cn .
\end{equation}
\end{enumerate}
\end{lemma}
\begin{proof}
For an interval $J$ write
\[
 \bar m_J=\sum_{i\in J}h_i^+ .
\]
For a positive step $v$ write $\bar m_v=\bar m_{I_v}$,
$\bar m_{vL}=\bar m_{I_{vL}}$ and $\bar m_{vR}=\bar m_{I_{vR}}$. The
derivation of \eqref{eq:integrated-split} uses no property of $m_I$.
For $\bar m_I$ the derivation of \eqref{eq:integrated-split} gives
\begin{multline}\label{sc:eq:integrated-split-plus}
 \sum_{s<\log h_v\le t}h_v^{1/p}
 \log\frac{h_v(\bar m_v+P(\log h_v))}
          {(\bar m_{vL}+P(\log h_v))(\bar m_{vR}+P(\log h_v))}
 +\frac1p\int_s^t e^{y/p}\sum_{h_j>e^y}\log\frac{\bar m_I+P(y)}{h_j}\,dy\\
 \le \left[e^{y/p}\sum_{h_j>e^y}\log\frac{\bar m_I+P(y)}{h_j}\right]_{y=s}^{y=t}.
\end{multline}

By
\eqref{sc:eq:count-bound} there are at most $ne^{-y/p}$ cuts at
threshold $e^y$. By \eqref{sc:eq:positive-tail},
\[
 \sum_{h_j>e^y}\bar m_I\le Cne^{y(1-1/p)} .
\]
Since $\tau\mapsto\tau^{1/p}$ is concave,
\begin{equation}\label{sc:eq:power-total}
 \sum_{h_j>e^y}\bigl(2\bar m_I+P\bigr)^{1/p}\le Cn .
\end{equation}
By \eqref{sc:eq:mass-product}, $m_I\le2\bar m_I$. By
\eqref{sc:eq:power-total} and $\log\tau\le p\tau^{1/p}$ for $\tau>0$,
\begin{gather*}
 \mathcal P(y)\le Cn,\qquad
 e^{y/p}\sum_{h_j>e^y}\log\frac{m_I+P}{h_j}\le Cn,\\
 \sum_{h_j>e^y}\bar m_I^{1/p}\le Cn,\qquad
 e^{y/p}\sum_{h_j>e^y}\log\frac{\bar m_I+P}{h_j}\le Cn .
\end{gather*}

For every cut $j$,
\[
 \log\frac{P+\bar m_I}{h_j}
 =\log\frac{P+m_I}{h_j}+\log\frac{P+\bar m_I}{P+m_I}.
\]
By \eqref{sc:eq:initial-log-floor} and
$(e^y/h_n)^{1/(2p)}\le y^{-3/(2p)}\le1/y$,
\[
 e^{y/p}\sum_{h_j>e^y}\log\frac{P+m_I}{h_j}\ge-\frac{Cn}{y} .
\]
By \eqref{sc:eq:mass-product} and \eqref{sc:eq:count-bound},
\[
 e^{y/p}\sum_{h_j>e^y}
 \Bigl(\log\frac{P+\bar m_I}{P+m_I}\Bigr)^-
 \le\frac{Cne^y}{\kappa} .
\]
We have $n/y\le Cn$ and $ne^y/\kappa\le Cn$ at $y=s$. Since $t\le4s$
and $e^t\le\kappa$,
\[
 \int_s^t\Bigl(\frac ny+\frac{ne^y}\kappa\Bigr)\,dy\le Cn .
\]
Hence $e^{y/p}\sum_{h_j>e^y}\log[(P+m_I)/h_j]$ is at least $-Cn$ at
$y=s$ and has integral at least $-Cn$ over $s\le y\le t$. Also
$e^{y/p}\sum_{h_j>e^y}\log[(P+\bar m_I)/h_j]$ is at least $-Cn$ at
$y=s$.

Let $v$ be a positive step with $s<\log h_v\le t$. We have
\[
 \bar m_v=h_v+\bar m_{vL}+\bar m_{vR} .
\]
Hence \eqref{sc:eq:split-hypothesis} holds for
$(m_v,m_{vL},m_{vR})=(\bar m_v,\bar m_{vL},\bar m_{vR})$ and
\cref{lem:max-plus-excess}\ref{it:split-real} applied with
$(h_v,P,m_v,m_{vL},m_{vR})=(h_v,P(\log h_v),\bar m_v,\bar m_{vL},\bar m_{vR})$
gives
\begin{multline}\label{sc:eq:split-plus}
 \log\frac{h_v(\bar m_v+P(\log h_v))}
          {(\bar m_{vL}+P(\log h_v))(\bar m_{vR}+P(\log h_v))}
 -\frac p{2-\sqrt2}\,
  \frac{\bar m_v^{1/p}-\bar m_{vL}^{1/p}-\bar m_{vR}^{1/p}}{h_v^{1/p}}\\
 \ge c\min\Bigl\{1,\sum_{w=vL,vR}\Bigl(\frac{\bar m_w}{h_v}-\rho\Bigr)^2\Bigr\}
 -C\Bigl(\frac{P(\log h_v)}{h_v}\Bigr)^{1/p}.
\end{multline}
Take the steps $v$ with $s<\log h_v\le t$ as cuts one at a time in the
order of removal. Each $v$ replaces $\bar m_v^{1/p}$ by
$\bar m_{vL}^{1/p}+\bar m_{vR}^{1/p}$ in the sum
$\sum_{h_j>e^y}\bar m_I^{1/p}$. Hence
\begin{equation}\label{sc:eq:dev-telescope}
 \sum_{s<\log h_v\le t}
 \bigl(\bar m_v^{1/p}-\bar m_{vL}^{1/p}-\bar m_{vR}^{1/p}\bigr)
 =\Bigl[\sum_{h_j>e^y}\bar m_I^{1/p}\Bigr]_{y=s}^{y=t} .
\end{equation}
Since $\sum_{h_j>e^y}\bar m_I^{1/p}\le Cn$ at $y=s$ and at $y=t$, the
right side of \eqref{sc:eq:dev-telescope} is at most $Cn$ in absolute
value. We have $P(\log h_v)/h_v=(\log h_v)^{-3}$. By
\eqref{sc:eq:count-bound} on intervals of $y$ of length one,
\begin{equation}\label{sc:eq:dev-P-sum}
 \sum_{s<\log h_v\le t}h_v^{1/p}\bigl(P(\log h_v)/h_v\bigr)^{1/p}
 \le Cn\int_{s-1}^t y^{-3/p}\,dy\le Cn .
\end{equation}
By \eqref{sc:eq:mass-product}, \eqref{sc:eq:count-bound} and
$\int_s^t(ne^y/\kappa)\,dy\le Cn$, we have
\begin{equation}\label{sc:eq:dev-minus}
 \int_s^t e^{y/p}\sum_{h_j>e^y}
 \Bigl(\log\frac{P(y)+\bar m_I}{P(y)+m_I}\Bigr)^-dy\le Cn .
\end{equation}
The sum $e^{y/p}\sum_{h_j>e^y}\log[(P+\bar m_I)/h_j]$ is at most $Cn$
at $y=t$ and at least $-Cn$ at $y=s$. Hence the right side of
\eqref{sc:eq:integrated-split-plus} is at most $Cn$. On the left side of
\eqref{sc:eq:integrated-split-plus} we bound each logarithm in the sum
over $v$ from below by \eqref{sc:eq:split-plus}. In the integral on
the left side of \eqref{sc:eq:integrated-split-plus} we write
\[
 \log\frac{P(y)+\bar m_I}{h_j}
 =\log\frac{P(y)+m_I}{h_j}
 +\Bigl(\log\frac{P(y)+\bar m_I}{P(y)+m_I}\Bigr)^+
 -\Bigl(\log\frac{P(y)+\bar m_I}{P(y)+m_I}\Bigr)^- .
\]
By \eqref{sc:eq:dev-telescope}, \eqref{sc:eq:dev-P-sum} and
\eqref{sc:eq:dev-minus} we obtain
\begin{multline}\label{sc:eq:dev-combined}
 c\sum_{s<\log h_v\le t}h_v^{1/p}\min\Bigl\{1,\sum_{w=vL,vR}
 \Bigl(\frac1{h_v}\sum_{i\in I_w}h_i^+-\rho\Bigr)^2\Bigr\}\\
 +\frac1p\int_s^t e^{y/p}\sum_{h_j>e^y}\log\frac{P(y)+m_I}{h_j}\,dy
 +\frac1p\int_s^t e^{y/p}\sum_{h_j>e^y}
 \Bigl(\log\frac{P(y)+\sum_{i\in I}h_i^+}{P(y)+m_I}\Bigr)^+dy\le Cn .
\end{multline}
The sum over $v$ and the integral of the positive part in
\eqref{sc:eq:dev-combined} are nonnegative. The integral
$\frac1p\int_s^t e^{y/p}\sum_{h_j>e^y}\log[(P(y)+m_I)/h_j]\,dy$ is at
least $-Cn$. Hence each of the three terms on the left side of
\eqref{sc:eq:dev-combined} is at most $Cn$. Thus
\eqref{sc:eq:first-log-integral}, \eqref{sc:eq:deviation-sum} and
\eqref{sc:eq:second-log-integral} hold.
\end{proof}

\Cref{sc:lem:excess-sum} follows from
\cref{sc:lem:deviation-integrals} and
\cref{sc:prop:working-realization}.

\begin{lemma}\label{sc:lem:excess-sum}\label{sc:lem:absolute}
Let $s$ be sufficiently large and $s<t\le4s$ with $e^t t^3\le h_n$.
Then:
\begin{enumerate}
\item\label{it:exc-var}
\begin{equation}\label{sc:eq:working-excess}
 \sum_{s<\log h_v\le t}\lvert\Delta_v\rvert\le Cn .
\end{equation}
\item\label{it:exc-unbal}
\begin{equation}\label{sc:eq:unbalanced-count}
 \sum_{v\in B'}\min\{h_v,e^t\}^{1/p}\le Cn\sqrt t .
\end{equation}
\end{enumerate}
\end{lemma}
\begin{proof}
For $s<\log h_v\le t$ we
have $h_v\le e^t\le h_n\le\kappa/4$. By \eqref{sc:eq:split-size},
\eqref{sc:eq:phi-defect} and subadditivity of
$\tau\mapsto\tau^{1/p}$,
\[
 \Delta_v^+\le Ch_v^{1/p}.
\]
By \eqref{sc:eq:count-bound}, $\sum h_v^{1/p}\le Cn$ over the steps
$v$ with $t-1<\log h_v\le t$. For each positive step $u$ define
\[
 \delta_u=\min\Bigl\{1,\sum_{w=uL,uR}
 \Bigl(\frac1{h_u}\sum_{i\in I_w}h_i^+-\rho\Bigr)^2\Bigr\}.
\]
Let $s<\log h_v\le t-1$. Since $h_v<h_n$, the step $v$ has a parent
$r$.

Suppose that $h_r<2h_v$. Let
$\tau=\max\{\delta_v,\delta_r\}\in[0,1]$. If $\tau<1$, then
$\sum_{i\in I_{vL}}h_i^++\sum_{i\in I_{vR}}h_i^+\ge(2\rho-\sqrt{2\tau})h_v$
and $\sum_{i\in I_v}h_i^+\le(\rho+\sqrt\tau)h_r$. Since
$\sum_{i\in I_v}h_i^+=h_v+\sum_{i\in I_{vL}}h_i^++\sum_{i\in I_{vR}}h_i^+$,
$\rho^2=2\rho+1$ and $h_r<2h_v$, we obtain
\[
 1\le\rho^2-\rho\frac{h_r}{h_v}\le\bigl(\sqrt2+2\bigr)\sqrt\tau
\]
and hence $\tau\ge c$. For $\tau=1$ the bound $\tau\ge c$ holds with
$c\le1$. Thus $\max\{\delta_v,\delta_r\}\ge c$. We have
$\log h_r\le t$ and $r$ has at most two children. Hence
\eqref{sc:eq:deviation-sum} gives $\sum h_v^{1/p}\le Cn$ over the steps
$v$ with $s<\log h_v\le t-1$ and $h_r<2h_v$. Thus
\[
 \sum\Delta_v^+\le Cn
\]
over the steps $v$ with $t-1<\log h_v\le t$ and over the steps $v$
with $s<\log h_v\le t-1$ and $h_r<2h_v$.

Suppose that $h_r\ge2h_v$. For $h_v<e^y<2h_v$ the step $r$ is a cut at
threshold $e^y$ and $v$ is not a cut. Hence $I_v$ is the interval $I$
before a cut at threshold $e^y$. For $s\le y\le t$ apply
\cref{sc:prop:working-realization} at threshold $e^y$ with $C_0=1$ and
$S=\varnothing$. Let $V$ be the set of steps $v$ with $h_v\ge e^y/2$
and $\Delta_v>0$ such that $v$ is the first step removed from the
interval $I$ before a cut $j$. By \eqref{sc:eq:bulk-split} with
$P=P(y)$,
\[
 \sum_{v\in V}\log\frac{h_v(m_v+P(y))}{(m_{vL}+P(y))(m_{vR}+P(y))}
 \ge\sum_{v\in V}\bigl(\Delta_v/h_v^{1/p}-Cy^{-2/p}\bigr).
\]
By \eqref{sc:eq:count-bound}, $\#\{j:h_j>e^y\}+\#V\le2ne^{-y/p}$. Also
$y^{-2/p}\le1/y$ and $(e^y/h_n)^{1/(2p)}\le y^{-3/(2p)}\le1/y$.
Multiply \eqref{sc:eq:working-quantitative} by $e^{y/p}$ and integrate
over $s\le y\le t$. By \eqref{sc:eq:first-log-integral} and
$\int_s^t dy/y\le\log 4$,
\[
 \int_s^t e^{y/p}\sum_{v\in V}\frac{\Delta_v}{h_v^{1/p}}\,dy\le Cn .
\]
Every $v$ with $h_r\ge2h_v$, $\Delta_v>0$ and $s<\log h_v\le t-1$
belongs to $V$ for $\log h_v<y<\log(2h_v)$ and
\[
 \int_{\log h_v}^{\log(2h_v)}e^{y/p}\frac{\Delta_v}{h_v^{1/p}}\,dy
 \ge(\log2)\Delta_v .
\]
We obtain $\sum\Delta_v^+\le Cn$ over all $v$ with $s<\log h_v\le t$.
By \eqref{sc:eq:band-sums} and
\cref{sc:lem:deviation-integrals}\ref{it:dev-P},
\[
 \sum_{s<\log h_v\le t}\lvert\Delta_v\rvert
 =2\sum_{s<\log h_v\le t}\Delta_v^+-\mathcal P(t)+\mathcal P(s)\le Cn .
\]

Thus \eqref{sc:eq:working-excess} holds. Let $v\notin B$ be a
positive step with $s<\log h_v\le t-1$ and $h_r\ge2h_v$ for the parent
$r$ of $v$. Suppose first that
\[
 \Bigl\lvert\frac1{h_v}\sum_{i\in I_w}h_i^+-\rho\Bigr\rvert
 \ge\frac{t^{-1/4}}4
\]
for $w=vL$ or for $w=vR$. Then $h_v^{1/p}\delta_v\ge ch_v^{1/p}/\sqrt t$.
Suppose now that
\[
 \Bigl\lvert\frac1{h_v}\sum_{i\in I_w}h_i^+-\rho\Bigr\rvert
 <\frac{t^{-1/4}}4
\]
for $w=vL$ and for $w=vR$. Since $\log h_v>s\ge y_0$,
$h_v(\log h_v)^3\le e^tt^3\le h_n$ and
$t^{-1/4}/4<(\log h_v)^{-1/4}$, the assumption $v\notin B$ forces
$m_v<(\rho^2-(\log h_v)^{-1/4})h_v\le(\rho^2-t^{-1/4})h_v$. Hence
\[
 \sum_{i\in I_v}h_i^+\ge m_v+h_vt^{-1/4}/2 .
\]
For $h_v<e^y<2h_v$ we have $I_v=I$ for a cut $j$ and $P(y)\le2h_v$.
Therefore
\[
 \log\frac{P(y)+\sum_{i\in I_v}h_i^+}{P(y)+m_v}\ge ct^{-1/4}
 \qquad(h_v<e^y<2h_v).
\]
Thus
\[
 \int_{\log h_v}^{\log(2h_v)}e^{y/p}
 \Bigl(\log\frac{P(y)+\sum_{i\in I_v}h_i^+}{P(y)+m_v}\Bigr)^+dy
 \ge ch_v^{1/p}t^{-1/4}\ge\frac{ch_v^{1/p}}{\sqrt t}.
\]
By \eqref{sc:eq:deviation-sum}, \eqref{sc:eq:second-log-integral} and
the bounds $\sum h_v^{1/p}\le Cn$ over the steps $v$ with
$t-1<\log h_v\le t$ and over the steps $v$ with $s<\log h_v\le t-1$
and $h_r<2h_v$,
\begin{equation}\label{eq:nonB-range}
 \sum_{\substack{v\notin B\\ s<\log h_v\le t}}h_v^{1/p}\le Cn\sqrt t .
\end{equation}
The sum in \eqref{eq:nonB-range} is over positive steps.

Fix an absolute cutoff $y_1\ge y_0$ sufficiently large for
\eqref{eq:nonB-range} to hold whenever $s\ge y_1$. Let $\nu\ge1$ be the smallest integer with $2^{-\nu}t\le2y_1$. Apply
\eqref{eq:nonB-range} to the dyadic ranges
$(2^{-i}t,2^{1-i}t]$ for $1\le i\le\nu$. The positive steps $v$ with
$\log h_v\le2^{-\nu}t$ satisfy $\log h_v\le2y_1$. Hence
$\sum_{\log h_v\le2y_1}h_v^{1/p}\le e^{2y_1/p}\#\{i:h_i>0\}\le Cn$.
Summing over the dyadic ranges gives $\sum h_v^{1/p}\le Cn\sqrt t$
over the positive steps $v\notin B$ with $h_v\le e^t$.

Let $r\in B$ with $h_r\le e^t$ have a child $b\notin B$. By
\cref{sc:lem:balance}\ref{it:bal-prop} and \eqref{sc:eq:mass-product},
\[
 \sum_{i\in I_b}h_i^+\ge ch_r .
\]
Let $T_b$ be the smallest set that contains $b$ and every child outside
$B$ of a step of $T_b$. A positive step of $I_b$ outside $T_b$ lies in
$I_u$ for a child $u\in B$ of a step $v\in T_b$. By \eqref{eq:B-def},
\[
 \sum_{i\in I_u}h_i^+\le Ch_u\le Ch_v .
\]
Every step has at most two children. Hence
\[
 h_r\le C\sum_{v\in T_b}h_v,\qquad h_r^{1/p}\le C\sum_{v\in T_b}h_v^{1/p}.
\]
Let $(r,b)$ and $(r',b')$ be distinct pairs. Then $r$ is the
nearest ancestor in $B$ of every step of $T_b$. Hence $T_b$ and
$T_{b'}$ are disjoint for $r\ne r'$. For $r=r'$ the sets
$T_b\subset I_b$ and $T_{b'}\subset I_{b'}$ are disjoint. Hence
$\sum h_r^{1/p}\le Cn\sqrt t$ over the steps $r\in B$ with
$h_r\le e^t$ and with a child not in $B$. By \eqref{sc:eq:count-bound},
the sum of $\min\{h_v,e^t\}^{1/p}$ over the steps $v$ with $h_v>e^t$ is
$e^{t/p}\#\{i:h_i>e^t\}\le n$. Thus \eqref{sc:eq:unbalanced-count}
holds.
\end{proof}

\section{Distance estimates for consecutive splits}\label{sec:consecutive-splits}
For sufficiently large $s$ and $s<t\le4s$ with $e^tt^3\le h_n$,
\cref{sc:lem:excess-sum} gives
$\sum_{s<\log h_v\le t}\lvert\Delta_v\rvert\le Cn$. The bound $Cn$
comes from the term $C/y$ for each cut in
\eqref{sc:eq:working-quantitative}. The number of ranges $(s,4s]$
between a bounded threshold and $\log h_n$ is unbounded. Hence
\eqref{sc:eq:band-sums} and \cref{sc:lem:excess-sum} do not give
$\mathcal P(y)\ge n/2$ at a bounded $y$. \Cref{sc:lem:four-phase}
constructs one solution on four consecutive intervals. For the
solution of \cref{sc:lem:four-phase} the term $O(y^{-2})$ replaces the
term $C/y$. For sufficiently large $t$ with $e^tt^3\le h_n$ and
$t/2\le y\le t$ let $N$ be the smallest integer with
$2^N\ge(\log t)^2$. \Cref{sc:lem:two-ranges} applies
\cref{sc:prop:working-realization} with a set $S$ of cuts $j$ at
threshold $e^y$ for which steps in $B$ divide $I$ into $2^N$ intervals.
Let $I_\omega$ be the interval with the largest indices among the
$2^N$ intervals. For every interval $J\ne I_\omega$ among the $2^N$
intervals, \cref{sc:lem:two-ranges} uses \cref{sc:lem:four-phase} on
$J$ and the step after $J$. The term $C/y$ remains for $I_\omega$ and
the cut $j$. For each $j\in S$ the terms $2^NO(y^{-2})$ of
\cref{sc:lem:four-phase} and the term $C/y$ of $I_\omega$ sum to at
most $C/y$. \Cref{sc:sec:proof-count} proves
\cref{sc:thm:count}.
\subsection{A trajectory on four consecutive intervals}\label{sc:sec:local-tests}\label{sc:sec:four-phase}

Let $y$ satisfy \eqref{sc:eq:band-B}. Let $a$ be a positive step with
$h_a\le e^y$ and children $aL,aR$ and write $I_a=\{k+1,\ldots,j-1\}$.
Here $j$ need not be a cut at threshold $e^y$. By
\eqref{sc:eq:ordinary-gain}, a solution given by
\cref{sc:lem:ordinary} on $(k,j]$ satisfies
\[
 \frac{\norm{x_j}}{\norm{x_k}}\ge\frac{h_j}{m_a+P}
 \exp\bigl(-C/\log A-m_a/\kappa\bigr).
\]
For solutions given by \cref{sc:lem:ordinary} on $(k,a]$ and on
$(a,j]$ the product of the two factors $h_j/(m_I+P)$ of
\eqref{sc:eq:ordinary-gain} is
\[
 \frac{h_ah_j}{(m_{aL}+P)(m_{aR}+P)} .
\]
For solutions given by \cref{sc:lem:ordinary} on $(k,aL]$, $(aL,a]$,
$(a,aR]$ and $(aR,j]$ the product of the four factors $h_j/(m_I+P)$ of
\eqref{sc:eq:ordinary-gain} is
\[
 \frac{h_{aL}h_ah_{aR}h_j}{(m_{aLL}+P)(m_{aLR}+P)(m_{aRL}+P)(m_{aRR}+P)}.
\]
The factor $\exp(-C/\log A)$ in \eqref{sc:eq:ordinary-gain} comes from
the bound $1+C/\log A$ of \cref{sc:lem:arc}\ref{sc:it:arc-size} and
gives the term $C/y$ in \eqref{sc:eq:working-quantitative}. For
$a,aL,aR,j\in B$ with $h_a\le h_j\le e^y$ and
$1\le e^y/h_a\le y^{1/12}$, \cref{sc:lem:four-phase} gives a solution
$\mathcal S$ on $(k,j]$ without \cref{sc:lem:arc}. The vectors $g_u$ of
$\mathcal S$ are equal for $k\le u<aL$, for $aL\le u<a$, for
$a\le u<aR$ and for $aR\le u<j$. By \eqref{sc:eq:four-gain},
$\log(\norm{x_j}/\norm{x_k})$ is a weighted mean of
$\log Q_0(\mu_7)$, $\log Q_1(\mu_7)$ and $\log Q_2(\mu_7)$ up to
$O(P/h_a)=O(y^{-2})$. Here $\mu_7$ is defined in
\eqref{sc:eq:four-mu} and $Q_0,Q_1,Q_2$ are defined in
\eqref{sc:eq:high-gain-score}.

Define
\[
 (k_0,j_0)=(k,aL),\qquad(k_1,j_1)=(aL,a),\qquad
 (k_2,j_2)=(a,aR),\qquad(k_3,j_3)=(aR,j).
\]
Every step of $I_a$ has stepsize at most $h_a\le e^y<\kappa$. For
$i=0,1,2,3$ choose $k_i\le u_i<j_i$ such that the prefix
$\{k_i+1,\ldots,u_i\}$ attains the maximum in \eqref{sc:eq:mass-def}
for the interval $\{k_i+1,\ldots,j_i-1\}$. Define
\begin{equation}\label{sc:eq:four-mu}
 \mu_{2i}=m_{\{k_i+1,\ldots,u_i\}}\quad(0\le i\le3),\qquad
 \mu_{2i+1}=m_{(u_i,j_i]}\quad(0\le i\le2),
\end{equation}
and $\mu_7=m_{(u_3,j]}$ when $h_j<\kappa$. Thus
$\mu_0,\mu_2,\mu_4,\mu_6$ equal $m_{aLL},m_{aLR},m_{aRL},m_{aRR}$. If
every stepsize in $\{k_i+1,\ldots,j_i-1\}$ is positive and
$h_{j_i}<\kappa$, then $u_i=j_i-1$ and
$\mu_{2i+1}=-\kappa\log(1-h_{j_i}/\kappa)$. For $b\in B$
with $e^y/y\le h_b\le e^y$, \cref{sc:lem:balance}\ref{it:bal-prop},
\eqref{sc:eq:split-size} and \eqref{sc:eq:phi-defect} give
\begin{align}
 m_b-m_{bL}-m_{bR}&=(1+O(y^{-1/4}))h_b,\label{sc:eq:balanced-unit-K}\\
 m_b-m_{bL}-m_{bR}&\le h_b+h_b^2/\kappa .\label{sc:eq:balanced-unit-upper}
\end{align}
Let $a,aL,aR\in B$. By \cref{sc:lem:balance}\ref{it:bal-prop} we have
$m_a,m_{aL},m_{aR},\mu_0,\mu_2,\mu_4,\mu_6>0$. Define
for $\tau>0$
\begin{equation}\label{sc:eq:high-gain-score}
 Q_0(\tau)=\frac{\tau}{m_a},\qquad
 Q_1(\tau)=\frac{\tau\mu_3}{m_{aL}m_{aR}},\qquad
 Q_2(\tau)=\frac{\mu_1\mu_3\mu_5\tau}{\mu_0\mu_2\mu_4\mu_6}.
\end{equation}

\begin{figure}[htbp]
\centering
\begin{tikzpicture}[font=\footnotesize,x=1cm,y=1cm]
 \draw[thick] (0,0)--(14.4,0);
 \foreach \x/\lab in {0.4/k, 3.8/aL, 7.2/a, 10.6/aR, 14.0/j} {
  \draw[very thick] (\x,-0.25)--(\x,0.25);
  \node[above=5pt] at (\x,0.18) {$\lab$};
 }
 \foreach \x/\lab in {2.25/u_0, 5.65/u_1, 9.05/u_2, 12.45/u_3} {
  \draw[dashed] (\x,-0.2)--(\x,0.2);
  \fill (\x,0) circle (1.2pt);
  \node[above=2pt] at (\x,0.16) {$\lab$};
 }
 \node[below=6pt,font=\scriptsize] at (1.32,0) {$\mu_0\approx h_a$};
 \node[below=6pt,font=\scriptsize] at (3.02,0) {$\mu_1\approx h_a/\rho$};
 \node[below=6pt,font=\scriptsize] at (4.72,0) {$\mu_2\approx h_a$};
 \node[below=6pt,font=\scriptsize] at (6.42,0) {$\mu_3\approx h_a$};
 \node[below=6pt,font=\scriptsize] at (8.12,0) {$\mu_4\approx h_a$};
 \node[below=6pt,font=\scriptsize] at (9.82,0) {$\mu_5\approx h_a/\rho$};
 \node[below=6pt,font=\scriptsize] at (11.52,0) {$\mu_6\approx h_a$};
 \node[below=6pt,font=\scriptsize] at (13.22,0) {$\mu_7$};
 \draw (0.5,-0.82)--(0.5,-0.92)--(3.7,-0.92)--(3.7,-0.82);
 \node[below] at (2.1,-0.92) {$(k_0,j_0]=(k,aL]$};
 \draw (3.9,-0.82)--(3.9,-0.92)--(7.1,-0.92)--(7.1,-0.82);
 \node[below] at (5.5,-0.92) {$(k_1,j_1]=(aL,a]$};
 \draw (7.3,-0.82)--(7.3,-0.92)--(10.5,-0.92)--(10.5,-0.82);
 \node[below] at (8.9,-0.92) {$(k_2,j_2]=(a,aR]$};
 \draw (10.7,-0.82)--(10.7,-0.92)--(13.9,-0.92)--(13.9,-0.82);
 \node[below] at (12.3,-0.92) {$(k_3,j_3]=(aR,j]$};
 \draw (0.5,0.85)--(0.5,0.95)--(7.1,0.95)--(7.1,0.85);
 \node[above] at (3.8,0.95) {$I_{aL},\; m_{aL}$};
 \draw (7.3,0.85)--(7.3,0.95)--(13.9,0.95)--(13.9,0.85);
 \node[above] at (10.6,0.95) {$I_{aR},\; m_{aR}$};
 \draw (0.5,1.5)--(0.5,1.6)--(13.9,1.6)--(13.9,1.5);
 \node[above] at (7.2,1.6) {$I_a,\; m_a$};
\end{tikzpicture}
\caption{The sets $(k_i,j_i]$ and the indices $u_i$ for $i=0,1,2,3$
with $\mu_0,\ldots,\mu_7$ of \eqref{sc:eq:four-mu}. The approximate values of $\mu_0,\ldots,\mu_6$ in the figure hold
for $a,aL,aR\in B$ by \eqref{sc:eq:balance}.}
\label{fig:four-intervals}
\end{figure}
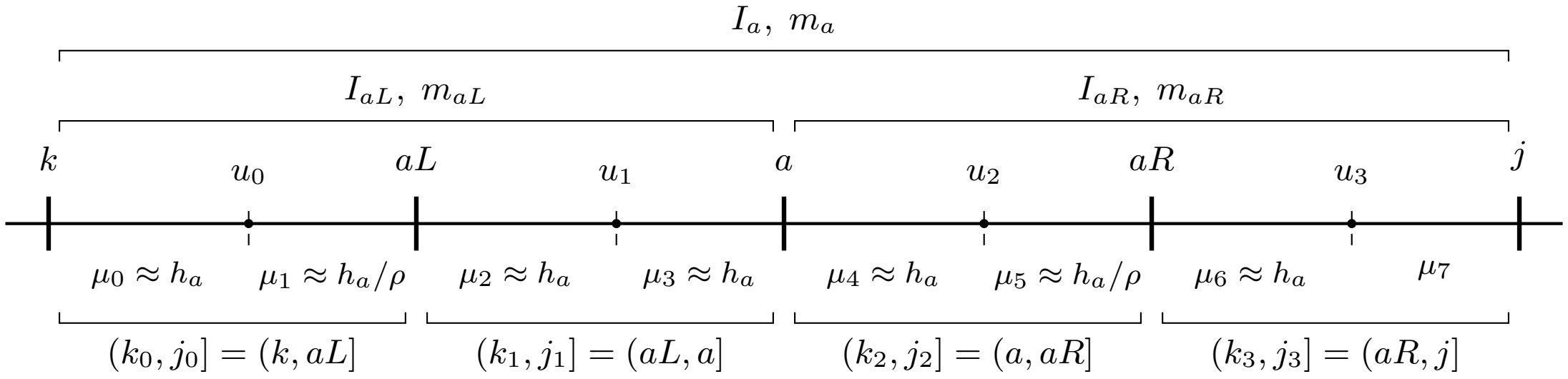

\begin{lemma}\label{sc:lem:four-masses}
Let $y$ be sufficiently large with \eqref{sc:eq:band-B}. Let
$a,aL,aR\in B$ and $1\le e^y/h_a\le y^{1/12}$. For $i=0,1,2$ write
$j_iL$ and $j_iR$ for the children of $j_i$. Then
\begin{equation}\label{sc:eq:four-masses}
 \mu_{2i+1}=-\kappa\sum_{u_i<w\le j_i}\log(1-h_w/\kappa)
 =m_{j_i}-m_{j_iL}-m_{j_iR}
\end{equation}
for $i=0,1,2$. Let also $j\in B$ and $h_a\le h_j\le e^y$. Write $j_3L$
and $j_3R$ for the children of $j$. Then
\eqref{sc:eq:four-masses} holds for $i=3$.
\end{lemma}
\begin{proof}
Let $i\in\{0,1,2\}$ or let $i=3$, $j\in B$ and $h_a\le h_j\le e^y$.
Write $b=j_i$ and $J=\{k_i+1,\ldots,j_i-1\}$. Then $J$ is the
interval of a child of $aL$ or of $aR$. Hence
$\mu_{2i}=(1+O(y^{-1/4}))h_a$ by
\cref{sc:lem:balance}\ref{it:bal-prop} and \ref{it:bal-parent}. For
$i=0$ we have $J=I_{aLL}=I_{bL}$. For $i=2$ we have
$J=I_{aRL}=I_{bL}$. For $i=1$ we have
$J=I_{aLR}\subset I_{aL}=I_{bL}$. For $i=3$ we have $I_a\subset I_{jL}$ since $j$ is removed before
every step of $I_a$. Hence $J=I_{aRR}\subset I_{jL}=I_{bL}$. For every
$i$ we have
$J=\{w\in I_{bL}:w>k_i\}$. Since
$b\in B$ and $h_b\le e^y$, \eqref{sc:eq:balanced-debt} and
$e^y/h_a\le y^{1/12}$ give
\[
 \sum_{w\in I_{bL}}h_w^-\le Cy^{-1/4}e^y\le Cy^{-1/6}h_a .
\]
Hence the sum $-\kappa\sum\log(1-h_w/\kappa)$ over a prefix of $I_{bL}$
contained in $\{w\in I_{bL}:w\le k_i\}$ exceeds the sum over
$\{w\in I_{bL}:w\le k_i\}$ by at most $Cy^{-1/6}h_a$. By the choice of
$u_i$, the sum $-\kappa\sum_{w=k_i+1}^{u_i}\log(1-h_w/\kappa)$ equals
$\mu_{2i}=(1+O(y^{-1/4}))h_a$ and is the maximum of
$-\kappa\sum_{w=k_i+1}^u\log(1-h_w/\kappa)$ over $k_i\le u<b$. Thus the
prefix $\{w\in I_{bL}:w\le u_i\}$ attains $m_{bL}$. Since $m_b>m_{bL}$
by \cref{sc:lem:balance}\ref{it:bal-prop}, \eqref{sc:eq:split-m} shows
that $-\kappa\sum\log(1-h_w/\kappa)$ over $I_{bL}\cup\{b\}$ equals
$m_b-m_{bR}$. Hence
\[
 -\kappa\sum_{u_i<w\le b}\log(1-h_w/\kappa)=m_b-m_{bL}-m_{bR}.
\]
For $u_i<u<b$ the sum $-\kappa\sum_{u_i<w\le u}\log(1-h_w/\kappa)$ is
at most zero by the choice of $u_i$. Also $m_b-m_{bL}-m_{bR}>0$ by
\eqref{sc:eq:balanced-unit-K}. Hence $m_{(u_i,b]}=m_b-m_{bL}-m_{bR}$.
Thus \eqref{sc:eq:four-masses} holds.
\end{proof}

Let $a,aL,aR\in B$ and $1\le e^y/h_a\le y^{1/12}$ with $y$ sufficiently
large. By \eqref{sc:eq:four-masses} for $i=0,1,2$,
\begin{equation}\label{sc:eq:four-sums}
 m_{aL}=\mu_0+\mu_1+\mu_2,\qquad m_{aR}=\mu_4+\mu_5+\mu_6,\qquad
 m_a=\mu_0+\cdots+\mu_6 .
\end{equation}
Hence for every $\tau>0$,
\begin{equation}\label{sc:eq:four-ratios}
\begin{aligned}
 \frac{Q_1(\tau)}{Q_0(\tau)}
 &=\frac{(m_a-m_{aL}-m_{aR})m_a}{m_{aL}m_{aR}},\\
 \frac{Q_2(\tau)}{Q_1(\tau)}
 &=\frac{(m_{aL}-m_{aLL}-m_{aLR})m_{aL}}{m_{aLL}m_{aLR}}
  \cdot\frac{(m_{aR}-m_{aRL}-m_{aRR})m_{aR}}{m_{aRL}m_{aRR}} .
\end{aligned}
\end{equation}
By \eqref{sc:eq:four-masses}, \eqref{sc:eq:balanced-unit-K} and
\cref{sc:lem:balance}\ref{it:bal-prop} and \ref{it:bal-parent},
\begin{equation}\label{sc:eq:balance}
 \frac1{h_a}(\mu_0,\ldots,\mu_6)
 =(1,1/\rho,1,1,1,1/\rho,1)+O(y^{-1/4}).
\end{equation}

Let $a$ be a positive step and $c_0,C_0>0$. We say that a solution
$\mathcal S'$ with first point $x'$ satisfies
\eqref{sc:eq:right-piece} for $a$ with constants $c_0,C_0$ if
\begin{equation}\label{sc:eq:right-piece}
 c_0h_a\le\norm{x'}\le C_0h_a,\qquad \norm{g}\le C_0,\qquad
 \ip{x'}g\ge0,\qquad
 0\le\ip{x'}z\le (\rho^2+2)h_a q+P^2
\end{equation}
for every triple $(z,g,q)$ of $\mathcal S'$. The constant $\rho^2+2$
in \eqref{sc:eq:right-piece} exceeds $\rho^2+1$ in
\eqref{sc:eq:four-boundary} and satisfies $(\rho^2+2)/(2\rho^2)<1$.

\begin{lemma}\label{sc:lem:four-phase}\label{lem:high-planes}
Let $y$ be sufficiently large with \eqref{sc:eq:band-B} and let $A$ and
$P$ be given by \eqref{eq:cut-scale}. Let $a$ be a positive step, let
$k$ be the index before $I_a$ and let $j$ be the index after $I_a$.
Assume $a,aL,aR,j\in B$, $1\le e^y/h_a\le y^{1/12}$ and
$h_a\le h_j\le e^y$. There are absolute constants
$c_{\mathrm{four}},C_{\mathrm{four}}>0$ and a four-dimensional solution
$\mathcal S$ on $(k,j]$ with nonzero first point $x_k$ and nonzero
last point $x_j\perp x_k$ such that:
\begin{enumerate}
\item\label{it:four-gain} there are numbers
$\alpha_0,\alpha_1,\alpha_2\ge0$ with $\alpha_0+\alpha_1+\alpha_2=1$,
$\alpha_1=1/2$ and $\alpha_2=1/(2\sqrt2)+O(y^{-1/4})$ such that
\begin{equation}\label{sc:eq:four-gain}
 \log\frac{\norm{x_j}}{\norm{x_k}}
 =\alpha_0\log Q_0(\mu_7)+\alpha_1\log Q_1(\mu_7)+\alpha_2\log Q_2(\mu_7)+O(P/h_a) .
\end{equation}
\item\label{it:four-boundary} $\mathcal S$ satisfies
\eqref{sc:eq:start-projections} with constant $C_{\mathrm{four}}$ and
every triple $(z_u,g_u,q_u)$ of $\mathcal S$ satisfies
\begin{equation}\label{sc:eq:four-boundary}
\begin{gathered}
 c_{\mathrm{four}}h_a\le\norm{x_k}\le C_{\mathrm{four}}h_a,\qquad
 \norm{g_u}\le C_{\mathrm{four}},\\
 \ip{x_k}{g_u}\ge0,\qquad
 0\le\ip{x_k}{z_u}\le(\rho^2+1)h_a q_u+P^2 .
\end{gathered}
\end{equation}
\item\label{it:four-join} Let $c_0,C_0>0$ and let $y$ be sufficiently
large depending on $c_0$ and $C_0$. Let $\tilde{\mathcal S}$ be the
image of $\mathcal S$ under a rotation and \eqref{eq:scaling}. Let
$j<b'\le n$ and let $\mathcal S'$ be a solution on $(j,b']$ that
satisfies \eqref{sc:eq:right-piece} for $a$ with constants $c_0,C_0$.
Let $\tilde{\mathcal S}'$ be the image of $\mathcal S'$ under a
rotation and \eqref{eq:scaling} such that the first point of
$\tilde{\mathcal S}'$ is the last point of $\tilde{\mathcal S}$. If
$\tilde{\mathcal S}$ and $\tilde{\mathcal S}'$ satisfy
\cref{sc:lem:assembly}\ref{it:asm-shared}--\ref{it:asm-perp}, then
\eqref{eq:tail-planes} holds for every triple of $\tilde{\mathcal S}$
together with every triple of $\tilde{\mathcal S}'$.
\end{enumerate}
\end{lemma}

From here on a solution given by \cref{sc:lem:four-phase} for $a$
denotes a four-dimensional solution $\mathcal S$ on $(k,j]$ with
nonzero first point $x_k$, nonzero last point $x_j\perp x_k$ and the
properties \ref{it:four-gain}--\ref{it:four-join} of
\cref{sc:lem:four-phase}.

\begin{proof}
In this proof $o(\cdot)$ refers to $y\to\infty$. Since $j$ is removed
before every step of $I_a$, we have $I_a\subset I_{jL}$. By
\cref{sc:lem:balance}\ref{it:bal-prop}, \eqref{sc:eq:mass-product} and
\eqref{eq:B-def},
\[
 (\rho^2-o(1))h_a\le m_a\le(1+h_a/\kappa)\sum_{w\in I_{jL}}h_w^+
 \le(\rho+o(1))h_j .
\]
Hence $h_j\ge(\rho-o(1))h_a$. By \eqref{sc:eq:four-masses} for $i=3$
and \eqref{sc:eq:balanced-unit-K}, $\mu_7\ge(\rho-o(1))h_a$.

Let $\bar g_0,\bar g_1,\bar g_2,\bar g_3$ be vectors. Take
$g_u=\bar g_i$ for $i=0,1,2,3$ and $k_i\le u<j_i$. For
$k_i\le u\le j_i$,
\eqref{sc:eq:tail-update} gives
\begin{equation}\label{sc:eq:four-constant}
 x_u+(\kappa-1)\bar g_i
 =\left[\prod_{w=k_i+1}^u(1-h_w/\kappa)\right]\bigl(x_{k_i}+(\kappa-1)\bar g_i\bigr).
\end{equation}
By the choice of $u_i$, the product over $k_i<w\le u_i$ is
$e^{-\mu_{2i}/\kappa}$. By \eqref{sc:eq:four-masses}, the product over
$u_i<w\le j_i$ is $e^{-\mu_{2i+1}/\kappa}$. Since
$x_{u_i}=z_{u_i}+\bar g_i$,
\begin{equation}\label{sc:eq:four-relations}
\begin{aligned}
 x_{k_i} &= e^{\mu_{2i}/\kappa}z_{u_i}+[\kappa(e^{\mu_{2i}/\kappa}-1)+1]\bar g_i,\\
 z_{u_i} &= e^{\mu_{2i+1}/\kappa}x_{j_i}+
       [(\kappa-1)(e^{\mu_{2i+1}/\kappa}-1)-1]\bar g_i .
\end{aligned}
\end{equation}
Put $z_{u_3}=P\bar g_3$. Then
$x_j=[(\kappa+P)e^{-\mu_7/\kappa}-(\kappa-1)]\bar g_3=-(\mu_7+O(P))\bar g_3$.

Put $E_{-1}=1$ and $E_r=\exp[(\mu_0+\cdots+\mu_r)/\kappa]$ for
$0\le r\le6$ and write $\bar z_i=E_{2i}z_{u_i}$ for $i=0,1,2$.
Multiplying $z_{u_i}$ by a positive number does not change the order of
$\ip g{z_{u_i}}$ for different $g$. Define
\begin{align*}
 \tilde\mu_{2i}
 &=E_{2i-1}\bigl[\kappa(e^{\mu_{2i}/\kappa}-1)+1\bigr]
  +\mathbf 1_{\{i=3\}}E_6P &&(0\le i\le3),\\
 \tilde\mu_{2i+1}
 &=E_{2i}\bigl[(\kappa-1)(e^{\mu_{2i+1}/\kappa}-1)-1\bigr] &&(0\le i\le2).
\end{align*}
By \eqref{sc:eq:four-relations} for $i=0,1,2,3$ and
$z_{u_3}=P\bar g_3$,
\begin{equation}\label{sc:eq:four-coefficients}
\begin{aligned}
 x_k-\bar z_0 &= \tilde\mu_0\bar g_0, &
 \bar z_0-\bar z_1 &= \tilde\mu_1\bar g_0+\tilde\mu_2\bar g_1,\\
 \bar z_1-\bar z_2 &= \tilde\mu_3\bar g_1+\tilde\mu_4\bar g_2, &
 \bar z_2 &= \tilde\mu_5\bar g_2+\tilde\mu_6\bar g_3 .
\end{aligned}
\end{equation}
Since $\mu_r\ge ch_a\to\infty$ and $\mu_r=O(e^y)$ for $r=0,\ldots,6$,
we have $\tilde\mu_0,\ldots,\tilde\mu_6>0$ for sufficiently large
$y$. We have
$\kappa(e^{\tau/\kappa}-1)=\tau+O(\tau^2/\kappa)$ for
$0\le\tau\le Ce^y$, $E_r=1+O(e^y/\kappa)$ and $e^{2y}/\kappa\le P$ by
\eqref{eq:cut-scale}. Hence
$\tilde\mu_r=\mu_r+O(1+P+e^{2y}/\kappa)=\mu_r+O(P)$ for
$r=0,\ldots,6$. By \eqref{sc:eq:balance} and $P/h_a\le y^{-2}$,
\begin{equation}\label{sc:eq:tilde-balance}
 \frac1{h_a}(\tilde\mu_0,\ldots,\tilde\mu_6)
 =(1,1/\rho,1,1,1,1/\rho,1)+O(y^{-1/4}).
\end{equation}
Write $\tilde\mu_{aL}=\tilde\mu_0+\tilde\mu_1+\tilde\mu_2$,
$\tilde\mu_{aR}=\tilde\mu_4+\tilde\mu_5+\tilde\mu_6$ and
$\tilde\mu_a=\tilde\mu_0+\cdots+\tilde\mu_6$. Adding the four equations
of \eqref{sc:eq:four-coefficients} gives
\begin{equation}\label{sc:eq:four-start}
 x_k=(\tilde\mu_0+\tilde\mu_1)\bar g_0+(\tilde\mu_2+\tilde\mu_3)\bar g_1
 +(\tilde\mu_4+\tilde\mu_5)\bar g_2+\tilde\mu_6\bar g_3 .
\end{equation}
The four coefficients in \eqref{sc:eq:four-start} are positive and sum
to $\tilde\mu_a=(\rho^2+O(y^{-1/4}))h_a$.

For numbers $\eta_1,\eta_2,\eta_3$ define the bilinear form
\[
 \ip zg=\sum_{r=1}^3\eta_r(z)_r(g)_r+(z)_4(g)_4\qquad(z,g\in\R^4).
\]
Take the rows of
\begin{equation}\label{eq:four-reference-directions}
 \begin{pmatrix}
 \tilde\mu_2/\tilde\mu_1&-1&1&1\\
 -1&-1&1&1\\
 \tilde\mu_3/(\tilde\mu_4+\tilde\mu_5)&\tilde\mu_3/(\tilde\mu_4+\tilde\mu_5)&1&1\\
 0&0&0&-(\tilde\mu_a-\tilde\mu_6)/\tilde\mu_6
 \end{pmatrix}
\end{equation}
as the coordinates of $\bar g_0,\bar g_1,\bar g_2,\bar g_3$. Since
$x_j=-(\mu_7+O(P))\bar g_3$, we have
$(x_j)_1=(x_j)_2=(x_j)_3=0<(x_j)_4$. By
\eqref{sc:eq:four-coefficients},
\begin{align*}
 \bar z_2&=\bigl(\tilde\mu_3\tilde\mu_5/(\tilde\mu_4+\tilde\mu_5),\;
                 \tilde\mu_3\tilde\mu_5/(\tilde\mu_4+\tilde\mu_5),\;
                 \tilde\mu_5,\;
                 -(\tilde\mu_{aL}+\tilde\mu_3+\tilde\mu_4)\bigr),\\
 \bar z_1&=\bigl(0,\;0,\;\tilde\mu_3+\tilde\mu_4+\tilde\mu_5,\;-\tilde\mu_{aL}\bigr),\\
 \bar z_0&=\bigl(0,\;-(\tilde\mu_1+\tilde\mu_2),\;\tilde\mu_a-\tilde\mu_0-\tilde\mu_6,\;-\tilde\mu_0\bigr),\\
 x_k&=\bigl(\tilde\mu_0\tilde\mu_2/\tilde\mu_1,\;-\tilde\mu_{aL},\;\tilde\mu_a-\tilde\mu_6,\;0\bigr).
\end{align*}
Since $(\bar z_0)_1=(\bar z_1)_1=(\bar z_1)_2=0$,
\[
 \ip{\bar g_0}{\bar z_0}=\ip{\bar g_1}{\bar z_0},\qquad
 \ip{\bar g_0}{\bar z_1}=\ip{\bar g_1}{\bar z_1}=\ip{\bar g_2}{\bar z_1}.
\]
Choose $\eta_1,\eta_2,\eta_3$ such that
$\ip{\bar z_1}{\bar g_2-\bar g_3}=\ip{\bar z_2}{\bar g_2-\bar g_3}=\ip{x_k}{\bar g_2-\bar g_3}=0$.
The equation $\ip{\bar z_1}{\bar g_2-\bar g_3}=0$ determines $\eta_3$.
Since $(\bar z_2)_1=(\bar z_2)_2>0$, the equation
$\ip{\bar z_2}{\bar g_2-\bar g_3}=0$ then determines $\eta_1+\eta_2$.
Since $(x_k)_1>0>(x_k)_2$, the equation
$\ip{x_k}{\bar g_2-\bar g_3}=0$ then determines $\eta_1$. Since
$(\bar z_1)_3>0>(\bar z_1)_4$, we have $\eta_3>0$.

Since $(\bar g_3)_1=(\bar g_3)_2=(\bar g_3)_3=0$,
$(\bar g_3)_4=-(\tilde\mu_a-\tilde\mu_6)/\tilde\mu_6<0$ and
$(\bar g_2)_4=1$, we have
$\ip{\bar g_3}{\bar g_2-\bar g_3}=-(\tilde\mu_a/\tilde\mu_6)\lvert(\bar g_3)_4\rvert$.
The inner products of $\bar z_2$, $\bar z_1-\bar z_2$ and $x_k-\bar z_1$
with $\bar g_2-\bar g_3$ are zero. By \eqref{sc:eq:four-coefficients},
\[
 (\ip{\bar g_i}{\bar g_2-\bar g_3})_{i=0}^2
 =\frac{\tilde\mu_a\lvert(\bar g_3)_4\rvert}{\tilde\mu_5}
 \left(\frac{\tilde\mu_2\tilde\mu_4}{(\tilde\mu_0+\tilde\mu_1)\tilde\mu_3},
 -\frac{\tilde\mu_4}{\tilde\mu_3},1\right).
\]
Thus $\ip{\bar g_0-\bar g_1}{\bar g_2-\bar g_3}>0$. Since only the
first coordinate of $\bar g_0-\bar g_1$ is nonzero, $\eta_1>0$. Since
$\ip{x_k}{\bar g_2-\bar g_3}=0$ and $(x_k)_4=0$, we have
$\ip{x_k}{\bar g_2}=\ip{x_k}{\bar g_3}=0$ and hence
\[
 0=\ip{x_k}{\bar g_2}
  =(\bar g_2)_1\bigl[(x_k)_1\eta_1+(x_k)_2\eta_2\bigr]
  +(x_k)_3\eta_3 .
\]
Because $(\bar g_2)_1>0$, $(x_k)_1>0$, $(x_k)_2<0$, $(x_k)_3>0$ and
$\eta_1,\eta_3>0$, we obtain $\eta_2>0$. Thus $\ip\cdot\cdot$ is an
inner product on $\R^4$. Let $\norm\cdot$ be the norm of
$\ip\cdot\cdot$. We have
$\norm{\bar g_3}=\lvert(\bar g_3)_4\rvert=(\tilde\mu_a-\tilde\mu_6)/\tilde\mu_6$,
$x_k\perp\bar g_2$ and $x_k\perp\bar g_3$. Since $x_j$ is a multiple
of $\bar g_3$, also $x_k\perp x_j$. We use the values
\begin{equation}\label{sc:eq:four-reference}
 (\tilde\mu_0,\ldots,\tilde\mu_6)=h_a(1,1/\rho,1,1,1,1/\rho,1) .
\end{equation}
At \eqref{sc:eq:four-reference} we have $(x_k)_1=\rho h_a=-(x_k)_2$.
Hence $\ip{x_k}{\bar g_2}=0$ gives at \eqref{sc:eq:four-reference}
\[
 (\bar g_2)_1(x_k)_1(\eta_2-\eta_1)=(x_k)_3\eta_3>0,\qquad
 \eta_2>\eta_1 .
\]

By $\bar z_1\perp\bar g_2-\bar g_3$ and
$\ip{\bar g_0}{\bar z_1}=\ip{\bar g_1}{\bar z_1}=\ip{\bar g_2}{\bar z_1}$,
the inner product $\ip{\bar g_i}{\bar z_1}$ is independent of
$i\in\{0,1,2,3\}$. Since $(\bar z_1)_4=-\tilde\mu_{aL}$, the common
value of $\ip{\bar g_i}{\bar z_1}$ is $\tilde\mu_{aL}\norm{\bar g_3}$.
The coefficients in \eqref{sc:eq:four-start} sum to $\tilde\mu_a$. We
have $x_k\perp\bar g_2$, $x_k\perp\bar g_3$ and
$\bar z_1=\tilde\mu_3\bar g_1+(\tilde\mu_4+\tilde\mu_5)\bar g_2+\tilde\mu_6\bar g_3$
by \eqref{sc:eq:four-coefficients}. Only the first coordinate of
$\bar g_0-\bar g_1$ is nonzero and $(\bar g_2-\bar g_3)_1=(\bar g_2)_1$.
Hence
\begin{equation}\label{eq:four-start-pairings}
 \tilde\mu_3\ip{x_k}{\bar g_1}=\ip{x_k}{\bar z_1}
 =\tilde\mu_a\tilde\mu_{aL}\norm{\bar g_3},\qquad
 \ip{x_k}{\bar g_0-\bar g_1}
 =\frac{(x_k)_1}{(\bar g_2)_1}\ip{\bar g_0-\bar g_1}{\bar g_2-\bar g_3}.
\end{equation}
We have $\tilde\mu_{aL}+\tilde\mu_3=\tilde\mu_a-\tilde\mu_{aR}$,
$\tilde\mu_4+\tilde\mu_5=\tilde\mu_{aR}-\tilde\mu_6$ and
\[
 \frac{\tilde\mu_2\tilde\mu_4}{(\tilde\mu_0+\tilde\mu_1)\tilde\mu_3}
 +\frac{\tilde\mu_4}{\tilde\mu_3}
 =\frac{\tilde\mu_4\tilde\mu_{aL}}{(\tilde\mu_0+\tilde\mu_1)\tilde\mu_3}.
\]
By \eqref{sc:eq:four-start}, $x_k\perp\bar g_2$ and
$x_k\perp\bar g_3$,
\[
 \norm{x_k}^2=(\tilde\mu_{aL}+\tilde\mu_3)\ip{x_k}{\bar g_1}
 +(\tilde\mu_0+\tilde\mu_1)\ip{x_k}{\bar g_0-\bar g_1} .
\]
Inserting \eqref{eq:four-start-pairings} gives
\begin{equation}\label{eq:high-norm-factorization}
 \frac{\tilde\mu_3\norm{x_k}^2}{\norm{\bar g_3}\tilde\mu_{aL}\tilde\mu_a}
 =\tilde\mu_a-\tilde\mu_{aR}
 +(\tilde\mu_{aR}-\tilde\mu_6)\frac{\tilde\mu_0\tilde\mu_2\tilde\mu_4}{\tilde\mu_1\tilde\mu_3\tilde\mu_5} .
\end{equation}
Define $\tilde Q_0=\mu_7/\tilde\mu_a$,
$\tilde Q_1=\mu_7\tilde\mu_3/(\tilde\mu_{aL}\tilde\mu_{aR})$ and
$\tilde Q_2=\tilde\mu_1\tilde\mu_3\tilde\mu_5\mu_7/(\tilde\mu_0\tilde\mu_2\tilde\mu_4\tilde\mu_6)$.
Multiplying \eqref{eq:high-norm-factorization} by
$\tilde\mu_{aL}\tilde\mu_a/(\tilde\mu_3\mu_7^2\norm{\bar g_3})$ and
substituting $\norm{\bar g_3}=(\tilde\mu_a-\tilde\mu_6)/\tilde\mu_6$
gives
\[
 \frac{\norm{x_k}^2}{\mu_7^2\norm{\bar g_3}^2}
 =\tilde Q_1^{-1}\left(
 \frac{\tilde\mu_6(\tilde\mu_a-\tilde\mu_{aR})}{(\tilde\mu_a-\tilde\mu_6)\tilde\mu_{aR}}\tilde Q_0^{-1}
 +\frac{\tilde\mu_a(\tilde\mu_{aR}-\tilde\mu_6)}{(\tilde\mu_a-\tilde\mu_6)\tilde\mu_{aR}}\tilde Q_2^{-1}
 \right).
\]
Let $\beta$ be the coefficient of $\tilde Q_2^{-1}$. The coefficient
of $\tilde Q_0^{-1}$ is $1-\beta$. By \eqref{sc:eq:tilde-balance},
$\beta=1/\sqrt2+O(y^{-1/4})$. Put
$\delta=\log(\tilde Q_0/\tilde Q_2)$. By \eqref{sc:eq:tilde-balance},
$\delta=O(y^{-1/4})$. By the mean value theorem for
$\tau\mapsto\log(1-\beta+\beta e^\tau)$ there is $\tau_0$ between $0$
and $\delta$ such that
\[
 \log\bigl((1-\beta)\tilde Q_0^{-1}+\beta\tilde Q_2^{-1}\bigr)
 =(1-\tilde\beta)\log(\tilde Q_0^{-1})+\tilde\beta\log(\tilde Q_2^{-1}),
 \qquad
 \tilde\beta=\frac{\beta e^{\tau_0}}{1-\beta+\beta e^{\tau_0}} .
\]
We have
$\tilde\beta=\beta+O(\lvert\delta\rvert)=1/\sqrt2+O(y^{-1/4})$ and
$\tilde\beta\in[0,1]$. Each of
$\mu_0,\ldots,\mu_7,m_{aL},m_{aR},m_a$ is at least $ch_a$ and
$\tilde\mu_r=\mu_r+O(P)$ for $r=0,\ldots,6$. Therefore
$\log\tilde Q_r=\log Q_r(\mu_7)+O(P/h_a)$ for $r=0,1,2$. Since
$\norm{x_j}=(\mu_7+O(P))\norm{\bar g_3}$, we have
$\log(\norm{x_k}^2/\norm{x_j}^2)=\log[\norm{x_k}^2/(\mu_7^2\norm{\bar g_3}^2)]+O(P/h_a)$.
Hence
$-\tfrac12\log(\norm{x_k}^2/\norm{x_j}^2)$ equals the right side of
\eqref{sc:eq:four-gain} with $\alpha_0=(1-\tilde\beta)/2$,
$\alpha_1=1/2$ and $\alpha_2=\tilde\beta/2$. Thus \ref{it:four-gain}
holds. The equations
$\ip{\bar z_1/h_a}{\bar g_2-\bar g_3}=\ip{\bar z_2/h_a}{\bar g_2-\bar g_3}=\ip{x_k/h_a}{\bar g_2-\bar g_3}=0$
form a linear system for $\eta_1,\eta_2,\eta_3$. By
\eqref{sc:eq:tilde-balance} the coefficients of the linear system
differ by $O(y^{-1/4})$ from the coefficients at
\eqref{sc:eq:four-reference}. At \eqref{sc:eq:four-reference} the
linear system is nonsingular. Hence $\eta_1,\eta_2,\eta_3$ and
$\norm{\bar g_0},\ldots,\norm{\bar g_3}$ are at most an absolute
constant. Since
$\norm{\bar g_3}=2\rho+o(1)$ and $\mu_7\ge(\rho-o(1))h_a$,
\eqref{eq:high-norm-factorization} gives
\begin{align}
 ch_a&\le\norm{x_k}\le Ch_a,\label{sc:eq:four-units}\\
 (2\rho^2-o(1))h_a&\le\norm{x_j}\le Ch_j,\label{sc:eq:four-end}\\
 \max_i\norm{\bar g_i}&\le C .\label{sc:eq:four-slopes}
\end{align}

At \eqref{sc:eq:four-reference} the rows of
\eqref{eq:four-reference-directions} are
\[
 \bar g_0=(\rho,-1,1,1),\quad
 \bar g_1=(-1,-1,1,1),\quad
 \bar g_2=(1/\sqrt2,1/\sqrt2,1,1),\quad
 \bar g_3=(0,0,0,-2\rho).
\]
Since $\rho-1/\sqrt2=1+1/\sqrt2$, the definition of $\ip\cdot\cdot$
gives at \eqref{sc:eq:four-reference}
\begin{equation}\label{sc:eq:four-products}
\begin{aligned}
 \ip{\bar g_0}{\bar g_0-\bar g_1}&=\rho(\rho+1)\eta_1, &
 \ip{\bar g_1}{\bar g_1-\bar g_0}&=(\rho+1)\eta_1,\\
 \ip{\bar g_0}{\bar g_0-\bar g_2}&=(1+1/\sqrt2)(\rho\eta_1+\eta_2), &
 \ip{\bar g_2}{\bar g_2-\bar g_0}&=(1+1/\sqrt2)(\eta_2-\eta_1)/\sqrt2,\\
 \ip{\bar g_1}{\bar g_1-\bar g_2}&=(1+1/\sqrt2)(\eta_1+\eta_2), &
 \ip{\bar g_2}{\bar g_2-\bar g_1}&=(1+1/\sqrt2)(\eta_1+\eta_2)/\sqrt2,\\
  \ip{\bar g_i}{\bar g_i-\bar g_3}
 &=(\bar g_i)_1^2\eta_1+(\bar g_i)_2^2\eta_2+\eta_3+1+2\rho, &
 \ip{\bar g_3}{\bar g_3-\bar g_i}&=2\rho(2\rho+1)\quad(i=0,1,2).
\end{aligned}
\end{equation}
Since $\eta_1,\eta_2,\eta_3>0$ and $\eta_2>\eta_1$ at
\eqref{sc:eq:four-reference}, every inner product in
\eqref{sc:eq:four-products} is positive. Thus at
\eqref{sc:eq:four-reference}
\begin{equation}\label{sc:eq:four-angles}
 \ip{\bar g_i}{\bar g_i-\bar g_r}>0\qquad(i,r\in\{0,1,2,3\},\ i\ne r).
\end{equation}
Let $r\in\{0,1\}$. We have $(\bar g_2-\bar g_r)_4=0$ and
$(\bar g_3)_1=(\bar g_3)_2=(\bar g_3)_3=0$. Hence
$\ip{\bar g_2-\bar g_r}{\bar g_3}=0$. Since
$\bar z_2=\tilde\mu_5\bar g_2+\tilde\mu_6\bar g_3$ by
\eqref{sc:eq:four-coefficients}, we have
$\ip{\bar g_2-\bar g_r}{\bar z_2}=\tilde\mu_5\ip{\bar g_2}{\bar g_2-\bar g_r}$.
Thus at \eqref{sc:eq:four-reference}
\begin{equation}\label{sc:eq:four-z2}
 \ip{\bar g_2}{\bar z_2}>\ip{\bar g_0}{\bar z_2},\qquad
 \ip{\bar g_2}{\bar z_2}>\ip{\bar g_1}{\bar z_2}.
\end{equation}
The inner products in \eqref{sc:eq:four-angles} and the inner products
in \eqref{sc:eq:four-z2} divided by $h_a$ are continuous functions of
$(\tilde\mu_0,\ldots,\tilde\mu_6)/h_a$ near
$(1,1/\rho,1,1,1,1/\rho,1)$. Hence by \eqref{sc:eq:tilde-balance},
\eqref{sc:eq:four-angles} and \eqref{sc:eq:four-z2} hold for
sufficiently large $y$.

The coordinates of $\bar z_0$, $\bar z_1$ and the rows of
\eqref{eq:four-reference-directions} give
\begin{gather*}
 \bar z_0-\bar z_1=(\tilde\mu_1+\tilde\mu_2)(0,-1,1,1),\qquad
 \bar g_1-\bar g_3=(-1,-1,1,\tilde\mu_a/\tilde\mu_6),\\
 \bar g_1-\bar g_2=-\Bigl(1+\frac{\tilde\mu_3}{\tilde\mu_4+\tilde\mu_5}\Bigr)(1,1,0,0).
\end{gather*}
Hence
\[
 \ip{\bar g_1-\bar g_2}{\bar z_0-\bar z_1}
 =(\tilde\mu_1+\tilde\mu_2)\Bigl(1+\frac{\tilde\mu_3}{\tilde\mu_4+\tilde\mu_5}\Bigr)\eta_2>0,
 \qquad
 \ip{\bar g_1-\bar g_3}{\bar z_0-\bar z_1}
 =(\tilde\mu_1+\tilde\mu_2)\Bigl(\eta_2+\eta_3+\frac{\tilde\mu_a}{\tilde\mu_6}\Bigr)>0 .
\]
Since $\ip{\bar g_r}{\bar z_1}$ is independent of $r\in\{0,1,2,3\}$, we
have
$\ip{\bar g_1-\bar g_r}{\bar z_0}=\ip{\bar g_1-\bar g_r}{\bar z_0-\bar z_1}>0$
for $r=2,3$. With $\ip{\bar g_0}{\bar z_0}=\ip{\bar g_1}{\bar z_0}$ we
obtain
\begin{equation}\label{sc:eq:four-z0}
 \ip{\bar g_0}{\bar z_0}=\ip{\bar g_1}{\bar z_0}
 >\max\{\ip{\bar g_2}{\bar z_0},\ip{\bar g_3}{\bar z_0}\}.
\end{equation}
We show
\begin{equation}\label{sc:eq:four-max}
 \ip{\bar g_i}{z_{u_i}}\ge\ip{\bar g_r}{z_{u_i}}\qquad(i,r\in\{0,1,2,3\}).
\end{equation}
We have $z_{u_i}=\bar z_i/E_{2i}$ for $i=0,1,2$. For $i=0$, \eqref{sc:eq:four-max} follows from
\eqref{sc:eq:four-z0}. For $i=1$, $\ip{\bar g_r}{\bar z_1}$ is
independent of $r$. For $i=2$, \eqref{sc:eq:four-max} follows from
\eqref{sc:eq:four-z2} and $\ip{\bar z_2}{\bar g_2-\bar g_3}=0$. For
$i=3$, \eqref{sc:eq:four-max} follows from $z_{u_3}=P\bar g_3$ and
\eqref{sc:eq:four-angles}.

For $i=0,1,2,3$ and $k_i\le u<j_i$ define
\[
 q_u=\ip{\bar g_i}{z_u}-P .
\]
Let $i\in\{0,1,2,3\}$. By \eqref{sc:eq:four-relations} and
\eqref{sc:eq:four-coefficients}, $x_{k_i}$ is a combination of
$\bar g_0,\ldots,\bar g_3$ with coefficients $O(h_a)$. The choice of $u_i$ gives
$\prod_{w=k_i+1}^u(1-h_w/\kappa)\ge\prod_{w=k_i+1}^{u_i}(1-h_w/\kappa)$
for $k_i\le u<j_i$. By \eqref{sc:eq:four-constant}, $z_u-z_{u_i}$ is a
nonnegative multiple of
\[
 (1-1/\kappa)\bar g_i+x_{k_i}/\kappa=\bar g_i+O(h_a/\kappa)
\]
for $k_i\le u<j_i$.
By \eqref{sc:eq:four-angles} and $h_a/\kappa=o(1)$, we have
$\ip{\bar g_i-\bar g_r}{z_u-z_{u_i}}\ge0$ for every
$r\in\{0,1,2,3\}$. Together with \eqref{sc:eq:four-max} we obtain
\begin{equation}\label{sc:eq:four-order}
 \ip{\bar g_i}{z_u}\ge\ip{\bar g_r}{z_u}
 \qquad(k_i\le u<j_i,\ r\in\{0,1,2,3\}).
\end{equation}
Let $i\in\{0,1,2\}$ and $k_i\le u<j_i$. The fourth coordinates of
$\bar z_0$, $\bar z_1$, $\bar z_2$ and $\bar g_3$ give
$\ip{\bar g_3}{\bar z_i}\ge ch_a$. Since $z_u-z_{u_i}$ is a nonnegative
multiple of $\bar g_i+O(h_a/\kappa)$, we have
$\ip{\bar g_i}{z_u}\ge\ip{\bar g_i}{z_{u_i}}$. By
\eqref{sc:eq:four-max},
$\ip{\bar g_i}{z_{u_i}}\ge\ip{\bar g_3}{z_{u_i}}\ge ch_a$. Hence
$q_u\ge ch_a$ for $k\le u<aR$. Since $z_{u_3}=P\bar g_3$, \eqref{sc:eq:four-relations} for $i=3$
shows that $x_{aR}$ is a positive multiple of $\bar g_3$. By \eqref{sc:eq:four-constant} every $z_u$ with
$aR\le u<j$ is $P\bar g_3$ plus a nonnegative multiple of
$\bar g_3$. Since $\norm{\bar g_3}>1$, we have $q_u\ge0$ for
$aR\le u<j$. Thus $0\le q_u\le\ip{g_u}{z_u}$ for every $k\le u<j$. Let
$i,r\in\{0,1,2,3\}$, $k_i\le u<j_i$ and $k_r\le u'<j_r$. By
\eqref{sc:eq:four-order}, $\ip{\bar g_r}{z_{u'}}\ge\ip{\bar g_i}{z_{u'}}$.
Hence
\[
 q_u+\ip{\bar g_i}{z_{u'}-z_u}=\ip{\bar g_i}{z_{u'}}-P
 \le\ip{\bar g_r}{z_{u'}}-P=q_{u'} .
\]
Hence \eqref{eq:tail-planes} holds for every pair of triples among
$(0,0,0)$ and the triples at $k\le u<j$. The points $x_k,\ldots,x_j$
and the triples at $k\le u<j$ form a four-dimensional solution
$\mathcal S$ on $(k,j]$.

Since $j\in B$ and $I_a\subset I_{jL}$, \eqref{sc:eq:balanced-debt} and
$e^y/h_a\le y^{1/12}$ give
\[
 \sum_{w\in I_a}h_w^-\le Cy^{-1/4}e^y\le Cy^{-1/6}h_a=o(h_a).
\]
Since $a\in B$, $\sum_{w\in I_a}h_w^+=O(h_a)$. Hence
\eqref{sc:eq:four-constant} gives $\norm{z_u}=O(h_a)$ for $k\le u<j$.
By \eqref{eq:four-start-pairings},
$\ip{x_k}{\bar g_0}>\ip{x_k}{\bar g_1}>0$. Also $x_k\perp\bar g_2$ and
$x_k\perp\bar g_3$. Thus $\ip{x_k}{g_u}\ge0$ for every $k\le u<j$. Let
$i\in\{0,1,2,3\}$ and $k_i\le u<j_i$. Then $z_u-z_{u_i}$ is a
nonnegative multiple of $(1-1/\kappa)\bar g_i+x_{k_i}/\kappa$. By
\eqref{sc:eq:four-relations} and \eqref{sc:eq:four-coefficients},
$x_{k_i}$ and $z_{u_i}$ are combinations of
$\bar g_0,\ldots,\bar g_3$ with nonnegative coefficients. Hence $z_u$
is a combination of $\bar g_0,\ldots,\bar g_3$ with nonnegative
coefficients and therefore $\ip{x_k}{z_u}\ge0$. Let $k\le u<aR$. By
\eqref{sc:eq:four-order} and $q_u\ge ch_a$, we have
$\ip{\bar g_r}{z_u}\le q_u+P=(1+o(1))q_u$ for every
$r\in\{0,1,2,3\}$. By \eqref{sc:eq:four-start}, $x_k$ is a combination
of $\bar g_0,\ldots,\bar g_3$ with positive coefficients that sum to
$(\rho^2+o(1))h_a$. Hence
$\ip{x_k}{z_u}\le(\rho^2+o(1))h_a q_u\le(\rho^2+1)h_a q_u$. Let
$aR\le u<j$. Then $z_u$ is a multiple of $\bar g_3$ and hence
$\ip{x_k}{z_u}=0$. Together with \eqref{sc:eq:four-units} and
\eqref{sc:eq:four-slopes}, \eqref{sc:eq:four-boundary} holds with
absolute constants $c_{\mathrm{four}}$ and $C_{\mathrm{four}}$. By
\eqref{sc:eq:band-B},
$A=e^yy^3$ and $P=e^y/y^3\le h_a$. By \eqref{sc:eq:four-units} and
\eqref{sc:eq:four-slopes},
\[
 \norm{z_u}+A\norm{g_u}\le Ch_a+CA\le C(A/P)^2\norm{x_k}.
\]
Since $\ip{\widehat{x_k}}{z_u}\ge0$, $\mathcal S$ satisfies
\eqref{sc:eq:start-projections} with constant $C_{\mathrm{four}}$
after an increase of $C_{\mathrm{four}}$. Thus \ref{it:four-boundary}
holds.

One rotation and \eqref{eq:scaling} with one $\sigma$ applied to
$\tilde{\mathcal S}$ and to $\tilde{\mathcal S}'$ preserve the
hypotheses and the conclusion of \ref{it:four-join}. Hence it suffices
to prove \ref{it:four-join} for $\tilde{\mathcal S}=\mathcal S$. Let
$c_0$, $C_0$, $y$, $\mathcal S'$ and $\tilde{\mathcal S}'$ satisfy the
hypotheses of \ref{it:four-join} with $\tilde{\mathcal S}=\mathcal S$.
Let $x'$ be the first point of $\mathcal S'$. We have $(g_u)_4=1$ for $k\le u<aR$,
$\widehat{x_j}=(0,0,0,1)$ and $(x_k)_4=0$. We have
$1-h_w/\kappa\le\exp(h_w^-/\kappa)$ for every $w\in I_a$ and
$\sum_{w\in I_a}h_w^-\le Cy^{-1/4}e^y$ with
$y^{-1/4}e^y/\kappa=o(1)$. Hence the fourth coordinate of
\eqref{sc:eq:four-constant} for $i=0,1,2$ gives
\begin{equation}\label{sc:eq:four-last-axis}
 \ip{\widehat{x_j}}{z_u}+1
 =(\kappa-1)\left(\prod_{w=k+1}^{u}(1-h_w/\kappa)-1\right)
 \le(\kappa-1)\left[\exp\Bigl(\frac1\kappa\sum_{w\in I_a}h_w^-\Bigr)-1\right]
 \le Cy^{-1/4}e^y
\end{equation}
for $k\le u<aR$.
For $aR\le u<j$ we have $\ip{\widehat{x_j}}{z_u}<0$. By
\eqref{sc:eq:four-end} and $h_j\le e^y$,
$(2\rho^2-o(1))h_a\le\norm{x_j}\le Ce^y$. Let $\sigma$ be the factor of
\eqref{eq:scaling} in the map from $\mathcal S'$ to
$\tilde{\mathcal S}'$. Then $\sigma=\norm{x_j}/\norm{x'}$. Since
$c_0h_a\le\norm{x'}\le C_0h_a$ by \eqref{sc:eq:right-piece}, we have
$\sigma\le Ce^y/h_a\le Cy^{1/12}$. Let $(z,g,q)$ be a triple of
$\tilde{\mathcal S}'$ and let $(z_u,g_u,q_u)$ be a triple of
$\mathcal S$. By \eqref{sc:eq:right-piece} and \eqref{eq:scaling},
$\norm g\le C_0\sigma\le Cy^{1/12}$ and $\ip{x_j}g\ge0$. By
\eqref{eq:adjacent-support} for $\mathcal S$ and
$\tilde{\mathcal S}'$ with shared point $x_j$,
\[
 \ip g{z_u}=\ip g{\widehat{x_j}}\ip{\widehat{x_j}}{z_u},\qquad
 \ip{g_u}z=\ip{g_u}{\widehat{x_j}}\ip{\widehat{x_j}}z .
\]
We have $q\le\ip gz$. For $k\le u<aR$ we have $q_u\ge ch_a$. Hence
\eqref{sc:eq:four-last-axis} gives
\[
 q+\ip g{z_u-z}\le\ip g{\widehat{x_j}}\ip{\widehat{x_j}}{z_u}
 \le Cy^{1/12}y^{-1/4}e^y\le Cy^{-1/12}h_a\le q_u .
\]
For $aR\le u<j$ we have $\ip{\widehat{x_j}}{z_u}<0$ and therefore
\[
 q+\ip g{z_u-z}\le0\le q_u .
\]
By \eqref{sc:eq:right-piece} and \eqref{eq:scaling},
$0\le\ip{x_j}z\le(\rho^2+2)h_aq+\sigma^2P^2$. Dividing by
$\norm{x_j}\ge(2\rho^2-o(1))h_a$ and using
$\sigma^2P/h_a\le Cy^{-11/4}=o(1)$ and $(\rho^2+2)/(2\rho^2)<1$ gives
\[
 0\le\ip{\widehat{x_j}}z\le q+P .
\]
Every triple of $\mathcal S$ satisfies $\ip{g_u}{\widehat{x_j}}\le1$
and $q_u-\ip{g_u}{z_u}=-P$. Hence
\[
 q_u+\ip{g_u}{z-z_u}=\ip{g_u}{\widehat{x_j}}\ip{\widehat{x_j}}z-P\le q .
\]
Thus \eqref{eq:tail-planes} holds for every triple of $\mathcal S$
together with every triple of $\tilde{\mathcal S}'$ and
\ref{it:four-join} holds.
\end{proof}

\subsection{The values of \texorpdfstring{$\mathcal P$}{P} on two ranges of thresholds}\label{sc:sec:conclusion}\label{sc:sec:pairs}

Let $t$ be sufficiently large with $e^t t^3\le h_n$ and let
$t/2\le s\le y\le t$. Let $j$ be a cut at threshold $e^y$ such that $I$
contains a positive step. The positive steps of $I$ form a binary tree
in which the children of $v$ are $vL$ and $vR$. The first step removed
from $I$ is the \emph{root} of $I$ and has \emph{depth} zero. The
children of a step of depth $\nu$ have depth $\nu+1$. By
\eqref{eq:tree-split} and the order of removal, the parent of $I$ is a
cut at threshold $e^y$.

Let $N$ be the smallest integer with $2^N\ge(\log t)^2$.
\Cref{sc:lem:two-ranges} applies \cref{sc:prop:working-realization} for
$\ell=N,N+1,N+2$. The set $S$ consists of the cuts $j$ with
\eqref{eq:three-cuts} for which the parent of $I$ and all positive
steps of $I$ of depth at most $N+2$ lie in $B$. Let $j\in S$. For a
step $a$ of depth $N$ in $I$ let $k_a$ be the index before $I_a$ and
let $j_a$ be the index after $I_a$. The sets $(k_a,j_a]$ form a
partition of $(k,j]$. Let $\omega$ be the largest step of depth $N$ in
$I$. Then $j_{\omega}=j$. For $a\ne \omega$ we take a solution given by
\cref{sc:lem:four-phase} for $a$ on $(k_a,j_a]$. Since $h_j>e^y$,
\cref{sc:lem:four-phase} does not apply to $\omega$. For $\ell=N$ we take a solution given by
\cref{sc:lem:ordinary} on $(k_{\omega},j]$. For $\ell=N+1$ we take
solutions given by \cref{sc:lem:ordinary} on two intervals with union
$(k_{\omega},j]$. For $\ell=N+2$ we take solutions given by
\cref{sc:lem:ordinary} on four intervals with union $(k_{\omega},j]$.
In \cref{sc:lem:two-ranges}, $\psi_N,\psi_{N+1},\psi_{N+2}$ equal
$\alpha_0,\alpha_1,\alpha_2$ of \eqref{sc:eq:four-gain} up to
$O(t^{-1/4})$. The sum of \eqref{sc:eq:working-quantitative} for
$\ell=N,N+1,N+2$ with weights $\psi_\ell$ gives
\eqref{eq:weighted-score}. The weights of
$\mathcal P(t-\theta)-\mathcal P(s-\theta)$ in \eqref{eq:weighted-score}
have integral $p\sum_{\ell=N}^{N+2}\psi_\ell(2^\ell-1)\ge c(\log t)^2$.

\begin{lemma}\label{sc:lem:two-ranges}
Let $t$ be sufficiently large with $e^t t^3\le h_n$ and let
$t/2\le s\le 2t/3$. Let $N$ be the smallest integer with
$2^N\ge(\log t)^2$ and define
\[
 \psi_N=\frac12-\frac1{2\sqrt2},\qquad
 \psi_{N+1}=\frac12,\qquad
 \psi_{N+2}=\frac1{2\sqrt2}.
\]
Then
\begin{equation}\label{eq:weighted-score}
 \sum_{\ell=N}^{N+2}\psi_\ell\int_0^{\ell\log\rho}e^{\theta/p}
 [\mathcal P(t-\theta)-\mathcal P(s-\theta)]\,d\theta\le Cn .
\end{equation}
\end{lemma}

\begin{proof}
For each $\ell\in\{N,N+1,N+2\}$, the substitution $\theta=y-\log h_b$
for each $b$ together with \eqref{sc:eq:band-sums} gives
\begin{equation}\label{sc:eq:test-score}
 \int_s^t e^{y/p}\sum_{e^y/\rho^\ell<h_b\le e^y}
                 \frac{\Delta_b}{h_b^{1/p}}\,dy
 =\int_0^{\ell\log\rho}e^{\theta/p}
       [\mathcal P(t-\theta)-\mathcal P(s-\theta)]\,d\theta .
\end{equation}
By \eqref{sc:eq:test-score}, \eqref{eq:weighted-score} is equivalent to
\begin{equation}\label{eq:weighted-delta-score}
 \sum_{\ell=N}^{N+2}\psi_\ell\int_s^t e^{y/p}
 \sum_{e^y/\rho^\ell<h_b\le e^y}\frac{\Delta_b}{h_b^{1/p}}\,dy\le Cn .
\end{equation}
In this proof $o(\cdot)$ refers to $t\to\infty$ and $\sum_\ell$ denotes
the sum over $\ell=N,N+1,N+2$. For $s\le y\le t$,
\cref{sc:lem:four-phase}\ref{it:four-gain} gives
\eqref{sc:eq:four-gain} with $\alpha_0=\psi_N+O(t^{-1/4})$, $\alpha_1=\psi_{N+1}$ and
$\alpha_2=\psi_{N+2}+O(t^{-1/4})$. Every fixed power of $2^N$ or of
$\rho^N$ is $t^{o(1)}$. Let $s\le y\le t$. For each $\ell\in\{N,N+1,N+2\}$ we
apply \cref{sc:prop:working-realization} at threshold $e^y$ with sets
$S$ and $V_\ell(y)$.

Let $S$ be the set of the cuts $j$ at threshold $e^y$ with
\eqref{eq:three-cuts} such that the parent of $I$ and all positive
steps of $I$ of depth at most $N+2$ lie in $B$. If $I$
contains no positive step, the parent of $I$ is not in $B$ by
\cref{sc:lem:balance}\ref{it:bal-prop}. Hence $I$ contains a positive
step for every $j\in S$. Let $V_\ell(y)$ be the set of the positive
steps $b$ with $h_b>e^y/\rho^\ell$ in the intervals $I$ of the cuts
$j\notin S$. For fixed $\ell$ and $y$ write $V=V_\ell(y)$. For $j\notin S$ the set
$V\cap I$ is an initial segment of the positive steps of $I$ in the
removal order. Every $b\in V$ satisfies
\[
 h_b>e^y/\rho^{N+2}>e^y/y>P(y) .
\]
By \eqref{sc:eq:bulk-split} with $P=P(y)$ and $y^{-2/p}\le1/y$,
\begin{equation}\label{eq:V-sum}
 \sum_{b\in V}\log\frac{h_b(m_b+P(y))}{(m_{bL}+P(y))(m_{bR}+P(y))}
 \ge\sum_{b\in V}\Bigl(\frac{\Delta_b}{h_b^{1/p}}-\frac Cy\Bigr).
\end{equation}

Let $j\in S$ and let $v$ be the root of $I$. Then $I=I_v$. The parent of $I$ has stepsize above
$e^y$ and $h_v\le e^y$. Since the parent of $I$ and $v$ lie in $B$,
\cref{sc:lem:balance}\ref{it:bal-parent} gives $h_v\ge ce^y$ and hence
$\log h_v\ge t/2-C$. By \cref{sc:lem:balance}\ref{it:bal-parent} the quotient of the
stepsizes of a parent in $B$ and a child in $B$ is at most an absolute
constant. Since $N=O(\log\log t)$, every step $b$ of $I$ of depth at
most $N+2$ has $t/3<\log h_b\le t$. By
\cref{sc:lem:balance}\ref{it:bal-parent} applied $N$ times, every step
$a$ of depth $N$ in $I$ satisfies
\[
 h_a=h_v\rho^{-N}e^{O(Nt^{-1/4})},\qquad
 1\le e^y/h_a\le C\rho^N\le y^{1/12}.
\]
The stepsizes of any two steps of depth $N$ in $I$ differ by a factor
$1+o(1)$. For a step $a$ of depth $N$ in $I$ let $k_a$ be the index
before $I_a$ and let $j_a$ be the index after $I_a$. The steps of $I$
of depth at most $N-1$ together with $j$ divide $(k,j]$ into the sets
$(k_a,j_a]$ with $a$ of depth $N$ in $I$. Let $\omega$ be the largest
step of depth $N$ in $I$. Then $j_{\omega}=j$.

Let $a\ne \omega$ be a step of depth $N$ in $I$. Then $j_a$ is a step
of $I$ of depth at most $N-1$. Hence $a,aL,aR,j_a\in B$ and
$h_a\le h_{j_a}\le e^y$. Let $\mathcal S_a$ be a solution given by
\cref{sc:lem:four-phase} for $a$. The interval of $\mathcal S_a$ is
$(k_a,j_a]$.

Since $h_j>e^y$, \cref{sc:lem:four-phase} does not apply to $\omega$.
Define
\[
 \Gamma_N=\{(k_\omega,j]\},\qquad
 \Gamma_{N+1}=\{(k_\omega,\omega],(\omega,j]\},\qquad
 \Gamma_{N+2}=\{(k_\omega,\omega L],(\omega L,\omega],(\omega,\omega R],(\omega R,j]\}.
\]
Let $\ell\in\{N,N+1,N+2\}$ and $(k',b]\in\Gamma_\ell$. Then
$(k',b]\subset(k,j]$ and
$\sum_{k'<w\le b}\lvert h_w\rvert\le\sum_{k<w\le j}\lvert h_w\rvert\le A$
by \eqref{eq:three-cuts}. Here $b$ is one of
$\omega L,\omega,\omega R,j$. Each
$u\in\{\omega L,\omega,\omega R\}$ satisfies
$h_u\ge ce^y/\rho^{N+1}>P$. Also $h_j>e^y>P$. Hence the hypotheses of
\cref{sc:lem:ordinary} hold on $(k',b]$. Let $\mathcal S_{(k',b]}$ be a
solution given by \cref{sc:lem:ordinary} on $(k',b]$.

Let $\ell\in\{N,N+1,N+2\}$. Let $\Gamma$ be the family of the sets
$(k_a,j_a]$ with $a\ne\omega$ of depth $N$ in $I$ and of the sets of
$\Gamma_\ell$. Then $\Gamma$ is a partition of $(k,j]$ into intervals.
Let $r_j=\#\Gamma$. The first and last points of each
$\mathcal S_a$ with $a\ne\omega$ and of each $\mathcal S_{(k',b]}$
with $(k',b]\in\Gamma_\ell$ are nonzero and orthogonal. Hence
rotations and \eqref{eq:scaling} in a Euclidean space of sufficiently
large dimension give images
$\mathcal S_{j,1},\ldots,\mathcal S_{j,r_j}$ of the $r_j$ solutions,
numbered in increasing order of the sets of $\Gamma$, that satisfy
\cref{sc:lem:assembly}\ref{it:asm-shared}--\ref{it:asm-perp}. Here
$\mathcal S_{j,1},\ldots,\mathcal S_{j,r_j}$ depend on $\ell$. For
$k\le u\le j$ write $x_u^{(\ell)}$ for the point at index
$u$ of $\mathcal S_{j,1},\ldots,\mathcal S_{j,r_j}$.

For a step $a$ of depth $N$ in $I$ write $\mu_{a,0},\ldots,\mu_{a,7}$
and $Q_{a,0},Q_{a,1},Q_{a,2}$ for $\mu_0,\ldots,\mu_7$ of
\eqref{sc:eq:four-mu} and $Q_0,Q_1,Q_2$ of
\eqref{sc:eq:high-gain-score} defined for $a$. By
\eqref{sc:eq:four-masses} for $\omega$ and $i=0,1,2$ and by
\eqref{sc:eq:balanced-unit-upper},
$\mu_{\omega,2r+1}\le h_u+h_u^2/\kappa$ for
$(u,r)=(\omega L,0),(\omega,1),(\omega R,2)$. Since
$h_u/\kappa\le y^{-3}$, we have $h_u\ge(1-y^{-3})\mu_{\omega,2r+1}$.
Each $m$ among
$m_\omega,m_{\omega L},m_{\omega R},\mu_{\omega,0},\mu_{\omega,2},\mu_{\omega,4},\mu_{\omega,6}$
satisfies
$P/m\le Cy^{-2}$ and $m/\kappa\le Cy^{-3}$. A rotation and
\eqref{eq:scaling} do not change the quotient of the norms of two
points of one solution. Hence \eqref{sc:eq:ordinary-gain} for
$\mathcal S_{(k',b]}$ gives
\[
 \frac{\norm{x_b^{(\ell)}}}{\norm{x_{k'}^{(\ell)}}}
 \ge\frac{h_b}{m_J+P}\exp\bigl(-C/\log A-m_J/\kappa\bigr),
 \qquad J=\{k'+1,\ldots,b-1\},
\]
for every $(k',b]\in\Gamma_\ell$. Here $J=I_\omega$ for $\ell=N$,
$J\in\{I_{\omega L},I_{\omega R}\}$ for $\ell=N+1$ and
$J\in\{I_{\omega LL},I_{\omega LR},I_{\omega RL},I_{\omega RR}\}$ for
$\ell=N+2$. Multiplying over $(k',b]\in\Gamma_\ell$ gives
\begin{equation}\label{sc:eq:last-gain}
  \frac{\norm{x_j^{(\ell)}}}{\norm{x_{k_\omega}^{(\ell)}}}
 \ge e^{-C/t}Q_{\omega,\ell-N}(h_j)\qquad(\ell=N,N+1,N+2).
\end{equation}

We verify \cref{sc:prop:working-realization}\ref{it:work-S} with
$C_0=C_{\mathrm{four}}$. Let $1\le i<i'\le r_j$ with $i'\ge i+2$. By
\cref{sc:lem:assembly}\ref{it:asm-shared}--\ref{it:asm-perp}, every
triple $(z,g,q)$ of $\mathcal S_{j,i}$ and every triple $(z',g',q')$ of
$\mathcal S_{j,i'}$ satisfy $\ip g{z'}=\ip{g'}z=0$. Since $q-\ip gz\le0\le q'$ and
$q'-\ip{g'}{z'}\le0\le q$, \eqref{eq:tail-planes} holds for $(z,g,q)$
together with $(z',g',q')$. Let $1\le i<r_j$ and let
$\mathcal S_{j,i}$ and $\mathcal S_{j,i+1}$ be images of
$\mathcal S_{(k',b]}$ and $\mathcal S_{(b,b']}$ with
$(k',b],(b,b']\in\Gamma_\ell$. By
\cref{sc:lem:ordinary}\ref{it:ord-gain}, $\mathcal S_{(b,b']}$
satisfies \eqref{sc:eq:start-projections} with constant
$C_{\mathrm{ord}}$. A rotation and \eqref{eq:scaling} preserve
\eqref{sc:eq:start-projections} with constant $C_{\mathrm{ord}}$.
Hence \cref{sc:lem:ordinary}\ref{it:ord-join} applied with
$C_0=C_{\mathrm{ord}}$ gives \eqref{eq:tail-planes} for every triple of
$\mathcal S_{j,i}$ together with every triple of $\mathcal S_{j,i+1}$.

Let $a'\ne\omega$ be a step of depth $N$ in $I$ and let $a$ be the
smallest step of depth $N$ in $I$ with $a>a'$. Then $k_a=j_{a'}$. Let
$a\ne\omega$. Since $h_a=(1+o(1))h_{a'}$, we have
$(\rho^2+1)h_a\le(\rho^2+2)h_{a'}$ and $h_{a'}/2\le h_a\le 2h_{a'}$ for
sufficiently large $t$. Let $x'$ be the first point of
$\mathcal S_a$. By \eqref{sc:eq:four-boundary} for $a$ we have
$c_{\mathrm{four}}h_a\le\norm{x'}\le C_{\mathrm{four}}h_a$ and every
triple $(z,g,q)$ of $\mathcal S_a$ satisfies
$\norm g\le C_{\mathrm{four}}$, $\ip{x'}g\ge0$ and
$0\le\ip{x'}z\le(\rho^2+1)h_aq+P^2$. Hence $\mathcal S_a$ satisfies
\eqref{sc:eq:right-piece} for $a'$ with constants
$c_{\mathrm{four}}/2$, $2C_{\mathrm{four}}$.

Let $(k_\omega,b]$ be the set of $\Gamma_\ell$ with smallest indices.
Then $b=j$ for $\ell=N$, $b=\omega$ for $\ell=N+1$ and $b=\omega L$
for $\ell=N+2$. Let $a'$ be the step of depth $N$ in $I$ with
$j_{a'}=k_\omega$. Write $\mathcal S'=\mathcal S_{(k_\omega,b]}$, let
$x'$ be the first point of $\mathcal S'$ and let
$J=\{k_\omega+1,\ldots,b-1\}$. Then $J=I_\omega$ for $\ell=N$,
$J=I_{\omega L}$ for $\ell=N+1$ and $J=I_{\omega LL}$ for $\ell=N+2$.
Since $h_\omega=(1+o(1))h_{a'}$, we have
$(\rho^2+1)h_\omega\le(\rho^2+2)h_{a'}$ and $h_\omega\ge h_{a'}/2$ for
sufficiently large $t$. Since
$m_J\le m_\omega=(\rho^2+O(y^{-1/4}))h_\omega$ by
\cref{sc:lem:balance}\ref{it:bal-prop} and
$m_J/\kappa+C/\log A+P/h_\omega=o(1)$, \eqref{sc:eq:ord-norm} and
\eqref{sc:eq:ord-boundary} give
$\norm{x'}\le(\rho^2+1)h_\omega\le(\rho^2+2)h_{a'}$, $\norm g\le C$,
$\ip{x'}g\ge0$ and $0\le\ip{x'}z\le(\rho^2+2)h_{a'}q+P^2$ for every
triple $(z,g,q)$ of $\mathcal S'$. Let $w=\omega$ for $\ell=N$,
$w=\omega L$ for $\ell=N+1$ and $w=\omega LL$ for $\ell=N+2$. Then
$w\in J\cap B$, $\{k_\omega+1,\ldots,w-1\}=I_{wL}$ and
$h_w\ge ch_\omega$ by \cref{sc:lem:balance}\ref{it:bal-parent}. By
\eqref{sc:eq:ord-lower} and \eqref{sc:eq:balanced-debt},
\[
 \norm{x'}\ge\frac{h_w}2-C\sum_{u\in I_{wL}}h_u^-
 \ge\frac{h_w}2-Cy^{-1/4}h_w\ge ch_\omega\ge(c/2)h_{a'}.
\]
Thus $\mathcal S'$ satisfies \eqref{sc:eq:right-piece} for $a'$ with
absolute constants $c_0,C_0$.

Let $1\le i<r_j$ and let $\mathcal S_{j,i}$ be the image of
$\mathcal S_{a'}$ with $a'\ne\omega$. Let $a$ be the smallest step of depth $N$ in $I$ with $a>a'$. Then
$\mathcal S_{j,i+1}$ is the image of $\mathcal S_a$ for $a\ne\omega$
and the image of $\mathcal S_{(k_\omega,b]}$ for $a=\omega$. Hence
\cref{sc:lem:four-phase}\ref{it:four-join} for $a'$ gives
\eqref{eq:tail-planes} for every triple of $\mathcal S_{j,i}$ together
with every triple of $\mathcal S_{j,i+1}$. Thus \eqref{eq:tail-planes}
holds for every two triples of
$\mathcal S_{j,1},\ldots,\mathcal S_{j,r_j}$. For sufficiently large
$t$ we have $N\ge1$. Hence $\mathcal S_{j,1}$ is the image of
$\mathcal S_a$ for the smallest step $a$ of depth $N$ in $I$ and
$a\ne\omega$. By \cref{sc:lem:four-phase}\ref{it:four-boundary},
$\mathcal S_{j,1}$ satisfies \eqref{sc:eq:start-projections} with
constant $C_{\mathrm{four}}$. Also $\mathcal S_{j,r_j}$ is the
image under a rotation and \eqref{eq:scaling} of
$\mathcal S_{(k',j]}$ with $(k',j]\in\Gamma_\ell$. Thus
\cref{sc:prop:working-realization}\ref{it:work-S} holds with
$C_0=C_{\mathrm{four}}$ for $\ell=N,N+1,N+2$.

Let $j\in S$ and $\ell\in\{N,N+1,N+2\}$. The sets $(k_a,j_a]$ with
$a$ of depth $N$ in $I$ form a partition of $(k,j]$. Hence
\begin{equation}\label{eq:S-term-split}
 \log\frac{(m_I+P)\norm{x_j^{(\ell)}}}{h_j\norm{x_k^{(\ell)}}}
 \ge\log\frac{m_I}{h_j}
 +\sum_{\text{depth of }a=N}
 \log\frac{\norm{x_{j_a}^{(\ell)}}}{\norm{x_{k_a}^{(\ell)}}} .
\end{equation}
Let $a\ne\omega$ be a step of depth $N$ in $I$. Since
$P/h_a\le y^{-35/12}\le Ct^{-2}$, \eqref{sc:eq:four-gain} for $a$
gives
\begin{multline*}
 \log\frac{\norm{x_{j_a}^{(\ell)}}}{\norm{x_{k_a}^{(\ell)}}}
 =(\psi_N+O(t^{-1/4}))\log Q_{a,0}(\mu_{a,7})
  +\psi_{N+1}\log Q_{a,1}(\mu_{a,7})\\
  +(\psi_{N+2}+O(t^{-1/4}))\log Q_{a,2}(\mu_{a,7})+O(t^{-2}).
\end{multline*}
By \eqref{sc:eq:last-gain},
$\sum_\ell\psi_\ell\log(\norm{x_j^{(\ell)}}/\norm{x_{k_\omega}^{(\ell)}})$
is at least
$\psi_N\log Q_{\omega,0}(h_j)+\psi_{N+1}\log Q_{\omega,1}(h_j)+\psi_{N+2}\log Q_{\omega,2}(h_j)-C/t$.

Let $a\ne\omega$ be a step of depth $N$ in $I$ and write $b=j_a$. Then
\eqref{sc:eq:four-masses} for $a$ and $i=3$ gives
$\mu_{a,7}=m_b-m_{bL}-m_{bR}$. By \eqref{sc:eq:high-gain-score},
$Q_{a,0}(\mu_{a,7})=(m_b-m_{bL}-m_{bR})/m_a$ and
$Q_{\omega,0}(h_j)=h_j/m_\omega$. Every step of $I$ of depth less than
$N$ equals $j_a$ for exactly one step $a\ne\omega$ of depth $N$ in $I$.
Every step of $I$ of depth between $1$ and $N$ is a child of one step
of depth less than $N$. Hence, with products over steps of $I$,
\[
 m_I\prod_{\text{depth of }a=N}\frac1{m_a}
 \prod_{\text{depth of }b<N}(m_b-m_{bL}-m_{bR})
 =\prod_{\text{depth of }b<N}
 \frac{(m_b-m_{bL}-m_{bR})m_b}{m_{bL}m_{bR}} .
\]
Therefore
\[
 \log\frac{m_I}{h_j}+\log Q_{\omega,0}(h_j)
 +\sum_{\substack{\text{depth of }a=N\\ a\ne\omega}}\log Q_{a,0}(\mu_{a,7})
 =\sum_{\text{depth of }b<N}
 \log\frac{(m_b-m_{bL}-m_{bR})m_b}{m_{bL}m_{bR}} .
\]
For a step $a$ of depth $N$ in $I$, $\tau>0$ and numbers
$\gamma_0,\gamma_1,\gamma_2$ with $\gamma_0+\gamma_1+\gamma_2=1$,
\[
 \gamma_0\log Q_{a,0}(\tau)+\gamma_1\log Q_{a,1}(\tau)+\gamma_2\log Q_{a,2}(\tau)
 =\log Q_{a,0}(\tau)+(\gamma_1+\gamma_2)\log\frac{Q_{a,1}(\tau)}{Q_{a,0}(\tau)}
  +\gamma_2\log\frac{Q_{a,2}(\tau)}{Q_{a,1}(\tau)} .
\]
There are at most $2^N$ steps of depth $N$ in $I$ and
$2^Nt^{-2}\le C/t$. By \eqref{eq:S-term-split} and
\eqref{sc:eq:four-ratios} for every step $a$ of depth $N$ in $I$,
\begin{equation}\label{eq:S-term-quotients}
 \sum_\ell\psi_\ell\log\frac{(m_I+P)\norm{x_j^{(\ell)}}}{h_j\norm{x_k^{(\ell)}}}
 \ge\sum_{\substack{b\in I\\ \text{depth of }b<N+2}}
 \Bigl(\sum_{\ell>\text{depth of }b}\psi_\ell+O(t^{-1/4})\Bigr)
 \log\frac{(m_b-m_{bL}-m_{bR})m_b}{m_{bL}m_{bR}}-\frac Ct .
\end{equation}
The term $O(t^{-1/4})$ in \eqref{eq:S-term-quotients} is zero for $b$
of depth less than $N$. Every $b\in B$ with $t/3<\log h_b\le t$
satisfies
\begin{equation}\label{sc:eq:remaining-log}
 \log\frac{(m_b-m_{bL}-m_{bR})m_b}{m_{bL}m_{bR}}
 \ge\frac{\Delta_b-Ct^{-1/4}|\Delta_b|}{h_b^{1/p}}.
\end{equation}
By \eqref{sc:eq:balanced-unit-K} applied with $y=\log h_b$,
$m_b-m_{bL}-m_{bR}=(1+O(t^{-1/4}))h_b>0$. By
\cref{sc:lem:balance}\ref{it:bal-prop} and \eqref{sc:eq:split-size},
$0<m_{bL}\le m_b$ and $m_{bR}>0$. Hence
\cref{lem:max-plus-excess}\ref{it:split-real} applied with
$(h_v,P,m_v,m_{vL},m_{vR})=(m_b-m_{bL}-m_{bR},0,m_b,m_{bL},m_{bR})$
gives
$\log[(m_b-m_{bL}-m_{bR})m_b/(m_{bL}m_{bR})]\ge\Delta_b/(m_b-m_{bL}-m_{bR})^{1/p}$.
Thus \eqref{sc:eq:remaining-log} holds. Every $b$ in
the sum of \eqref{eq:S-term-quotients} lies in $B$ and has
$t/3<\log h_b\le t$. The coefficients of the logarithms in
\eqref{eq:S-term-quotients} are nonnegative. By
\eqref{eq:S-term-quotients} and \eqref{sc:eq:remaining-log},
\begin{equation}\label{sc:eq:prescribed-term}
 \sum_\ell\psi_\ell\log\frac{(m_I+P)\norm{x_j^{(\ell)}}}{h_j\norm{x_k^{(\ell)}}}
 \ge\sum_\ell\psi_\ell
 \sum_{\substack{b\in I\\ \text{depth of }b<\ell}}
 \frac{\Delta_b}{h_b^{1/p}}
 -\frac Ct
 -Ct^{-1/4}\sum_{\substack{b\in I\\ \text{depth of }b<N+2}}
 \frac{\lvert\Delta_b\rvert}{h_b^{1/p}} .
\end{equation}

Let $b$ and the parent of $b$ lie in $B$ with $\log h_b>t/3$. By
\cref{sc:lem:balance}\ref{it:bal-parent} the stepsize of the parent of
$b$ is $\rho e^{O(t^{-1/4})}h_b$. Also $h_v\le e^y$ and the
parent of $I$ has stepsize above $e^y$. Hence a step $b$ of $I$ that has depth less than $\ell$ and
$h_b\le e^y/\rho^\ell$ or has depth at least $\ell$ and
$h_b>e^y/\rho^\ell$ satisfies
\begin{equation}\label{eq:depth-exception}
 \lvert y-\log h_b-\ell\log\rho\rvert\le C(\ell+1)t^{-1/4}.
\end{equation}
If $b$ has depth at least $\ell$ and $h_b>e^y/\rho^\ell$, then
$b$ lies in $I_w$ for a step $w$ of depth $\ell$ in $I$ and
$e^y/\rho^\ell<h_b\le h_w\le e^y\rho^{-\ell}e^{C\ell t^{-1/4}}$.
If $b$ has depth $\nu<\ell$ and $h_b\le e^y/\rho^\ell$,
then $e^y/\rho^\ell\ge h_b\ge e^y\rho^{-(\nu+1)}e^{-C(\nu+1)t^{-1/4}}$.
Hence $\nu=\ell-1$ and \eqref{eq:depth-exception} holds. Every step $b$
of $I$ of depth less than $N+2$ satisfies
$y-\log h_b\le(N+2)\log\rho+1$ for sufficiently large $t$. For fixed
$b$,
\[
 \int_{\{y:\eqref{eq:depth-exception}\}}e^{(y-\log h_b)/p}\,dy
 \le C2^N(N+1)t^{-1/4},\qquad
 \int_{\log h_b}^{\log h_b+(N+2)\log\rho+1}e^{(y-\log h_b)/p}\,dy
 \le C2^N .
\]
By \eqref{sc:eq:working-excess} for the range $(t/3,t]$,
\begin{equation}\label{sc:eq:depth-replace-int}
\begin{aligned}
 &\int_s^t e^{y/p}\sum_{j\in S}\sum_\ell\psi_\ell
 \Biggl|\sum_{\substack{b\in I\\ \text{depth of }b<\ell}}\frac{\Delta_b}{h_b^{1/p}}
 -\sum_{\substack{b\in I\\ h_b>e^y/\rho^\ell}}\frac{\Delta_b}{h_b^{1/p}}\Biggr|\,dy\\
 &\qquad+\int_s^t e^{y/p}\sum_{j\in S}Ct^{-1/4}
 \sum_{\substack{b\in I\\ \text{depth of }b<N+2}}\frac{\lvert\Delta_b\rvert}{h_b^{1/p}}\,dy
 \le C2^N(N+1)t^{-1/4}\sum_{t/3<\log h_b\le t}\lvert\Delta_b\rvert\le Cn .
\end{aligned}
\end{equation}
By \eqref{sc:eq:count-bound} there are at most $ne^{-y/p}$ cuts at
threshold $e^y$ and hence $\int_s^t e^{y/p}(C/t)\#S\,dy\le Cn$.

Let $j\notin S$ be a cut at threshold $e^y$. The stepsizes of the steps of
$V\cap I\cap B$ lie in $(e^y/\rho^\ell,e^y]$ and decrease by a factor
$\rho e^{O(t^{-1/4})}$ from a parent in $B$ to a child in $B$. Hence a
subset of $V\cap I\cap B$ that is connected in the binary tree of $I$
spans at most $\ell+O(1)$ depth levels and has at most $C2^N$ elements.
Because $V\cap I$ is closed under taking parents in $I$, each connected
component of $V\cap I\cap B$ is rooted at the root of $I$ or at a child
of a step in $V\cap I\setminus B$. Hence $V\cap I\cap B$ has at most
$1+2\#(V\cap I\setminus B)$ connected components. Thus
\[
 1+\#(V\cap I)\le C2^N\bigl(1+\#(V\cap I\setminus B)\bigr).
\]
Summing over the cuts $j\notin S$ gives
\[
 \frac Cy\bigl(\#\{j\notin S:h_j>e^y\}+\#V\bigr)
 \le\frac{C2^N}y\bigl(\#\{j\notin S:h_j>e^y\}+\#(V\setminus B)\bigr).
\]
By \cref{sc:prop:working-realization} and
$(e^y/h_n)^{1/(2p)}\le y^{-3/(2p)}\le1/y$, the number of cuts without
\eqref{eq:three-cuts} is at most $Cne^{-y/p}/y$. Hence the integral
over $s\le y\le t$ of $e^{y/p}C2^N/y$ times the number of cuts without
\eqref{eq:three-cuts} is at most $C2^Nn/t\le Cn$.

To each $b\in V\setminus B$ and to each cut $j\notin S$ with
\eqref{eq:three-cuts} we assign a step $r\in B'$. For
$b\in V\setminus B$ let $r=b$. Then $h_r>e^y/\rho^{N+2}$. For a cut
$j\notin S$ with \eqref{eq:three-cuts} such that the parent of $I$ is
not in $B$, let $r$ be the parent of $I$. Then $h_r>e^y$. Let
$j\notin S$ be a cut with \eqref{eq:three-cuts} such that the parent of
$I$ lies in $B$. Then some positive step of $I$ of depth at most $N+2$
lies outside $B$. Among the positive steps of $I$ outside $B$, take the
step $u$ of smallest depth and, among steps of equal depth, of least
index. Let $r$ be the parent of $u$. Then $r\in B$. Either
$r$ is the parent of $I$ with $h_r>e^y$ or $r\in I$ has depth at most
$N+1$ and every step on the path from the parent of $I$ to $r$ in the
binary tree lies in $B$. Let $r\in I$. By \cref{sc:lem:balance}\ref{it:bal-parent} the quotient
of the stepsizes of a parent in $B$ and a child in $B$ is at most an
absolute constant. Since $N=O(\log\log t)$ and $y\ge t/2$, every step
on the path from the parent of $I$ to $r$ has $\log h>t/3$. Hence
\cref{sc:lem:balance}\ref{it:bal-parent} along the path gives
$h_r\ge ce^y\rho^{-N-2}$. Thus the step $r\in B'$ assigned to
$b\in V\setminus B$ or to a cut $j\notin S$ with \eqref{eq:three-cuts}
satisfies
\begin{equation*}
 e^y\le C\rho^Nh_r .
\end{equation*}
For fixed $y$, a step $r\in B'$ lies in at most one interval $I$ and is
the parent of at most two intervals $I$. Hence one step $r\in B'$ is
assigned to at most one $b\in V\setminus B$ and to at most three cuts
$j\notin S$ with \eqref{eq:three-cuts}. For fixed $r\in B'$,
\[
 \int_{\{s\le y\le t:\,e^y\le C\rho^Nh_r\}}e^{y/p}\,\frac{C2^N}y\,dy
 \le\frac{C4^N}t\min\{h_r,e^t\}^{1/p}.
\]
By \eqref{sc:eq:unbalanced-count},
\[
  \sum_{r\in B'}\frac{C4^N}t\min\{h_r,e^t\}^{1/p}
 \le\frac{C4^Nn}{\sqrt t}\le Cn .
\]
Together with the bound $C2^Nn/t\le Cn$ for the cuts without
\eqref{eq:three-cuts} we obtain
\begin{equation}\label{sc:eq:error-integrated}
 \int_s^t e^{y/p}\,\frac Cy\bigl(\#\{j\notin S:h_j>e^y\}+\#V_\ell(y)\bigr)\,dy
 \le Cn
\end{equation}
for each $\ell\in\{N,N+1,N+2\}$.

For each fixed $s\le y\le t$, multiply \eqref{sc:eq:working-quantitative}
by $\psi_\ell$ and sum over $\ell=N,N+1,N+2$. Substituting
\eqref{sc:eq:prescribed-term} for $j\in S$ and \eqref{eq:V-sum} for
$b\in V_\ell(y)$, absorbing the terms $C/y$ of \eqref{eq:V-sum} into
$(C/y)(\#\{j\notin S:h_j>e^y\}+\#V_\ell(y))$ and using
$\sum_\ell\psi_\ell=1$ and $(e^y/h_n)^{1/(2p)}\le y^{-3/(2p)}\le1/y$
eliminates $\norm{x_j^{(\ell)}}/\norm{x_k^{(\ell)}}$ and gives
\begin{multline}\label{sc:eq:pointwise-weighted}
 \sum_\ell\psi_\ell\biggl(
 \sum_{j\in S}\sum_{\substack{b\in I\\ \text{depth of }b<\ell}}
 \frac{\Delta_b}{h_b^{1/p}}
 +\sum_{b\in V_\ell(y)}\frac{\Delta_b}{h_b^{1/p}}
 -\frac Cy\bigl(\#\{j\notin S:h_j>e^y\}+\#V_\ell(y)\bigr)\biggr)\\
 -\frac Ct\#S
 -Ct^{-1/4}\sum_{j\in S}\sum_{\substack{b\in I\\ \text{depth of }b<N+2}}
 \frac{\lvert\Delta_b\rvert}{h_b^{1/p}}
 \le\sum_{h_j>e^y}\log\frac{m_I+P(y)}{h_j}+\frac{Cn}{ye^{y/p}} .
\end{multline}
Both sides of \eqref{sc:eq:pointwise-weighted} are Borel functions of
$y$. Every positive step $b$ with
$e^y/\rho^\ell<h_b\le e^y$ lies in the interval $I$ of one cut $j$. For
$j\notin S$ the positive steps $b$ of $I$ with $h_b>e^y/\rho^\ell$ are
the elements of $V_\ell(y)\cap I$. Multiply
\eqref{sc:eq:pointwise-weighted} by $e^{y/p}$ and integrate over
$s\le y\le t$. By \eqref{sc:eq:depth-replace-int},
$\int_s^t e^{y/p}(C/t)\#S\,dy\le Cn$, \eqref{sc:eq:error-integrated},
\eqref{sc:eq:first-log-integral} and $\int_s^t dy/y\le\log4$ we obtain
\eqref{eq:weighted-delta-score}. Thus \eqref{eq:weighted-score} holds.
\end{proof}

\subsection{The number of large steps}\label{sc:sec:proof-count}

\begin{proof}[Proof of \cref{sc:thm:count}]
Let $t_0>0$ be the number with $e^{t_0}t_0^3=h_n$. Then
$\log h_n=t_0+3\log t_0$. Along the tail of the sequence where
$t_0>1$, we have $\log t_0>0$, $0<t_0-\log t_0<t_0<\log h_n$ and
$(t_0-\theta)/\log h_n\to1$ uniformly for $0\le\theta\le\log t_0$.
Because $\mathcal P$ is a step function taking finitely many values on
$[t_0-\log t_0,t_0]$, for each schedule with $t_0>1$ choose
$\theta_0\in[0,\log t_0]$ attaining
$\inf_{0\le\theta\le\log t_0}\mathcal P(t_0-\theta)$. Applying
\cref{sc:lem:initial-power} along $s=t_0-\theta_0$ gives
$\mathcal P(t_0-\theta_0)\ge(p-o(1))n\ge n$ along the tail of the
sequence and hence $\mathcal P(t_0-\theta)\ge n$ for every
$0\le\theta\le\log t_0$.

Let $t$ be sufficiently large with $e^t t^3\le h_n$ and let
$t/2\le s\le2t/3$. By \eqref{sc:eq:band-sums} two values of
$\mathcal P$ on $[s-\log(2s),s]$ differ by at most
\begin{equation}\label{eq:window-sum}
 \sum_{s-\log(2s)<\log h_b\le s}\lvert\Delta_b\rvert .
\end{equation}
Each positive step $b$ satisfies $s-\log(2s)<\log h_b\le s$ for a set
of $s\in[t/2,2t/3]$ of length at most $\log(2t)$. By
\eqref{sc:eq:working-excess} for the range $(t/3,t]$, the mean of
\eqref{eq:window-sum} over $t/2\le s\le2t/3$ is therefore at most
$Cn(\log t)/t$. Choose $s=s(t)\in[t/2,2t/3]$ with
\eqref{eq:window-sum} at most $Cn(\log t)/t$.

Let $N$ be the smallest integer with $2^N\ge(\log t)^2$ and let
$\psi_N,\psi_{N+1},\psi_{N+2}$ be the weights of
\cref{sc:lem:two-ranges}. By \cref{sc:lem:two-ranges},
\[
 \sum_{\ell=N}^{N+2}\psi_\ell\int_0^{\ell\log\rho}e^{\theta/p}
 [\mathcal P(t-\theta)-\mathcal P(s-\theta)]\,d\theta\le Cn .
\]
For every $\ell\in\{N,N+1,N+2\}$ we have $\ell\log\rho\le\log t\le\log(2s)$
and
\[
 \sum_{\ell=N}^{N+2}\psi_\ell\int_0^{\ell\log\rho}e^{\theta/p}\,d\theta
 =p\sum_{\ell=N}^{N+2}\psi_\ell(2^\ell-1)\ge c(\log t)^2 .
\]
For every $\ell\in\{N,N+1,N+2\}$ and $0\le\theta\le\ell\log\rho$, we
have $s-\theta\in[s-\log(2s),s]$ and
$[s-\log s,s]\subset[s-\log(2s),s]$. Hence
\[
 \mathcal P(t-\theta)-\mathcal P(s-\theta)
  \ge\inf_{0\le\theta'\le\log t}\mathcal P(t-\theta')
   -\inf_{0\le\theta'\le\log s}\mathcal P(s-\theta')-Cn(\log t)/t .
\]
Substituting the lower bound for $\mathcal P(t-\theta)-\mathcal P(s-\theta)$
into \eqref{eq:weighted-score}, dividing by
$p\sum_{\ell=N}^{N+2}\psi_\ell(2^\ell-1)\ge c(\log t)^2$ and using
$(\log t)/t\le(\log t)^{-2}$ gives
\begin{equation}\label{sc:eq:descent}
 \inf_{0\le\theta\le\log s}\mathcal P(s-\theta)
 \ge \inf_{0\le\theta\le\log t}\mathcal P(t-\theta)-Cn/(\log t)^2 .
\end{equation}

Fix a sufficiently large absolute constant $t_\ast$. Consider the tail
of the sequence where $t_0\ge t_\ast$ and
$\mathcal P(t_0-\theta)\ge n$ for every $0\le\theta\le\log t_0$. For
$i\ge0$ with $t_i\ge t_\ast$ define $t_{i+1}=s(t_i)$. Then
$t_i/2\le t_{i+1}\le2t_i/3$ and \eqref{sc:eq:descent} with $t=t_i$ and
$s=t_{i+1}$ gives
\[
 \inf_{0\le\theta\le\log t_{i+1}}\mathcal P(t_{i+1}-\theta)
 \ge\inf_{0\le\theta\le\log t_i}\mathcal P(t_i-\theta)-Cn/(\log t_i)^2 .
\]
Let $r$ be the smallest index with $t_r<t_\ast$. Then $r\ge1$,
$t_{r-1}\ge t_\ast$ and $t_r\ge t_{r-1}/2\ge t_\ast/2>1$. Since
$t_{i+1}\le2t_i/3$, we have $\log t_i-\log t_{i+1}\ge\log(3/2)$. A comparison with
$\int_{\log t_\ast}^\infty\tau^{-2}\,d\tau$ gives
$\sum_{i=0}^{r-1}Cn/(\log t_i)^2\le Cn/\log t_\ast\le n/2$ for
sufficiently large $t_\ast$. Thus $t_r\in(1,t_\ast)$ and
$\mathcal P(t_r)\ge\inf_{0\le\theta\le\log t_r}\mathcal P(t_r-\theta)\ge n/2$.
By \eqref{sc:eq:power-small-threshold},
\[
 n/2\le\mathcal P(t_r)\le Ce^{t_\ast/p}\#\{i:h_i>0\},
\]
contrary to \eqref{sc:eq:standing-sparse}.
\end{proof}

\section*{AI Use Statement}
\addcontentsline{toc}{section}{AI Use Statement}
AI tools were used to assist in writing and checking the mathematical
statements and proofs.

\begin{thebibliography}{99}\small
\setlength{\itemsep}{1pt plus .5pt minus .5pt}
\setlength{\parsep}{0pt}
\Urlmuskip=0mu plus 1mu\relax
\expandafter\def\expandafter\UrlBreaks\expandafter{\UrlBreaks\do\-}
\addcontentsline{toc}{section}{References}
\bibitem{AGZ21} N.~Agarwal, S.~Goel, and C.~Zhang.
Acceleration via fractal learning rate schedules.
\emph{Proceedings of the 38th International Conference on Machine
Learning}, Proceedings of Machine Learning Research 139, 87--99 (2021).

\bibitem{Alt18} J.~M.~Altschuler.
Greed, hedging, and acceleration in convex optimization.
Master's thesis, Massachusetts Institute of Technology (2018).

\bibitem{AP24} J.~M.~Altschuler and P.~A.~Parrilo.
Acceleration by random stepsizes: Hedging, equalization, and the arcsine stepsize schedule.
\href{https://arxiv.org/abs/2412.05790v2}{arXiv:2412.05790v2} (2026).

\bibitem{AP-I} J.~M.~Altschuler and P.~A.~Parrilo.
Acceleration by stepsize hedging: Multi-step descent and the silver stepsize schedule.
\emph{Journal of the ACM}, 72(2), Article~12 (2025).

\bibitem{altschuler2025silver} J.~M.~Altschuler and P.~A.~Parrilo.
Acceleration by stepsize hedging: Silver stepsize schedule for smooth convex optimization.
\emph{Mathematical Programming}, 213(1--2), 1105--1118 (2025).

\bibitem{AP26} J.~M.~Altschuler and P.~A.~Parrilo.
Stepsize hedging: An alternative mechanism for accelerating gradient descent.
\emph{INFORMS Computing Society Newsletter}, 19--27 (2026).

\bibitem{bell2026predetermined} L.~Bell.
Predetermined gradient descent: Horizon, dimension, and step signs.
Manuscript (2026).\newline
\url{https://www.primalspace.me/papers/gradient-descent-dimension-gaps/horizon-dimension-bounds.pdf}
(accessed 18 November 2026).

\bibitem{BA25} J.~Bok and J.~M.~Altschuler.
Accelerating proximal gradient descent via silver stepsizes.
\emph{Proceedings of the 38th Conference on Learning Theory},
Proceedings of Machine Learning Research 291, 421--453 (2025).

\bibitem{DVR24} S.~Das~Gupta, B.~P.~G.~Van~Parys, and E.~K.~Ryu.
Branch-and-bound performance estimation programming: A unified methodology for constructing optimal optimization methods.
\emph{Mathematical Programming}, 204(1--2), 567--639 (2024).

\bibitem{Dro17} Y.~Drori.
The exact information-based complexity of smooth convex minimization.
\emph{Journal of Complexity}, 39, 1--16 (2017).

\bibitem{DT22} Y.~Drori and A.~Taylor.
On the oracle complexity of smooth strongly convex minimization.
\emph{Journal of Complexity}, 68, Article~101590 (2022).

\bibitem{DT14} Y.~Drori and M.~Teboulle.
Performance of first-order methods for smooth convex minimization: A novel approach.
\emph{Mathematical Programming}, 145(1--2), 451--482 (2014).

\bibitem{GSDTP22} B.~Goujaud, D.~Scieur, A.~Dieuleveut, A.~B.~Taylor,
and F.~Pedregosa.
Super-acceleration with cyclical step-sizes.
\emph{Proceedings of the 25th International Conference on Artificial
Intelligence and Statistics}, Proceedings of Machine Learning Research
151, 3028--3065 (2022).

\bibitem{Gri24} B.~Grimmer.
Provably faster gradient descent via long steps.
\emph{SIAM Journal on Optimization}, 34(3), 2588--2608 (2024).

\bibitem{GSW23} B.~Grimmer, K.~Shu, and A.~L.~Wang.
Accelerated gradient descent via long steps.
\href{https://arxiv.org/abs/2309.09961v2}{arXiv:2309.09961v2} (2023).

\bibitem{GSW25} B.~Grimmer, K.~Shu, and A.~L.~Wang.
Accelerated objective gap and gradient norm convergence for gradient descent via long steps.
\emph{INFORMS Journal on Optimization}, 7(2), 156--169 (2025).

\bibitem{GSW24} B.~Grimmer, K.~Shu, and A.~L.~Wang.
Composing optimized stepsize schedules for gradient descent.
\emph{Mathematics of Operations Research}, Articles in Advance (2025).

\bibitem{he2026optimal} H.~He and G.~Cheng.
Optimal polynomial exponents of gradient descent with arbitrary real stepsizes.
\href{https://ssrn.com/abstract=7569138}{SSRN:7569138} (2026).

\bibitem{JCY26} M.~Jung, H.~Cho, and C.~Yun.
Stronger lower bounds for (non-)anytime acceleration of gradient descent.
\href{https://arxiv.org/abs/2609.04032v2}{arXiv:2609.04032v2} (2026).

\bibitem{KMSSVY22} J.~Kelner, A.~Marsden, V.~Sharan, A.~Sidford,
G.~Valiant, and H.~Yuan.
Big-step-little-step: Efficient gradient methods for objectives with
multiple scales.
\emph{Proceedings of the 35th Conference on Learning Theory},
Proceedings of Machine Learning Research 178, 2431--2540 (2022).

\bibitem{KF16} D.~Kim and J.~A.~Fessler.
Optimized first-order methods for smooth convex minimization.
\emph{Mathematical Programming}, 159(1--2), 81--107 (2016).

\bibitem{KS24} G.~Kornowski and O.~Shamir.
Open problem: Anytime convergence rate of gradient descent.
\emph{Proceedings of the 37th Conference on Learning Theory},
Proceedings of Machine Learning Research 247, 5335--5339 (2024).

\bibitem{LF71} V.~I.~Lebedev and S.~A.~Finogenov.
Ordering of the iterative parameters in the cyclical Chebyshev
iterative method.
\emph{USSR Computational Mathematics and Mathematical Physics}, 11(2),
155--170 (1971).

\bibitem{LY26} K.~Liu and Y.~Ye.
Silver rate is (almost) optimal for gradient descent: The strongly convex case.
\href{https://arxiv.org/abs/2609.26873v2}{arXiv:2609.26873v2} (2026).

\bibitem{MC26} J.~Ma and Y.~Chen.
A lower bound for stepsize-based acceleration of gradient descent.
\href{https://arxiv.org/abs/2608.10418v1}{arXiv:2608.10418v1} (2026).

\bibitem{NY83} A.~S.~Nemirovsky and D.~B.~Yudin.
\emph{Problem Complexity and Method Efficiency in Optimization}.
John Wiley \& Sons (1983).

\bibitem{Nes83} Yu.~E.~Nesterov.
A method of solving a convex programming problem with convergence rate $O(1/k^2)$.
\emph{Doklady Akademii Nauk SSSR}, 269(3), 543--547 (1983). In Russian.

\bibitem{Nes18} Yu.~Nesterov.
\emph{Lectures on Convex Optimization}, 2nd edition.
Springer Optimization and Its Applications 137. Springer, Cham (2018).

\bibitem{Pol64} B.~T.~Polyak.
Some methods of speeding up the convergence of iteration methods.
\emph{USSR Computational Mathematics and Mathematical Physics}, 4(5),
1--17 (1964).

\bibitem{TD23} A.~Taylor and Y.~Drori.
An optimal gradient method for smooth strongly convex minimization.
\emph{Mathematical Programming}, 199(1--2), 557--594 (2023).

\bibitem{THG17} A.~B.~Taylor, J.~M.~Hendrickx, and F.~Glineur.
Smooth strongly convex interpolation and exact worst-case performance of first-order methods.
\emph{Mathematical Programming}, 161(1--2), 307--345 (2017).

\bibitem{Tsai26note} C.-E.~Tsai.
An improved lower bound for non-anytime gradient descent.
Online note (2026).\newline
\url{https://chungentsai.github.io/gd-lower-bounds.html}
(accessed 18 November 2026).

\bibitem{TFZH26} C.-E.~Tsai, I.~Fatkhullin, L.~Zhang, and N.~He.
Lower bounds for anytime acceleration of gradient descent.
\href{https://arxiv.org/abs/2607.02053v1}{arXiv:2607.02053v1} (2026).

\bibitem{ye2026beyond} Y.~Ye and K.~Liu.
Improved gradient descent lower bounds beyond Nesterov.
\href{https://arxiv.org/abs/2609.02855v2}{arXiv:2609.02855v2} (2026).

\bibitem{ye2026almost} Y.~Ye and K.~Liu.
Silver rate is (almost) optimal for gradient descent.
\href{https://arxiv.org/abs/2609.09152v3}{arXiv:2609.09152v3} (2026).

\bibitem{You53} D.~Young.
On Richardson's method for solving linear systems with positive definite matrices.
\emph{Journal of Mathematics and Physics}, 32(1--4), 243--255 (1953).

\bibitem{ZJ26sc} Z.~Zhang and R.~Jiang.
Accelerated concatenation of stepsize schedules for smooth strongly convex gradient descent.
\href{https://arxiv.org/abs/2611.05242v1}{arXiv:2611.05242v1} (2026).

\bibitem{ZJ26} Z.~Zhang and R.~Jiang.
Accelerated gradient descent by concatenation of stepsize schedules.
\emph{SIAM Journal on Optimization}, 36(2), 1182--1210 (2026).

\bibitem{ZLDC25} Z.~Zhang, J.~D.~Lee, S.~S.~Du, and Y.~Chen.
Anytime acceleration of gradient descent.
\emph{Proceedings of the 38th Conference on Learning Theory},
Proceedings of Machine Learning Research 291, 5991--6013 (2025).

\end{thebibliography}
\end{document}